\documentclass[reqno, 11pt]{amsart}
\usepackage[left=2cm, right=2cm, top=2cm, bottom=2cm]{geometry}
\usepackage[T1]{fontenc}
\usepackage[utf8]{inputenc}	
\usepackage[english]{babel}
\usepackage[babel]{csquotes}
\usepackage{amssymb}
\usepackage{amsmath}
\usepackage{amsthm}
\usepackage{latexsym}
\usepackage{xcolor}
\usepackage{mathrsfs}
\usepackage{bbm}
\usepackage{verbatim}
\usepackage[colorlinks=true, pdfstartview=FitV, linkcolor=blue, citecolor=blue, urlcolor=blue]{hyperref}
\usepackage[style=numeric, sorting=nty, sortcites=true, backend=bibtex, maxbibnames=10, giveninits=true]{biblatex}
\usepackage{faktor}
\usepackage[shortlabels]{enumitem}

\numberwithin{equation}{section}
\newtheorem{thm}{Theorem}[section]

\newtheorem{lem}[thm]{Lemma}
\newtheorem{prop}[thm]{Proposition}
\theoremstyle{definition}
\newtheorem{defin}[thm]{Definition}
\newtheorem{remark}[thm]{Remark}

\renewcommand{\d}{{\mathrm d}} 
\renewcommand{\div}{\operatorname{div}} 
\newcommand{\norm}[1]{\left\|#1\right\|} 
\newcommand{\ip}[2]{\left\langle #1,#2 \right\rangle} 
\newcommand{\ii}[4]{\int_{#1}^{#2} #3 \: \d#4} 
\def\laweq{\stackrel{\mathcal L}{=}}
\def\supp{\sup_{\tau \in [0,t]}}
\def\sups{\sup_{s \in [0,t]}}
\def\tom{\widetilde{\Omega}}
\def\tF{\widetilde{\mathscr{F}}}
\def\tP{\widetilde{\mathbb{P}}}
\def\tE{\widetilde{\mathbb{E}}}
\def\X{\mathcal{X}}
\def\hom{\widehat{\Omega}}
\def\hF{\widehat{\mathscr{F}}}
\def\hP{\widehat{\mathbb{P}}}
\def\hE{\widehat{\mathbb{E}}}
\def\hbu{\widehat{\bu}}
\def\hphi{\widehat\varphi}
\def\hmu{\widehat\mu}
\def\ia{\b{A}^{-1}}
\def\bu{{\boldsymbol u}}
\def\bv{{\boldsymbol v}}
\def\bn{{\boldsymbol n}}

\def\bH{{\boldsymbol H}}
\def\bHs{{\boldsymbol H_\sigma}}
\def\bVs{{\boldsymbol V_\sigma}}
\def\bVsd{{\boldsymbol{V}_\sigma^*}}
\def\b #1{{\boldsymbol #1}}
\def\enne{\mathbb{N}}

\def\erre{\mathbb{R}}

\def\P{\mathbb{P}}

\def\E{\mathop{{}\mathbb{E}}}

\def\cL{\mathscr{L}}
\def\cF{\mathscr{F}}

\def\cP{\mathscr{P}}

\def\OO{\mathcal{O}}
\def\embed{\hookrightarrow}

\def\andrea #1{{\color{black} #1}}
\def\mau #1{{\color{red} #1}}

\makeatletter
\DeclareFontFamily{OMX}{MnSymbolE}{}
\DeclareSymbolFont{MnLargeSymbols}{OMX}{MnSymbolE}{m}{n}
\SetSymbolFont{MnLargeSymbols}{bold}{OMX}{MnSymbolE}{b}{n}
\DeclareFontShape{OMX}{MnSymbolE}{m}{n}{
	<-6>  MnSymbolE5
	<6-7>  MnSymbolE6
	<7-8>  MnSymbolE7
	<8-9>  MnSymbolE8
	<9-10> MnSymbolE9
	<10-12> MnSymbolE10
	<12->   MnSymbolE12
}{}
\DeclareFontShape{OMX}{MnSymbolE}{b}{n}{
	<-6>  MnSymbolE-Bold5
	<6-7>  MnSymbolE-Bold6
	<7-8>  MnSymbolE-Bold7
	<8-9>  MnSymbolE-Bold8
	<9-10> MnSymbolE-Bold9
	<10-12> MnSymbolE-Bold10
	<12->   MnSymbolE-Bold12
}{}

\let\llangle\@undefined
\let\rrangle\@undefined
\DeclareMathDelimiter{\llangle}{\mathopen}%
{MnLargeSymbols}{'164}{MnLargeSymbols}{'164}
\DeclareMathDelimiter{\rrangle}{\mathclose}%
{MnLargeSymbols}{'171}{MnLargeSymbols}{'171}
\makeatother
\AtEveryBibitem{%
\clearfield{doi}%
\clearfield{url}%
\clearfield{issn}%
\clearfield{isbn}
} 
\renewbibmacro*{publisher+location+date}{%
\ifentrytype{book}
{\printlist{publisher}%
\iflistundef{location}
{\setunit*{\addcomma\space}}
{\setunit*{\addcomma\space}}%
\printlist{location}}
{\printlist{location}%
\iflistundef{publisher}
{\setunit*{\addcomma\space}}
{\setunit*{\addcomma\space}}%
\printlist{publisher}}
\setunit*{\addcomma\space}%
\usebibmacro{date}%
\newunit} 

\usepackage[]{amsaddr}
\begin{document}	
	\title[A stochastic Allen-Cahn-Navier-Stokes system
	with singular potential]{A stochastic Allen-Cahn-Navier-Stokes system\\
	with singular potential}
	
	\author{Andrea Di Primio, Maurizio Grasselli \and Luca Scarpa}
	\address{Dipartimento di Matematica,
		Politecnico di Milano, Via E.~Bonardi 9, 20133 Milano, Italy}
	\email{andrea.diprimio@polimi.it}
 	\email{maurizio.grasselli@polimi.it}
	\email{luca.scarpa@polimi.it}
	\subjclass[2020]{35Q35, 35R60, 60H15, 76T06}
	\keywords{Allen-Cahn-Navier-Stokes system; stochastic two-phase flow;
	logarithmic potential; martingale solutions; pathwise uniqueness.}
	
	\begin{abstract}
		We investigate a stochastic version of the Allen--Cahn--Navier--Stokes system in
		a smooth two- or three-dimensional domain with random initial data.
		The system consists of a Navier--Stokes equation coupled with a convective Allen--Cahn
		equation, with two independent sources of randomness given
		by general multiplicative-type Wiener noises.
		In particular, the Allen--Cahn equation is characterized by a singular potential of logarithmic type
		as prescribed by the classical thermodynamical derivation of the model.
		The problem is endowed with a no-slip boundary condition for the (volume averaged) velocity field,
		as well as a homogeneous Neumann condition for the order parameter.
		We first prove the existence of analytically weak martingale solutions in two and
		three spatial dimensions. Then, in two dimensions, we also estabilish
		pathwise uniqueness and the existence of a unique probabilistically-strong solution.
		Eventually, by exploiting a suitable generalisation of the
		classical De Rham theorem to stochastic processes, existence and uniqueness of a
		pressure is also shown.
	\end{abstract}
	\maketitle

	\section{Introduction} \label{sec:intro} \noindent
Modeling the behavior of immiscible (or partially miscible) binary fluids is a very active area of research because of its importance, for instance, in Biology and Materials Science. A well-known and effective approach is the so-called diffuse interface method (see, e.g., \cite{AMW}). This approach is based on the introduction of an order parameter
(or phase field) which accounts for the presence of the fluid components in a sufficiently smooth way, that is, there is no sharp interface separating them but a sufficiently thin region where there is some mixing. More precisely, denoting by $\varphi$ the relative difference between the (rescaled) concentrations of the two components, the regions $\{\varphi=1\}$ and $\{\varphi=-1\}$ represent the pure phases. However, they are separated by diffuse interfaces where $\varphi$ can take any intermediate value, i.e. $\varphi\in(-1,1)$. The interaction between the two components is a competition between the mixing entropy and demixing effects and can be represented by a potential energy density of the form
\begin{equation} \label{eq:fhpot}
		F(\varphi) = \dfrac{\theta}{2}\left[(1+\varphi)\log(1+\varphi)+(1-\varphi)\log(1-\varphi) \right]
		- \dfrac{\theta_0}{2}\varphi^2,
	\end{equation}
for some $0 < \theta < \theta_0$. This is known as the Flory--Huggins potential (see \cite{Flory42,Huggins41}). Letting
$\OO$ be a (sufficiently) smooth domain of $\mathbb{R}^d$, $d=2,3$, the Helmholtz free energy associated with the order parameter is then given by
$$
\mathcal{E}(\varphi) = \int_\OO \left(\frac{\varepsilon^2}{2} \vert \nabla \varphi \vert^2 + F(\varphi) \right)dy
$$
where the first term accounts for the surface energy separating the phases. Here $\varepsilon >0$ is related to the thickness of the diffuse interface.
Then, the functional derivative of $\mathcal{E}(\varphi)$ is called the chemical potential and usually denoted by $\mu$, that is,
$$
\mu = -\varepsilon^2 \Delta \varphi + F^\prime(\varphi).
$$
We can now introduce the two basic equations which govern the evolution of $\varphi$ in some time interval $(0,T)$: the Cahn-Hilliard equation
(see \cite{cahn-hill,CH1961})
$$
\partial_t \varphi = \Delta \mu
$$
and the Allen-Cahn equation (see \cite{AC1979})
$$
\partial_t \varphi = -\mu.
$$
Here we have assumed that the mobility is constant and equal to the unity. We also recall that, due to the singular behavior of the mixing entropy,
the Flory-Huggins potential \eqref{eq:fhpot} is often approximated with a regular potential like
	\begin{equation}
    \label{smoothpot}
	F(x) = \dfrac{1}{4}(x^2-1)^2, \qquad x \in \mathbb{R}.
	\end{equation}	
This choice simplifies the mathematical treatment. However, when the total mass of $\varphi$ is conserved (e.g. in \eqref{CH})
one cannot ensure that $\varphi$ takes its values in the physical range $[-1,1]$. Here we choose to keep the physically relevant potential
also in view of extending our analysis to conserved Allen-Cahn equations where in \eqref{AC} or in \eqref{eq:ACNS3} $\mu$ is replaced
by $\mu - \bar{\mu}$, $\bar{\mu}$ being the spatial average of $\mu$ (see \cite{RS1992}, see also \cite{GGW20} and references therein).

When we deal with a two-component fluid mixture, the equation for the phase variable is coupled with an equation for a suitably averaged velocity $\bu$ of the fluid mixture itself. A well known choice is the Navier--Stokes system subject to a capillary force, known as Korteweg force, which can be represented as $\mu\nabla\varphi$. More precisely, in the case of an incompressible mixture and taking $\varepsilon=1$,  constant density equal to the unity and constant viscosity $\nu>0$, we have
\begin{align}
\label{NS}
		&\partial_t \bu + (\bu\cdot\nabla)\bu -\nu \Delta\bu
		 + \nabla \pi = \mu\nabla\varphi\\
\label{inc}
		&\nabla \cdot \bu = 0
\end{align}
coupled with
\begin{equation}
\label{CH}
\partial_t\varphi + \bu\cdot\nabla\varphi = \Delta \mu
\end{equation}
or
\begin{equation}
\label{AC}
\partial_t\varphi + \bu\cdot\nabla\varphi = -\mu
\end{equation}		
in $(0,T) \times \OO$, for some given $T>0$. Here $\bu$ represents the volume averaged velocity and $\pi$ stands for the pressure. System \eqref{NS}-\eqref{inc} coupled with \eqref{CH} is known as Cahn--Hilliard--Navier--Stokes system, if \eqref{CH} is replaced by \eqref{AC} then the system is known as Allen--Cahn--Navier--Stokes system.
We recall that the standard boundary conditions are no-slip for $\bu$ and no-flux for \eqref{CH} or \eqref{AC}.

Starting from the pioneering contribution \cite{HH}, two-phase flow models have then been developed in several works. In particular, we refer to \cite{GPV} for the Cahn--Hilliard--Navier--Stokes system and to \cite{Blesgen} for the Allen--Cahn--Navier--Stokes system (see \cite{AGG2012,GKL2018,LT1998}
for more refined models with unmatched densities and \cite{HMR2012,HMR2012-2} for general thermodynamic derivations).
The corresponding mathematical analysis of such models has also experienced a remarkable development in the last decades.
Concerning the Cahn--Hilliard--Navier--Stokes system with matched densities see \cite{abels,GMT2019} and references therein
(see also \cite{ADG2013,AF2008,BBG2011,Giorgini2021,GT2020,GTV2021} for more general models).
Regarding the Allen--Cahn--Navier--Stokes system, we refer to \cite{GG2010,GG2010-2} for the matched case (see also references therein) and to \cite{FPP2019,FRPS2010,HM2019,JLL2017,GGW20,K2012,LDH2016,LH2018} for more refined models.

The deterministic description fails in rendering possible unpredictable oscillations at the microscopic level.
These include, for example, the environmental noise due to temperature and configurational effects.
The most natural way to take into account such factors was first proposed in \cite{cook} where a stochastic version of the Cahn--Hilliard equation
was introduced (see also \cite{BMW08, BGW10} for nucleation effects). That equation has been analyzed under various assssumptions in a number
of contributions (see, for instance, \cite{daprato-deb,deb-zamb, deb-goud, EM1991,goud} and, more recently, \cite{scar-SCH, scar-SVCH,scarpa21}).
We also refer to \cite{orr-roc-scar, scar-OVCSCH, scar-OCSHC} for related stochastic optimal control problems.
The stochastic Allen-Cahn equation has  been investigated in the framework of regular potentials
(see, for instance, \cite{BBP17, HRW12, HLR22, orr-scar} for examples of well-posedness analysis, see also \cite{BG19,RNT12} for numerical schemes and simulations).
The singular potential has been analyzed in \cite{Bertacco21}  (see also \cite{bauz-bon-leb} for the double obstacle potential), while in \cite{BOS21} the separation property from the pure phases has been established.

Here we analyze a stochastic version of the Allen--Cahn--Navier--Stokes system
characterized by two independent sources of randomness, the former acting on the fluid velocity and the latter acting on the order parameter dynamic
to incorporate thermal fluctuations.
More precisely, on account of \eqref{NS}-\eqref{inc} and \eqref{AC}, taking $\nu=1$ for the sake of simplicity,
we consider the following system of stochastic partial differential equations
	\begin{align}
		\label{eq:ACNS1}
		\d\bu +\left[-\Delta\bu
		+ (\bu\cdot\nabla)\bu + \nabla \pi - \mu\nabla\varphi\right]\,\d t
		= G_1(\bu)\,\d W_1
		\qquad&\text{in } (0,T)\times\OO,\\
		\label{eq:ACNS2}
		\nabla \cdot \bu = 0 \qquad & \text{in } (0,T) \times \OO,\\
		\label{eq:ACNS3}
		\d\varphi + \left[\bu\cdot\nabla\varphi + \mu\right]\,\d t =
		G_2(\varphi)\,\d W_2
		\qquad&\text{in } (0,T)\times\OO,\\
		\label{eq:ACNS4}
		\mu = -\Delta\varphi + F'(\varphi)
		\qquad&\text{in } (0,T)\times\OO,\\
		\label{eq:ACNS5}
		\bu=0, \quad \nabla\varphi\cdot\bn=0  \qquad&\text{on } (0,T)\times\partial\OO,\\
		\label{eq:ACNS6}
		\bu(0)=\bu_0, \quad \varphi(0)=\varphi_0 \qquad&\text{in } \OO.
	\end{align}
    Here  $W_1$ and $W_2$ are two independent cylindrical
	Wiener processes on some (possibly different) separable Hilbert spaces,
	and $G_i$ is a suitable stochastically integrable process with respect to $W_i$, for $i = 1,2$. Moreover,
    $\bn$ stands for the outward normal unit vector to $\partial\OO$.

The presence of random terms in both the equations has been considered in \cite{feir-petc2} in the case of a smooth potential like \eqref{smoothpot}
(see also \cite{GM21,TTM1,Medjo21} for modified models and \cite{DJMTT21,TTM2,TTM19,TTM22} for random terms only in the Navier-Stokes system).
We also remind that the case of Cahn--Hilliard--Navier--Stokes system for a compressible fluid has been studied in \cite{feir-petc} (see, e.g., \cite{TTM17,DNNTM21} for
random terms only in the Navier-Stokes system in the case of regular potential).

Here, for system \eqref{eq:ACNS1}--\eqref{eq:ACNS6} with a potential like \eqref{eq:fhpot}, we prove the existence of martingale solutions in dimension two and three,
as well as pathwise uniqueness and existence of probabilistically-strong solutions in dimension two.
The main difficulties on the mathematical side are two. The former is
the presence of noise also in the Allen-Cahn equation with singular potential:
this requires some ad hoc ideas based on a suitable compensation between the degeneracy of
the noise and the blow up of $F''$ at the endpoints (see \eqref{hyp:G2} below).
The latter is the coupling term $\mu\nabla\varphi$ in the Navier-Stokes equation.
Indeed, for the Allen--Cahn equation
one can recover only a $L^2(0,T; L^2(\OO))$-regularity for $\mu$, while
for the Cahn--Hilliard equation one gets $ \mu \in L^2(0,T;H^1(\OO))$.
This results in the necessity to reformulate the first equation for the fluid
in an alternative fashion, i.e., without employing $\mu$ explicitly.

We recall that, in \cite{feir-petc2}, the authors proved the existence of a (dissipative) martingale solution for a similar problem with a smooth potential.
Then, taking advantage of the smooth potential, they used the maximum principle to show that the range of the order parameter remains
confined in $[-1,1]$. Thus the global Lipschitz continuity of the potential and its derivatives holds. This fact was exploited to prove
the weak-strong pathwise and in law uniqueness in dimension three. However, if the potential is given by \eqref{eq:fhpot}, then no global Lipschitz continuity can be
achieved unless one can prove that the solution stays uniformly away from the pure phases, but this is not straightforward in the stochastic case (see \cite{BOS21} for the single
stochastic Allen--Cahn equation).

Besides the existence and uniqueness of solutions,
there are still a number of issues to investigate, which will be object of future work.
For example, regularity properties of the solution
and existence of analytically-strong solutions are open issues.
The low regularity of the chemical potential $\mu$ in the Allen-Cahn equation that we have mentioned above
seems to make the analysis challenging. Yet, some higher regularity properties have been shown in the deterministic case (see \cite{GGW20}). Their extension to the stochastic case is currently under investigation. Moreover, in the spirit of \cite{BOS21}, it would be interesting to
establish some random strict separation property from the pure phases. To do this, suitable regularity results might be needed.
It also worth pointing out that  system \eqref{eq:ACNS1}--\eqref{eq:ACNS6} is the
non-conserved version of the model, meaning that the spatial average of $\varphi$ is not preserved during the evolution.
The deterministic conserved version is now well-understood (see \cite{GGW20}). Its stochastic counterpart will also be the subject of further analysis. This issue will require a tuning of the diffusion coefficient $G_2$ (see for example \cite{ABDK16}). Finally, we point out that also more general versions of the stochastic Cahn--Hilliard-Navier--Stokes system might be considered on account of the recent advances in the analysis of the stochastic Cahn--Hilliard equation.

The content of this work is structured as follows. In Section \ref{sec:main}, we introduce the notation used throughout the work and state the main results. Sections \ref{sec:proof1} and \ref{sec:proof2} are devoted to the proof of existence of a martingale solutions and, in dimension two, of a probabilistically-strong solution, respectively.

	\section{Preliminaries and main results} \label{sec:main}
	\subsection{Functional setting and notation}
	For any (real) Banach space $E$, its (topological) dual is denoted by $E^*$
	and the duality pairing between $E^*$ and $E$ by $\ip{\cdot}{\cdot}_{E^*,E}$.
	If $E$ is a Hilbert space, then the scalar product of $E$ is denoted by $(\cdot,\cdot)_E$.
	For every couple of separable Hilbert spaces $E,F$ the space of Hilbert-Schmidt operators
	from $E$ to $F$ is denoted by the symbol $\cL^2(E,F)$
	and endowed with its canonical norm $\norm{\cdot}_{\cL^2(E,F)}$.
	Let $(\Omega,\cF,(\cF_t)_{t\in[0,T]},\P)$ be a filtered probability space
	satisfying the usual conditions (namely it is saturated and right-continuous),
	with $T>0$ being a prescribed final time.
 	We will use the symbol $\laweq$ to denote identity in law for
	random variables.
	Throughout the paper, $W_1$ and $W_2$ are independent
	cylindrical Wiener process on some separable Hilbert spaces $U_1$ and $U_2$,
	respectively. For convenience, we fix
	once and for all two complete orthonormal systems
	$\{u^1_j\}_{j\in\enne}$ on $U_1$ and $\{u^2_j\}_{j\in\enne}$ on $U_2$.
	We denote by $\cP$ the progressive sigma algebra on $\Omega\times[0,T]$.
	For every $s,r\in[1,+\infty]$ and for every Banach space $E$
	the symbols $L^s(\Omega; E)$ and $L^r(0,T; E)$
	indicate the usual spaces of strongly measurable Bochner-integrable functions
	on $\Omega$ and $(0,T)$, respectively.
	For all $s,r\in[1,+\infty)$ we write
	$L^s_\cP(\Omega;L^r(0,T; E))$ to stress
	that measurability is intended with respect to $\cP$.
	For all $s\in(1,+\infty)$ and for every separable and reflexive Banach space $E$
	we also define
	\[
	L^s_w(\Omega; L^\infty(0,T; E^*)):=
	\left\{v:\Omega\to L^\infty(0,T; E^*) \text{ weakly* measurable}\,:\,
	\norm{v}_{L^\infty(0,T; E^*)}\in L^s(\Omega)
	\right\}\,,
	\]
	which yields by
	\cite[Thm.~8.20.3]{edwards} the identification
	\[
	L^s_w(\Omega; L^\infty(0,T; E^*))=
	\left(L^{s/(s-1)}(\Omega; L^1(0,T; E))\right)^*\,.
	\]
	In case of distribution-valued processes, for every
	$s\in[1,+\infty)$, $r\in(0,+\infty)$, and $q\in(1,+\infty]$
	we set
	\begin{align*}
	L^s_\cP(\Omega; W^{-r,q}(0,T; E^*))
	&:=\left\{v:\Omega\to W^{-r,q}(0,T; E^*) \text{ weakly* measurable}\,:\right.\\
	&\left. \qquad \qquad \qquad v\in L^s(\Omega,\cF_t; W^{-r,q}(0,t; E))\quad\forall\,t\in[0,T]\right\}.
	\end{align*}
	Let $d = 2,3$ and consider a bounded domain $\OO\subset\erre^d$	with smooth boundary $\partial\OO$ and outward
	normal unit vector $\bn$. The spatiotemporal domains generated by $\OO$ are denoted by $Q:=(0,T)\times\OO$ and $Q_t:=(0,t)\times\OO$
	for all $t\in(0,T]$. Moreover, we employ the classical notation
	$W^{s,p}(\OO)$, where $s\in\erre$ and $p\in[1,+\infty]$, for the real Sobolev spaces
	and we denote by $\norm{\cdot}_{W^{s,p}(\OO)}$ their canonical norms.
	We define the Hilbert space $H^s(\OO):=W^{s,2}(\OO)$, $s\in\erre$,
	endowed with its canonical norm $\norm{\cdot}_{H^s(\OO)}$,
	and indicate by $H^1_0(\OO)$ the closure of $C^\infty_0(\andrea{\OO})$
	in $H^1(\OO)$.
	We now define the functional spaces
	\[
	H:=L^2(\OO)\,, \qquad V_1:=H^1(\OO)\,, \qquad V_2:=\left\{\psi\in H^2(\OO):\;\partial_\bn\psi=0\text{ a.e.~on } \partial\OO\right\}\,,
	\]
	endowed with their standard norms $\norm{\cdot}_H$,
	$\norm{\cdot}_{V_1}$, and $\norm{\cdot}_{V_2}$, respectively.
	As usual, we identify the Hilbert space $H$ with its dual through
	the Riesz isomorphism, so that we have the variational structure
	\[
	V_2\embed V_1\embed H \embed V_1^* \embed V_2^*\,,
	\]
	with dense and compact embeddings (both in the cases $d = 2$ and $d = 3$). We will also denote by $A:V_1\to V_1^*$
	the variational realization of the $-\Delta$ with homogeneous Neumann boundary
	condition, namely
	\[
	\ip{A\psi}{\phi}_{V_1^*,V_1}=\int_\OO\nabla\psi\cdot\nabla\phi\,,
	\qquad\psi,\phi\in V_1\,.
	\]
	For any Banach space $E$, we use the symbol $\boldsymbol E$ for the product space $E^d$. We also need to define the following solenoidal vector-valued spaces
	\andrea{
	\begin{align*}
		\bHs := \overline{\{ \bv \in  \b{C}^\infty_0(\OO): \nabla \cdot \bv = 0 \text{ in } \OO\}}^{\b{L}^2(\OO)}, \quad
		\bVs := \overline{\{ \bv \in \b{C}^\infty_0(\OO) : \nabla \cdot \bv = 0 \text{ in } \OO\}}^{\b{H}^1(\OO)}.
	\end{align*}}The space $\bHs$ is endowed with the norm $\norm{\cdot}_{\bHs}$
	of $\b H$ and its respective scalar product $(\cdot,\cdot)_{\bHs}$.
	\andrea{By means of the Poincaré inequality}, on the space $\bVs$
	we can use the norm $\norm{\bv}_{\bVs}:=\norm{\nabla\bv}_{\b L^2(\OO)}$, $\bv\in\bVs$,
	induced by the scalar product $(\cdot,\cdot)_{\bVs}$. The $d$-dimensional realisation of the \andrea{$-\Delta$ with homogeneous Dirichlet boundary condition
	$\b L: \b{H}^1_0(\OO) \to  (\b{H}_0^1(\OO))^*$ is defined as
	\[
	\ip{\b L\bv}{\b w}_{(\b{H}^1_0(\Omega))^*, \b{H}^1_0(\Omega)}:=
	(\nabla\bv, \nabla\b w)_{\b H}\,, \qquad \bv,\b w\in \b{H}^1_0(\Omega)\,.
	\]
	Furthermore, we also point out that for any $\b{u} \in (\b{H}^1_0(\Omega))^*$ and $\bv \in \left[ \mathcal{C}^\infty_0(\Omega) \right]^d$ we have
	\[
	\left\langle \bu, \bv \right\rangle_{(\b{H}^1_0(\Omega))^*, \b{H}_0^1(\Omega)} = \left\langle \bu, \bv \right\rangle_{\left( \left[ \mathcal{C}^\infty_0(\Omega) \right]^d \right)^*, \left[ \mathcal{C}^\infty_0(\Omega) \right]^d}.
	\]}The Stokes operator $\b A:\bVs\to\b V_\sigma^*$ is defined as
	the canonical Riesz isomorphism of $\bVs$, i.e.
	\[
	\ip{\b A\bv}{\b w}_{\b V^*_\sigma, \b V_\sigma}:=
	(\nabla\bv, \nabla\b w)_{\b H}\,, \qquad \bv,\b w\in\bVs\,.
	\]
	Employing the spectral properties of the operator $\b{A}$,  as customary, we also define the family of operators $\b{A}^s$ for any $s \in \mathbb{R}$. In particular, if $\{ \beta_k, \b{e}_k\}_{k \in \enne_+} \subset \erre \times \bVs$ denote the eigencouples of $\b{A}$, where $\{\b{e}_k\}_{k \in \enne_+}$ is an orthonormal basis of $\bHs$ and an orthogonal basis of $\bVs$, we introduce for any $s \geq 0$ the family of Hilbert spaces
	\[
	D(\b{A}^s) := \left\{ \bv \in \bHs: \bv = \sum_{i=1}^\infty c_i \b{e}_i \text{ and }
	\sum_{i=1}^\infty \beta_i^{2s}|c_i|^2 < +\infty \right\},
	\]
	and we set $D(\b{A}^{-s}) = D(\b{A}^{s})^*$.
	Next, for all $s \geq0$, we define the operators
	\[
	\b{A}^s: D(\b{A}^s) \to \bHs, \qquad \b{v} = \sum_{i=1}^\infty c_i \b{e}_i \mapsto
	\b{A}^s\bv := \sum_{i=1}^\infty \beta_i^s c_i \b{e}_i.
	\]
	Hereafter, we recall a number of standard facts:
	\begin{enumerate}[(i)] \itemsep0.2em
		\item if $s = 1$, then the Hilbert space $D(\b{A}) =\{\bv \in \bHs: \b{A}\bv \in \bHs\} = \b{H}^2(\OO) \cap \bVs$ denotes the so-called part of $\b{A}$ in $\bHs$;
		\item if $s = \frac{1}{2}$, then we have $D(\b{A}^\frac{1}{2}) = \bVs$ and $D(\b{A}^{-\frac{1}{2}}) = \bVsd$;
		\item if $s = 0$, then $\b{A}^0$ is the identity operator in $\bHs$ so that $D(\b{A}^0) = \bHs$;
		\item if $s = -1$, then $\b{A}^{-1}$ coincides with the inverse of the Stokes operator on $\bVsd$ and extends it on $D(\b{A}^{-1})$.
	\end{enumerate}
	In light of the previous considerations, using $\bHs$ as pivot space, we also have the general variational structure
	\[
	D(\b{A}^s) \embed D(\b{A}^t) \embed \bHs \equiv D(\b{A}^0) \embed  D(\b{A}^{-t}) \embed D(\b{A}^{-s})
	\]
	for any $s > t > 0$, with dense and compact embeddings in two and three spatial dimensions. Finally, we remind that, owing to the Korn inequality, we have
	\[
	\norm{\bv}_{\bVs}=\norm{\nabla\bv}_{\b H}\leq\sqrt2\norm{D\bv}_{\b H}
	\leq\sqrt2\norm{\nabla\bv}_{\b H} \qquad\forall\,\bv\in\bVs\,,
	\]
	where $D\bv=\frac12(\nabla\bv+(\nabla\bv)^t)$ denotes the symmetric gradient.
	Furthermore, we define the usual Stokes trilinear form
	$b$ on $\bVs \times \bVs \times \bVs$
	\[
	b(\bu, \b v, \b w):=\int_{\OO}(\b u\cdot\nabla)\b v\cdot\b w
	=\sum_{i,j=1}^d\int_\OO u_i\frac{\partial v_i}{\partial x_j}w_j\,,
	\qquad \b u, \b v, \b w \in \bVs\,,
	\]
	and the associated bilinear form $B:\bVs\times\bVs\to \b V_\sigma^*$ as
	\[
	\ip{B(\bu, \bv)}{\b w}_{\b V_\sigma^*, \bVs}:=b(\bu, \b v, \b w)\,,
	\qquad \b u, \b v, \b w \in \bVs\,.
	\]
	Let us recall that $b(\bu, \bv, \b w)=-b(\bu, \b w, \bv)$ for all
	$\b u, \b v, \b w \in \bVs$, from which it follows in particular
	that $b(\bu, \bv, \bv)=0$ for all $\b u, \b v\in \bVs$.
	Moreover, we point out that thanks to the usual functional embeddings it holds that
	$B:\bVs\times\bVs\to \b L^{\andrea{\frac65}}(\OO)$, hence, in particular, that
	$B:\bVs\times\bVs\to \b V_1^*$.\\
	We now report for the reader's convenience a basic embedding result and its proof. This will be useful in the forthcoming analysis.
	\begin{lem} \label{lem:gn}
		Let $r > 1$ and let $X$ be a Banach space. For every $p > 1$, there exists
		$\alpha = \alpha(p,r)\in(0,1)$ such that
		$W^{1,r}(0,T;X) \embed W^{\alpha,p}(0,T;X)$.
		In particular, if $p \leq r$ then $\alpha$ is any quantity in
		$ (0,1)$, and if $p > r$ then $\alpha = \frac{r}{p}$.
	\end{lem}
	\begin{proof}
		The embedding holds trivially for every $\alpha \in (0,1)$ if $p = r$. The same follows in the case $1 < p < r$ from the chain of embeddings
		\[
		W^{1,r}(0,T;X) \embed W^{1,p}(0,T;X) \embed W^{\alpha,p}(0,T;X).
		\]
		Let now $p > r$. If $\alpha \in (0,1)$, $q \in [1,+\infty]$ satisfy
		\[
		\dfrac{1}{p} = \dfrac{1-\alpha}{q} + \dfrac{\alpha}{r},
		\]
		then the fractional Gagliardo-Nirenberg inequality (see \cite[Theorem 1]{Brezis18}) entails that
		\[
		\|f\|_{W^{\alpha,p}(0,T)} \leq C\|f\|_{L^q(0,T)}^{1-\alpha}\|f\|_{W^{1,r}(0,T)}^\alpha,
		\]
		for any $f \in W^{1,r}(0,T)$.
		Taking into account the embedding
		\[
		W^{1,r}(0,T) \embed C^0([0,T]) \embed L^\infty(0,T),
		\]
		valid for every $r > 1$, we infer that the right hand side of the inequality is finite for every $q \in [1,+\infty]$. Moreover, we also get
		\[
		\|f\|_{W^{\alpha,p}(0,T)} \leq C\|f\|_{W^{1,r}(0,T)}.
		\]
		If we set $q = +\infty$, then we have
		\[
		\alpha = \frac{r}{p}.
		\]
		If $u \in W^{1,r}(0,T;X)$ the claim follows applying the proved inequality to $t\mapsto
		f(t) = \|u(t)\|_X$.
	\end{proof} \noindent
	Finally, we shall make precise the rigorous interpretation of the stochastic terms (see \eqref{eq:ACNS1} and \eqref{eq:ACNS3}).
    As a cylindrical process on $U_i$, $i=1,2$, $W_i$ admits the following representation
	\begin{equation} \label{eq:representation}
		W_i = \sum_{k=0}^{+\infty} b_k u^i_k,
	\end{equation}
	where $\{b_k\}_{k \in \enne}$ is a family of real and independent Brownian motions. However, it is well known that \eqref{eq:representation} does not converge in $U_i$, in general. That being said, it always exists a larger Hilbert space $U_0^i$, such that
	$U_i \embed U_0^i$ with Hilbert-Schmidt embedding $\iota_i$, such that we can identify $W_i$ as a $Q^0_i$-Wiener process on $U_0^i$, for some trace-class operator $Q^0_i$ (see \cite[Subsections 2.5.1]{LiuRo}). Actually, it holds that $Q^0_i = \iota_i \circ \iota_i^*$. In the following, we may implicitly assume this extension by simply saying that $W_i$ is a cylindrical process on $U_i$. This holds also for stochastic integration with respect to $W_i$. The symbol
	\[
	\int_0^\cdot B(s)\,\d W_i(s) := \int_0^\cdot B(s) \circ \iota_i^{-1}(s)\,\d W_i(s),
	\]
	for every progressively measurable process $B \in L^2(\Omega; L^2(0,T;\cL^2(U,K)))$, where $K$ is any (real) Hilbert space. The definition is well posed and does not depend on the choice of $U_i^0$ or $\iota_i$ (see \cite[Subsection 2.5.2]{LiuRo}).
	\subsection{Structural assumptions}
	The following assumptions are in order throughout the paper.
	\begin{enumerate}[label = \textbf{(A\arabic*)}, ref = \textbf{(A\arabic*)}]
		\item\label{hyp:potential} The potential $F:[-1,1]\to[0,+\infty)$
		is of class $C^0([-1,1]) \cap C^2(-1,1)$ with $F'(0)=0$ and satisfies
		\[
		\lim_{x\to(\pm1)^{\mp}}F'(x)=\pm\infty.
		\]
		Furthermore, there exists $c_F > 0$ such that
		\[
		F''(x) \geq -c_F, \quad x \in (-1,1).
		\]
		\item\label{hyp:G1} The operator
		$G_1: \bHs \to \cL^2(U_1, \bHs)$ is linearly bounded in $\bHs$,
		namely there exists $C_{G_1} > 0$ such that
		\[
		\|G_1(\bv)\|_{\cL^2(U_1,\bHs)} \leq C_{G_1}\left( 1 + \|\bv\|_\bHs \right)
		\]
		for any $\bv \in \bHs$. Moreover, taking $Y = \bHs$ or $Y = \bVsd$,
		we assume that $G_1: Y \to \cL^2(U_1, Y)$ is $L_1$-Lipschitz-continuous
		for some positive constant $L_1$.\vspace{0.2cm}
		\item\label{hyp:G2} Setting $\mathcal B$ as the closed unit ball in $L^\infty(\OO)$,
		the operator $G_2:\mathcal B\to \cL^2(U_2,H)$ satisfies
		\[
		G_2(\psi)[u_k^2]=g_k(\psi) \qquad\forall\,k\in\enne_+\quad\forall\,\psi\in\mathcal B\,,
		\]
		where the sequence $\{g_k\}_{k\in\enne_+}\subset W^{1,\infty}(-1,1)$ is such that
		\[
		g_k(\pm1)=0\,, \quad
		F''g_k^2\in L^\infty(-1,1) \qquad\forall\,k\in\enne_+\,,
		\]
		and
		\[
		L_2^2:=\sum_{k=1}^\infty\left(\norm{g_k}_{W^{1,\infty}(-1,1)}^2 +
		\norm{F''g_k^2}_{L^\infty(-1,1)}\right)<+\infty\,.
		\]
		In particular, note that this implies that $G_2:\mathcal B\to\cL^2(U_2,H)$
		is $L_2$-Lipschitz-continuous with respect to the $H$-metric on $\mathcal B$,
		and also $G_2(\mathcal B\cap V_1)\subset\cL^2(U_2,V_1)$.
		With a slight abuse of notation, we will use the symbol
		\[
		\nabla G_2:\mathcal B\cap V_1\to \cL^2(U_2,\bH)
		\]
		to indicate the operator
		\[
		\nabla G_2(\psi)[u_k^2]:=\nabla g_k(\psi)=g_k'(\psi)\nabla\psi\,,
		\quad k\in\enne_+\,,\quad \psi\in\mathcal B\cap V_1\,.
		\]
	\end{enumerate}
	\begin{remark}
		Let us point out that the physically relevant choice of $F$ (see \eqref{eq:fhpot})
		satisfies \ref{hyp:potential} and the compatibility condition in \ref{hyp:G2}, up to a suitable extension by right (or left) continuity at the boundary of $[-1,1]$ 
        and some additive constant to grant positivity (see, e.g., \cite[Remark 2.3]{scarpa21}).
	\end{remark} \noindent
	\begin{remark}
		If $Y = \bHs$ in \ref{hyp:G1}, then linear boundedness is directly implied by Lipschitz continuity.
	\end{remark}
	\subsection{Main results} We first introduce suitable notions of solution for problem \eqref{eq:ACNS1}-\eqref{eq:ACNS5}. The first is a martingale
solution, the second is a probabilistically-strong solution.
	\begin{defin}
		\label{def:mart_sol}
		Let $p \geq 1$ and let $(\bu_0, \varphi_0)$ satisfy
		\begin{align}
			\label{eq:u0}
			\bu_0&\in L^p(\Omega,\cF_0; \bHs)\,,\\
			\label{eq:phi0}
			\varphi_0&\in L^p(\Omega,\cF_0; V_1)\,, \quad
			F(\varphi_0)\in L^{p/2}(\Omega,\cF_0; L^1(\OO))\,.
		\end{align}
		A martingale solution to problem \eqref{eq:ACNS1}-\eqref{eq:ACNS5}
		with respect to the initial datum $(\bu_0, \varphi_0)$
		is a family
		\[
		\left(\left(\widehat\Omega, \widehat\cF, (\widehat\cF_t)_{t\in[0,T]}, \widehat\P\right),
		 \widehat W_1, \widehat W_2,
		\widehat\bu, \widehat\varphi\right)\,,
		\]
		where: $(\widehat\Omega, \widehat\cF, (\widehat\cF_t)_{t\in[0,T]}, \widehat\P)$ is a filtered
		probability space satisfying the usual conditions; $\widehat W_1, \widehat W_2$
		are two independent cylindrical Wiener processes on $U_1$ and $U_2$, respectively; the pair of processes $(\widehat \bu, \widehat\varphi)$ satisfies
		\begin{align}
			\label{eq:u_hat}
			&\widehat\bu\in L^p_w(\widehat\Omega; L^\infty(0,T; \bHs)) \cap
			L^p_\cP(\widehat\Omega; L^2(0,T; \bVs))\,,\\
			\label{eq:phi_hat}
			&\widehat\varphi \in L^p_\cP(\widehat\Omega; C^0([0,T]; H))\cap
			L^p_w(\widehat\Omega; L^\infty(0,T; V_1)) \cap
			L^p_\cP(\widehat\Omega; L^2(0,T; V_2))\,,\\
			& |\widehat{\varphi}(\omega, x, t)| < 1 \text{ for a.a. } (\omega, x,t) \in \hom \times \OO \times (0,T)\,, \\
			\label{eq:mu_hat}
			&\widehat\mu:=-\Delta\widehat\varphi+F'(\widehat\varphi) \in
			L^{p}_\cP(\widehat\Omega; L^2(0,T; H))\,,\\
			\label{eq:initial_hat}
			&(\widehat \bu(0), \widehat \varphi(0))\stackrel{\mathcal L}{=}(\bu_0, \varphi_0)
			\quad\text{on } \bHs\times V_1\,;
		\end{align}
		and, for every $\bv \in \bVs$ and $\psi\in V_1$, it holds that
		\begin{align}
			\nonumber
			&(\widehat \bu(t),\bv)_{\bHs} +
			\int_0^t\left[\ip{\b A\widehat \bu(s)}{\bv}_{\b V_\sigma^*, \b V_\sigma}
			+\ip{\b B(\widehat \bu(s), \widehat \bu(s))}{\bv}_{\b V_\sigma^*, \b V_\sigma}
			-\int_\OO\widehat \mu(s)\nabla\widehat \varphi(s)\cdot\bv
			\right]\,\d s\\
			\label{eq:var1_hat}
			&\qquad= (\widehat \bu(0),\bv)_{\bHs} +
			\left(\int_0^t G_1(\widehat \bu(s))\,\d \widehat W_1(s), \bv\right)_{\bHs}
			\qquad\forall\,t\in[0,T]\,,\:\:\widehat \P\text{-a.s.}\\
			\nonumber
			&(\widehat \varphi(t),\psi)_H +
			\int_0^t\!\int_\OO\left[\widehat \bu(s)\cdot\nabla\widehat \varphi(s) + \widehat \mu(s)\right]\psi\,\d s\\
			\label{eq:var2_hat}
			&\qquad= (\widehat \varphi(0),\psi)_{H} +
			\left(\int_0^t G_2(\widehat \varphi(s))\,\d \widehat W_2(s), \psi\right)_{H}
			\qquad\forall\,t\in[0,T]\,,\: \:\widehat \P\text{-a.s.}
		\end{align}
	\end{defin}
	
	\begin{defin}
		\label{def:strong_sol}
		Let $p \geq 1$ and let $(\bu_0, \varphi_0)$ satisfy
		\eqref{eq:u0}-\eqref{eq:phi0}.
		A probabilistically-strong solution to problem \eqref{eq:ACNS1}--\eqref{eq:ACNS5}
		with respect to the initial datum $(\bu_0, \varphi_0)$
		is a pair of processes $(\bu, \varphi)$ such that
		\begin{align}
			\label{u}
			&\bu\in L^p_w(\widehat\Omega; L^\infty(0,T; \bHs))\cap
			L^p_\cP(\widehat\Omega; L^2(0,T; \bVs))\,,\\
			\label{phi}
			&\varphi \in L^p_\cP(\Omega; C^0([0,T]; H))\cap
			L^p_w(\Omega; L^\infty(0,T; V_1)) \cap
			L^p_\cP(\Omega; L^2(0,T; V_2))\,,\\
			& |{\varphi}(\omega, x, t)| < 1 \text{ for a.a. } (\omega, x,t) \in \Omega \times \OO \times (0,T)\,, \\
			\label{mu}
			&\mu:=-\Delta\varphi+F'(\varphi) \in L^{p}_\cP(\Omega; L^2(0,T; H))\,,\\
			\label{initial}
			&(\bu(0), \varphi(0))=(\bu_0, \varphi_0)\,,
		\end{align}
		and
		\begin{align}
			\nonumber
			&(\bu(t),\bv)_{\bHs} +
			\int_0^t\left[\ip{\b A \bu(s)}{\bv}_{\b V_\sigma^*, \b V_\sigma}
			+\ip{\b B( \bu(s),  \bu(s))}{\bv}_{\b V_\sigma^*, \b V_\sigma}
			-\int_\OO \mu(s)\nabla \varphi(s)\cdot\bv
			\right]\,\d s\\
			\label{var1}
			&\qquad= (\bu_0,\bv)_{\bHs} +
			\left(\int_0^t G_1( \bu(s))\,\d  W_1(s), \bv\right)_{\bHs}
			\qquad\forall\,t\in[0,T]\,,\quad \P\text{-a.s.}\\
			\nonumber
			&( \varphi(t),\psi)_H +
			\int_0^t\!\int_\OO\left[ \bu(s)\cdot\nabla \varphi(s) +  \mu(s)\right]\psi\,\d s\\
			\label{var2}
			&\qquad= ( \varphi_0,\psi)_{H} +
			\left(\int_0^t G_2( \varphi(s))\,\d  W_2(s), \psi\right)_{H}
			\qquad\forall\,t\in[0,T]\,,\quad \P\text{-a.s.}
		\end{align}
		for every $\bv \in \bVs$ and $\psi\in V_1$.
	\end{defin}

	\begin{remark}
	Note that in Definitions~\ref{def:mart_sol} and \ref{def:strong_sol} one has in particular
	that $\widehat\bu\in \andrea{C^0_\text{w}([0,T]; \bH_{\sigma})}$ $\widehat\P$-almost surely and
	$\bu\in \andrea{C^0_\text{w}([0,T]; \bH_{\sigma})}$ $\P$-almost surely, respectively. \andrea{Here, the subscript ``w'' stands for weak continuity in time.} Thus the initial conditions \eqref{eq:initial_hat} and \eqref{initial} make sense.
	\end{remark}

	\noindent
	The first main result is the existence of a martingale solution.
	\begin{thm}
		\label{th:martingalesols}
		Assume \ref{hyp:potential}-\ref{hyp:G2} and let $p>2$. Then,
		for every initial datum $(\bu_0,\varphi_0)$ satisfying \eqref{eq:u0}-\eqref{eq:phi0}
		there exists a martingale solution
		$((\widehat\Omega, \widehat\cF, (\widehat\cF_t)_{t\in[0,T]}, \widehat\P), \widehat W_1, \widehat W_2,
		\widehat\bu, \widehat\varphi)$ to problem \eqref{eq:ACNS1}--\eqref{eq:ACNS6}
		satisfying the energy inequality
		\begin{align}
		\nonumber
		&\dfrac{1}{2}\supp\hE\|\widehat\bu(\tau)\|^2_{\bHs}
		+ \dfrac{1}{2}\supp\hE\|\nabla \widehat\varphi(\tau)\|_\b{H}^2
		+ \supp\hE\|F(\widehat\varphi(\tau))\|_{L^1(\OO)} \\
		\nonumber
		&\qquad+ \hE\int_0^t \left[ \|\nabla\widehat\bu(s)\|^2_{\bHs}
		+ \|\widehat\mu(s)\|^2_H\right]\mathrm{d}s \\
		\nonumber
		&\leq \left(C_{G_1}^2 + \frac{L_2^2}2|\OO|\right)t
		+ \dfrac{1}{2}\hE\|\widehat\bu_{0}\|^2_{\bHs}
		+ \dfrac{1}{2}\hE\|\nabla \widehat\varphi_{0}\|_\b{H}^2
		+ \hE\|F(\widehat\varphi_{0})\|_{L^1(\OO)} \\
		&\qquad+ C_{G_1}^2\hE\int_0^t \|\widehat\bu(\tau)\|^2_{\bHs} \: \mathrm{d}\tau
		+ \frac{L_2^2}2\hE\int_0^t  \norm{\nabla \widehat\varphi(\tau)}_{\bH}^2 \, \mathrm{d}\tau
	\label{eq:energy1}
	\end{align}
	for every $t\in[0,T]$, $\widehat \P$-almost surely. Here $|\OO|$ stands for the Lebesgue measure of $\OO$.
	Furthermore, there exists $\widehat{\pi} \in L^{\frac p2}_\cP(\hom; W^{-1, \infty}(0,T; H))$ such that
	\andrea{
		\begin{align}
		\nonumber
			&\int_0^T\ip{\partial_t(\hbu - G_1(\hbu)\cdot\widehat W_1)(t)
			+\b L \hbu(t)
			+\b B( \hbu(t),  \hbu(t))}{\bv(t)}_{\andrea{(\b{H}^1_0(\OO))^*, \b{H}^1_0(\OO)}} \,\d t\\
		\label{eq:press}
			& \quad = (\hbu_0,\bv)_{\bHs} +
			\int_0^T\int_\OO \hmu(t)\nabla \hphi(t) \cdot\bv(t)\,\d t
			+\ip{\widehat{\pi}}{\operatorname{div}\bv}_{\left( \left [ C^\infty_0((0,T)\times\OO)\right]^d\right)^* ,\left[ C^\infty_0((0,T)\times\OO) \right]^d}
		\end{align}
	}
	for every $\bv \in \andrea{\left[ C^\infty_0((0,T)\times\OO) \right]^d}$, $\hP$-almost surely.
	Finally, the following estimate holds:
	\begin{align}
	\nonumber
		\|\widehat{\pi}\|_{L^{\frac p2}(\hom;W^{-1,\infty}(0,T;H))} &\leq
		C\left(1 + \|\hbu\|_{L^{\frac p2}(\hom;L^\infty(0,T;\bHs))}
		+ \|\hbu\|_{L^{\frac p2}(\hom;L^2(0,T;\bVs))}
		+ \|\hbu\|_{L^p(\hom;L^2(0,T;\bVs))}^2 \right. \\
	\label{ineq:press}
		&\qquad\quad\left. + \|\hphi\|_{L^p(\hom;L^2(0,T;V_2))}^2
		+ \|F'(\hphi)\|_{L^p(\hom;L^2(0,T;H))}^2\right).
	\end{align}
	\end{thm} \noindent
\begin{remark}
The above result still holds if the viscosity depends on $\varphi$ in a smooth way and it is bounded from below by a positive constant.
Moreover, we recall that, in \cite{feir-petc2}, the energy inequality is written $\P\text{-a.s.}$ in a distributional sense.
\end{remark}
	The second is a stronger result in dimension two, namely, the existence of a (unique) probabilistically-strong solution.
	\begin{thm}
		\label{th:probstrongsol}
		Assume \ref{hyp:potential}-\ref{hyp:G2}, let $d=2$, $p>2$, and $Y = \bVsd$ in \ref{hyp:G1}. Then,
		for every initial datum $(\bu_0,\varphi_0)$ satisfying \eqref{eq:u0}--\eqref{eq:phi0},
		there exists a unique probabilistically-strong solution
		$(\bu, \varphi)$ for problem \eqref{eq:ACNS1}--\eqref{eq:ACNS6}
		and a pressure $\pi \in L^{\frac p2}_\cP(\Omega; W^{-1, \infty}(0,T; H))$,
		which satisfy on the original probability space $(\Omega, \cF,\P)$
		the analogous of the energy inequality \eqref{eq:energy1},
		the pressure-variational formulation \eqref{eq:press}, and
		the estimate \eqref{ineq:press}.
	\end{thm}

\begin{remark}
Referring to \cite{feir-petc2}, we observe that a more general $G_1(\bu,\varphi)$ can be considered. Instead, considering $G_2(\bu,\varphi)$
would require appropriate assumptions on account of \ref{hyp:G2}. For instance, in place of $g_k(\psi)$ we could suppose to
have $g_k(\psi)h_k(\bu)$ for a suitable $\{h_k\}_{k\in\enne_+}$.
\end{remark}

	\section{Proof of Theorem~\ref{th:martingalesols} }
	\label{sec:proof1}
	\noindent
	Here we prove the existence of martingale solutions to problem \eqref{eq:ACNS1}--\eqref{eq:ACNS6}.
    For the sake of clarity, the proof is split into several steps.
	\subsection{Regularization of the singular potential} \label{ssec:lambda}
	First of all, note that assumption \ref{hyp:potential} implies that the function
	\[
	\gamma:(-1,1)\to\erre\,, \qquad \gamma(x):=F'(x) + c_Fx\,, \quad x\in(-1,1)\,,
	\]
	can be identified with a maximal monotone graph in $\erre\times\erre$.
	Consequently, one can consider, for every $\lambda\in(0,1)$,
	the resolvent operator and the Yosida approximation of $\beta$, defined as follows
	\[
	J_\lambda, \: \gamma_\lambda:\erre\to\erre\,,\qquad
	J_\lambda(x):=(I + \lambda\gamma)^{-1}(x)\,, \qquad
	\gamma_\lambda(x):=\lambda^{-1}(x-J_\lambda(x))\,, \qquad x\in\erre\,.
	\]
	For notation and general properties of monotone operators we refer the reader to \cite{barbu-monot}.
	For every $\lambda\in(0,1)$, we define an approximation of $F$ as follows
	\begin{equation} \label{eq:F_lam}
		F_\lambda:\erre\to[0,+\infty)\,, \qquad
		F_\lambda(x):=F(0) + \ii{0}{x}{\gamma_\lambda(s)}{s}- \frac{c_F}{2}x^2\,, \quad x\in\erre\,.
	\end{equation}
	Thus it holds
	\begin{equation} \label{eq:derivativereg}
		F_\lambda'(x)=\gamma_\lambda(x) - c_Fx \quad\forall\,x\in\erre\,.
	\end{equation}
	In order to preserve the scaling of the Yosida-approximation on $F'$,
	we analogously define the $\lambda$-approximation of $G_2$ by setting
	\begin{equation} \label{eq:G2_lam}
		G_{2,\lambda}:=G_2\circ J_\lambda:H\to\cL^2(U_2,H).
	\end{equation}
	Notice that, by assumption \ref{hyp:G2} and the non-expansivity of $J_\lambda$,
	the operator $G_{2,\lambda}$ is $L_2$-Lipschitz-continuous
	(therefore uniformly in $\lambda$), and converges pointwise to $G_2$ as $\lambda\to0^+$.
	Now, we consider the $\lambda$-approximated (formal) problem
	\begin{align}
		\label{eq:1_lam}
		\d\bu_\lambda +\left[-\Delta\bu_\lambda + (\bu_\lambda\cdot\nabla)\bu_\lambda + \nabla p_\lambda
		- \mu_\lambda\nabla\varphi_\lambda\right]\,\d t
		= G_1(\bu_\lambda)\,\d W_1
		\qquad&\text{in } (0,T)\times\OO,\\
		\label{eq:2_lam}
		\nabla \cdot \bu_\lambda = 0 \qquad&\text{in } (0,T)\times\OO,\\
		\label{eq:3_lam}
		\d\varphi_\lambda + \left[\bu_\lambda\cdot\nabla\varphi_\lambda + \mu_\lambda\right]\,\d t =
		G_{2,\lambda}(\varphi_\lambda)\,\d W_2
		\qquad&\text{in } (0,T)\times\OO,\\
		\label{eq:4_lam}
		\mu_\lambda = -\Delta\varphi_\lambda + F_\lambda'(\varphi_\lambda)
		\qquad&\text{in } (0,T)\times\OO,\\
		\label{eq:5_lam}
		\bu_\lambda=0, \quad \partial_\bn\varphi_\lambda=0  \qquad&\text{in } (0,T)\times\partial\OO,\\
		\label{eq:6_lam}
		\bu_\lambda(0)=\bu_0\,, \quad \varphi_\lambda(0)=\varphi_0 \qquad&\text{in } \OO.
	\end{align}
	\subsection{Faedo-Galerkin approximation} \label{ssec:galerkin}
	A discretization scheme is now applied to problem \eqref{eq:1_lam}-\eqref{eq:6_lam}. Let us consider the (countably many) eigencouples of the negative Laplace operator with homogeneous Neumann boundary condition, namely the couples $\{(\alpha_j, e_j)\}_{j \in \mathbb{N_+}} \subset \erre \times V_2 $ such that
	\[
	\begin{cases}
		-\Delta e_j = \alpha_j e_j, & \quad \text{in } \OO,\\
		\partial_\bn e_j = 0 & \quad \text{on } \partial\OO,
	\end{cases} \qquad j \in \enne_+.
	\]
	Analogously, we also consider the (countably many) eigencouples of the Stokes operator, namely the couples $\{(\beta_k, \b{e}_k)\}_{k \in \mathbb{N}_+} \subset \mathbb{R} \times \bVs$, and $\{\pi_k\}_{k \in \mathbb{N}_+} \subset L^2_0(\mathbb{\OO})$ such that
	\[
	\begin{cases}
		-\Delta \b{e}_k + \nabla \pi_k = \beta_k \b{e}_k & \quad \text{in } \OO,\\
		\nabla \cdot \b{e}_k = 0 & \quad \text{in } \OO,\\
		\b{e}_k = 0 & \quad \text{on } \partial \OO,
	\end{cases} \qquad k \in \enne_+.
	\]
	It is well known that, up to a renormalization, the set $\{e_j\}_{j \in \mathbb{N}_+}$ (resp. $\{\b{e}_k \}_{k \in \mathbb{N}_+}$) is an orthonormal system in $H$ (resp. $\bHs$) and an orthogonal system in $V_1$ (resp. $\bVs$).  Let $n \in \enne_+$ and consider the finite-dimensional spaces
	$Z_n := \operatorname{span}\{e_1, \dots, e_n\}$ and
	$\b{Z}_n := \operatorname{span}\{\b{e}_1, \dots, \b{e}_n\}$,
	both endowed with the $L^2$-norm. In order to approximate the stochastic perturbation, we define the operators $G_{1,n}$ and $G_{2,\lambda,n}$ as
	\[
	G_{1,n}: \b{Z}_n \to \cL^2(U_1,\bHs), \qquad G_{2,\lambda,n}: Z_n \to \cL^2(U_2,H)
	\]
	and such that
	\[
	G_{1,n}(\b{v})u_k^1 := \sum_{j=1}^n (G_1(\b{v})u_k^1, \b{e}_j)_\bHs\b{e}_j, \qquad
	G_{2,\lambda,n}(v)u_k^2 := \sum_{j=1}^n (G_{2, \lambda}(v)u_k^2, e_j)_He_j
	\]
	for any $k \in \enne$, $\b{v} \in \b{Z}_n$ and
	$v \in Z_n$. Notice that, fixed any
	$\b{v} \in \b{Z}_n$ and $v \in Z_n$,
	$G_{1,n}(\b{v})$ and $G_{2, \lambda, n}(v)$ are actually well defined as elements of
	$\cL^2(U_1, \b{Z}_n)$ and $\cL^2(U_2, Z_n)$, respectively. Indeed, for instance,
	\begin{equation} \label{eq:lipschitz}
		\begin{split}
			\|G_{1,n}(\b{v})\|^2_{\cL^2(U_1, \b{Z}_n)} =
			\|G_{1,n}(\b{v})\|^2_{\cL^2(U_1, \bHs)}& =
			\sum_{k=1}^{+\infty} \|G_{1,n}(\b{v})u_k^1\|_\bHs^2 \\
			& = \sum_{k=1}^{+\infty} \sum_{j=1}^n |(G_1(\b{v})u_k^1, \b{e_j})_\bHs|^2 \\ & \leq
			 \sum_{k=1}^{+\infty} \sum_{j=1}^{+\infty} |(G_1(\b{v})u_k^1, \b{e_j})_\bHs|^2
			 = \|G_1(\b{v})\|^2_{\cL^2(U_1, \bHs)}.
		\end{split}	
	\end{equation}
	Moreover, since $G_1$ is $L_1$-Lipschitz continuous in the sense of assumption \ref{hyp:G1}
	and the orthogonal projection on $\b Z_n$ is non-expansive
	as an operator from $\bVs$ to itself, we can deduce by the same argument that
	$G_{1,n}$ is also $L_1$-Lipschitz continuous as an operator from $Y$ to $\cL^2(U_1,Y)$.
	Similar considerations also apply to $G_{2, \lambda, n}$. More precisely, we have
	\begin{prop} \label{prop:lipschitz}
		Let $\lambda \in (0,1)$ and $n \in \enne_+$.
		The operators
		\[
		G_{1,n}: \b{Z}_n \to \cL^2(U_1,\b{Z}_n), \qquad G_{2,\lambda,n}: Z_n \to \cL^2(U_2,Z_n)
		\]
		are well defined
		and uniformly Lipschitz continuous with respect to $n$ and $\lambda$.
		In particular, $G_{1,n}$ is $L_1$-Lipschitz continuous from $Y$ to $\cL^2(U_1,Y)$
		and $G_{2,\lambda,n}$ is $L_2$-Lipschitz continuous from $H$ to $\cL^2(U_1,H)$.
	\end{prop} \noindent
	Next, we define suitable projections (orthogonal with respect to the $L^2$-inner products)
	of initial data (evaluated at some point in $\Omega$) on the discrete spaces
	$Z_n$ and $\b{Z}_n$, namely, for all $n \in \enne_+$, we set
	\[
	\bu_{0,n} = \sum_{j=1}^n (\bu_0, \b{e}_j)_\bHs \b{e}_j, \qquad
	\varphi_{0,n} = \sum_{j=1}^n (\varphi_0, e_j)_He_j.
	\]
	It is now possible to formulate the discretized problem, which reads
	\begin{align}
		\label{eq:1_disc}
		\d\bu_{\lambda,n} + \left[-\Delta\bu_{\lambda,n} + (\bu_{\lambda,n}\cdot\nabla)\bu_{\lambda,n} + \nabla p_{\lambda,n}
		- \mu_{\lambda,n}\nabla\varphi_{\lambda,n}\right]\,\d t = G_{1,n}(\bu_{\lambda,n})\,\d W_1
		\qquad&\text{in } (0,T)\times\OO,\\
		\label{eq:2_disc}
		\nabla \cdot \bu_{\lambda,n} = 0 \qquad&\text{in } (0,T)\times\OO,\\
		\label{eq:3_disc}
		\d\varphi_{\lambda,n} + \left[\bu_{\lambda,n}\cdot\nabla\varphi_{\lambda,n} + \mu_{\lambda,n}\right]\,\d t =
		G_{2,\lambda,n}(\varphi_{\lambda,n})\,\d W_2
		\qquad&\text{in } (0,T)\times\OO,\\
		\label{eq:4_disc}
		\mu_{\lambda,n} = -\Delta\varphi_{\lambda,n} + F_{\lambda}'(\varphi_{\lambda,n})
		\qquad&\text{in } (0,T)\times\OO,\\
		\label{eq:5_disc}
		\bu_{\lambda,n}=0, \quad \partial_\bn\varphi_{\lambda,n}=0  \qquad&\text{in } (0,T)\times\partial\OO,\\
		\label{eq:6_disc}
		\bu_{\lambda,n}(0)=\bu_{0,n}\,, \quad \varphi_{\lambda,n}(0)=\varphi_{0,n} \qquad&\text{in } \OO.
	\end{align}
	The variational formulation of problem \eqref{eq:1_disc}-\eqref{eq:6_disc} is given by
	{\small
	\begin{align}
		\nonumber
		&(\bu_{\lambda,n}(t),\bv)_{\bHs} +
		\int_0^t\left[\ip{\b A \bu_{\lambda,n}(s)}{\bv}_{\b V_\sigma^*, \b V_\sigma}
		+\ip{\b B( \bu_{\lambda,n}(s),  \bu_{\lambda,n}(s))}{\bv}_{\b V_\sigma^*, \b V_\sigma}
		-\int_\OO \mu_{\lambda,n}(s)\nabla \varphi_{\lambda,n}(s)\cdot\bv
		\right]\,\d s\\
		\label{eq:vardisc1}
		&\qquad= (\bu_{0,n},\bv)_{\bHs} +
		\left(\int_0^t G_{1,n}( \bu_{\lambda,n}(s))\,\d  W_1(s), \bv\right)_{\bHs}
		\qquad\forall\,t\in[0,T]\,,\quad \P\text{-a.s.}\\
		\nonumber
		&( \varphi_{\lambda,n}(t),\psi)_H +
		\int_0^t\!\int_\OO\left[ \bu_{\lambda,n}(s)\cdot\nabla \varphi_{\lambda,n}(s) +  \mu_{\lambda,n}(s)\right]\psi\,\d s\\
		\label{eq:vardisc2}
		&\qquad= ( \varphi_{0,n},\psi)_{H} +
		\left(\int_0^t G_{2,\lambda,n}( \varphi_{\lambda,n}(s))\,\d  W_2(s), \psi\right)_{H}
		\qquad\forall\,t\in[0,T]\,,\quad \P\text{-a.s.}
	\end{align}}\noindent
	for every $\bv \in \b{Z}_n$ and $\psi\in Z_n$.
	Fixed any $\lambda \in (0,1)$ and $n \in \enne_+$, we search
	for a weak solution to \eqref{eq:vardisc1}-\eqref{eq:vardisc2} of the form
	\begin{equation} \label{eq:galerkin}
		\bu_{\lambda,n} = \sum_{j=1}^n a^j_{\lambda,n}\b{e}_j, \qquad
		\varphi_{\lambda,n} = \sum_{j=1}^n b^j_{\lambda,n}e_j, \qquad
		\mu_{\lambda,n} = \sum_{j=1}^n c^j_{\lambda,n}e_j,
	\end{equation}
	where
	\begin{align*}
		\b{a}_{\lambda,n} = (a^1_{\lambda,n}, a^2_{\lambda,n}, ..., a^n_{\lambda,n}) :
		\Omega \times [0,T] \to \mathbb{R}^n, \\
		\b{b}_{\lambda,n} = (b^1_{\lambda,n}, b^2_{\lambda,n}, ..., b^n_{\lambda,n}) :
		\Omega \times [0,T] \to \mathbb{R}^n, \\
		\b{c}_{\lambda,n} = (c^1_{\lambda,n}, c^2_{\lambda,n}, ..., c^n_{\lambda,n}) :
		\Omega \times [0,T] \to \mathbb{R}^n,
	\end{align*}
	are suitable stochastic processes.
	Inserting \eqref{eq:galerkin} into \eqref{eq:vardisc1}-\eqref{eq:vardisc2} and choosing as test functions $\psi = e_i$ and $\b{v} = \b{e}_i$ for each $i \in \{1, \dots, n\}$, we deduce that the three processes $\b{a}_{\lambda,n}$, $\b{b}_{\lambda,n}$ and $\b{c}_{\lambda,n}$ satisfy the system of $3n$ ordinary stochastic differential equations
	\begin{align}
		\nonumber
		\mathrm{d}a^i_{\lambda,n} + \beta_ia^i_{\lambda,n}
		+ \sum_{j=1}^n \sum_{k = 1}^n a^j_{\lambda,n}a^k_{\lambda,n}b(\b{e}_j, \b{e}_k, \b{e}_i)
		- \sum_{j=1}^n \sum_{k = 1}^n c^j_{\lambda,n}b^k_{\lambda,n}\int_\OO e_j\nabla e_k \cdot \b{e}_i&\\
		\label{eq:gal1}
		= \left( G_{1,n}\left( \sum_{j=1}^n a^j_{\lambda,n}\b{e}_j \right)\mathrm{d}W_1, \b{e}_i \right)_\bHs&\\
		\label{eq:gal2}
		\mathrm{d}b^i_{\lambda,n} + c^i_{\lambda,n} = \left( G_{2,\lambda,n}\left( \sum_{j=1}^n b^j_{\lambda,n}e_j \right)\mathrm{d}W_2, e_i \right)_H\\
		\label{eq:gal3}
		c^i_{\lambda,n} = \alpha_ib^i_{\lambda,n} + \int_\OO F'_\lambda\left( \sum_{j=1}^n b^j_{\lambda,n}e_j\right)e_i, \\
		\label{eq:gal4}
		a^i_{\lambda,n}(0) = (\b{u}_0, \b{e}_i)_\bHs \\
		\label{eq:gal5}
		b^i_{\lambda,n}(0) = (\varphi_0, e_i)_H
	\end{align}
	Let us point out that, in order to derive \eqref{eq:gal1}-\eqref{eq:gal5}, we exploited the fact that, for every choice of integers $j$ and $k$ between $1$ and $n$,
	\begin{equation} \label{eq:basis1}
		\ii{\OO}{}{\b{e}_j \cdot \nabla e_k}{x}  = -\ii{\OO}{}{e_k \nabla \cdot \b
		{e}_j}{x} +\ii{\partial \OO}{}{e_k \b
		{e}_j \cdot \b{n}}{\sigma} = 0,
	\end{equation}
	as well as the orthogonality in $\bVs$ of the basis $\{\b{e}_j\}_{j \in \mathbb{N}}$.
	The stochastic integrals in \eqref{eq:gal1}-\eqref{eq:gal2} have to be regarded as $G^i_{1,\lambda,n}\,\mathrm{d}W_1$ and $G^i_{2, \lambda,n}\,\mathrm{d}W_2$ for every $i = 1, \dots, n$, where
	\[
	G^i_{1,\lambda,n} : \b{Z}_n \to \cL^2(U_1, \mathbb{R}), \qquad 	
	G^i_{1,\lambda,n}(\b{u}_{\lambda,n})u_k^1 :=
	\left( G_{1,n}\left( \sum_{j=1}^n a^j_{\lambda,n}\b{e}_j \right)u^1_k, \b{e}_i \right)_\bHs
	\]
	and
	\[
	G^i_{2,\lambda,n} : Z_n \to \cL^2(U_2, \mathbb{R}), \qquad 	
	G^i_{2,\lambda,n}(\varphi_{\lambda,n})u_k^2 :=
	\left( G_{2,\lambda,n}\left( \sum_{j=1}^n b^j_{\lambda,n}e_j \right)u_k^2, e_i \right)_H,
	\]
	for every $k \in \enne$. By Lipschitz continuity of all the nonlinearities appearing in  \eqref{eq:gal1}-\eqref{eq:gal5}, the standard theory of abstract stochastic evolution equations applies. Therefore, we are able to infer that
	\begin{prop} \label{prop:galerkin}
		For every $\lambda \in (0,1)$ and $n \in \enne_+$,
		there exists a unique triplet of $(\cF_t)_t$-adapted processes
		$\b{a}_{\lambda,n}$, $\b{b}_{\lambda,n}$, $\b{c}_{\lambda,n}$
		satisfying problem \eqref{eq:gal1}-\eqref{eq:gal5}. Furthermore, for every $r \geq 2$, we have
		\[
		\b{a}_{\lambda,n}, \b{b}_{\lambda,n}, \b{c}_{\lambda,n} \in L^r(\Omega; C^0([0,T]; \erre^n)),
		\]
		implying
		\[
		\b{u}_{\lambda,n} \in L^r(\Omega; C^0([0,T]; \b{Z}_n)),
		\quad\varphi_{\lambda,n}, \mu_{\lambda,n} \in L^r(\Omega; C^0([0,T]; Z_n)).
		\]
	\end{prop}
	\subsection{Uniform estimates with respect to $n$}
	\label{ssec:est}
	First of all, we prove some uniform estimates with respect to the Galerkin parameter $n$, keeping $\lambda \in (0,1)$ fixed. Hereafter, the symbol $C$ (possibly numbered) denote positive constants whose special dependencies are explicitly pointed out when necessary. `In some cases, in order to ease notation, we may use the same symbol to denote different constants throughout the same argument. In any case, such constants are always independent of $n$.
	\paragraph{\textit{First estimate}}
		We exploit the It\^{o} formula for the $H$-norm of $\varphi_{\lambda,n}$ given in \cite[Theorem 4.2.5]{LiuRo}. This gives
		\begin{multline} \label{eq:unif10}
			\dfrac{1}{2}\|\varphi_{\lambda,n}(t)\|_{H}^2 +\int_0^t \left[ \|\nabla \varphi_{\lambda,n}(\tau)\|^2_\b{H} + \left( \varphi_{\lambda,n}(\tau), F'_\lambda(\varphi_{\lambda,n}(\tau))\right)_H \right] \: \mathrm{d}\tau \\ = \dfrac{1}{2}\|\varphi_{0,n}\|_{H}^2 + \int_0^t \left(\varphi_{\lambda,n}(\tau), G_{2, \lambda, n} (\varphi_{\lambda,n}(\tau))\,\mathrm{d}W_2(\tau)\right)_H + \dfrac{1}{2}\int_0^t \|G_{2,\lambda,n}(\varphi_{\lambda,n}(\tau))\|^2_{\cL^2(U_2, H)} \: \mathrm{d}\tau.
		\end{multline}
		Let us now address the above equality term by term. First of all, recalling
		\eqref{eq:derivativereg} and that $F_\lambda'(0)=0$, we find
		\begin{equation} \label{eq:unif11}
				\left( \varphi_{\lambda,n}(\tau), F'_\lambda(\varphi_{\lambda,n}(\tau))\right)_H
				\geq  - c_F\|\varphi_{\lambda,n}(\tau)\|^2_H.
		\end{equation}
		Next, owing to \eqref{eq:lipschitz} and \ref{hyp:G2}, we have
		\begin{equation} \label{eq:unif12}
			\begin{split}
				\|G_{2,\lambda,n}(\varphi_{\lambda,n}(\tau))\|^2_{\cL^2(U_2, H)} \leq \|G_{2,\lambda}(\varphi_{\lambda,n}(\tau))\|^2_{\cL^2(U_2, H)} & = \sum_{k=1}^{+\infty} \|g_k(J_\lambda(\varphi_{\lambda,n}(\tau)))\|^2_H \\ & \leq \sum_{k=1}^{+\infty} \|g_k\|^2_{W^{1,\infty}(-1,1)}|\OO|
				\\ & \leq L_2^2|\OO|.
			\end{split}			
		\end{equation}
		Finally, by 1-Lipschitz-continuity of the projection $\Pi_n : H \to H$, it follows
		\begin{equation} \label{eq:unif13}
			\|\varphi_{0,n}\|_{H}^2 \leq \|\varphi_{0}\|_{H}^2.
		\end{equation}
		Thus, combining \eqref{eq:unif11}-\eqref{eq:unif13} with \eqref{eq:unif10},
		letting $p \in [2,+\infty)$, multiplying the resulting inequality by two,
		taking $\frac{p}{2}$-powers, the supremum on the interval $[0,t]$ and expectations, we arrive at
		\begin{align*}
			&\E\supp\|\varphi_{\lambda,n}(\tau)\|_{H}^p +
			\E \left| \int_0^t \|\nabla \varphi_{\lambda,n}(\tau)\|^2_\b{H}  \: \mathrm{d}\tau \right|^\frac{p}{2} \\
			&\leq C \left[ 1+\E\|\varphi_{0}\|_{H}^p +
			\E \int_0^t \|\varphi_{\lambda,n}(\tau)\|_{H}^p  \: \mathrm{d}\tau
			+\E\sups\left| \int_0^s \left(\varphi_{\lambda,n}(\tau),
			G_{2, \lambda, n} (\varphi_{\lambda,n}(\tau))\right)_H\,\mathrm{d}W_2(\tau)
			 \right|^\frac{p}{2}  \right],
		\end{align*}
		where $C$ depends on $p$ and also on $c_F$, $L_2$, $|\OO|$, $T$.
		The Burkholder-Davis-Gundy and H\"older inequalities jointly with \eqref{eq:unif12} entail
		\begin{align}
		\nonumber
		& \E\sups\left| \int_0^s \left(\varphi_{\lambda,n}(\tau),
		G_{2, \lambda, n} (\varphi_{\lambda,n}(\tau))\,\mathrm{d}W_2(\tau)\right)_H
		\right|^\frac{p}{2} \\
		\nonumber
		&\leq C
		\E\left| \int_0^t \|\varphi_{\lambda,n}(\tau)\|^2_H
		\|G_{2, \lambda, n} (\varphi_{\lambda,n}(\tau))\|^2_{\cL^2(U_2,H)}\, \mathrm{d}\tau \right|^\frac{p}{4} \\
		\nonumber
		&\leq C\E\left| \supp \|\varphi_{\lambda,n}(\tau)\|^2_H
		\int_0^t\|G_{2, \lambda, n}  (\varphi_{\lambda,n}(\tau))\|^2_{\cL^2(U_2,H)}\,\mathrm{d}\tau
		\right|^\frac{p}{4} \\
		\label{eq:unif14}
		&\leq CL_2^\frac{p}{2}|\OO|^\frac{p}{4}t^\frac{p}{4}
		\E\supp \|\varphi_{\lambda,n}(\tau)\|^\frac{p}{2}_H,	
		\end{align}
		where $C$ only depends on $p$.
		In turn, thanks to \eqref{eq:unif14} and the Young inequality, we can refine the estimate and get
		\begin{align*}
			&\E\supp\|\varphi_{\lambda,n}(\tau)\|_{H}^p +
			\E \left| \int_0^t \|\nabla \varphi_{\lambda,n}(\tau)\|^2_\b{H}  \: \mathrm{d}\tau
			\right|^\frac{p}{2}
			\leq C \left[ 1+\E\|\varphi_{0}\|_{H}^p +
			\E \int_0^t \|\varphi_{\lambda,n}(\tau)\|_{H}^p  \: \mathrm{d}\tau  \right].
		\end{align*}
		The Gronwall lemma entails that there exists $C_1$,
		independent of $n$ and $\lambda$, but depending on $p$ and
		the structural data of the problem, such that
		\begin{equation} \label{eq:estimate1}
			\|\varphi_{\lambda,n}\|_{L^p_\cP(\Omega;C^0([0,T];H))} + \|\varphi_{\lambda,n}\|_{L^{p}_\cP(\Omega;L^2([0,T];V_1))} \leq C_1,
		\end{equation}
		for every fixed $p \geq 2$.
		\paragraph{\textit{Second estimate.}} We devise a similar argument for the $\bHs$-norm of $\b{u}_{\lambda,n}$. Still exploiting \cite[Theorem 4.2.5]{LiuRo}, the It\^{o} formula implies
		\begin{multline} \label{eq:unif20}
			\dfrac{1}{2}\|\b{u}_{\lambda,n}(t)\|^2_{\bHs} + \int_0^t \left[ \|\nabla\b{u}_{\lambda,n}(\tau)\|^2_{\bHs} - \mu_{\lambda,n}(\tau)\nabla\varphi_{\lambda,n}(\tau)\cdot\b{u}_{\lambda,n}(\tau)\right]\mathrm{d}\tau \\ = \dfrac{1}{2}\|\b{u}_{0,n}\|^2_{\bHs} +\int_0^t\left(\bu_{\lambda,n}(\tau), G_{1,n}( \bu_{\lambda,n}(\tau))\,\d  W_1(\tau)\right)_{\bHs} + \dfrac{1}{2}\int_0^t \|G_{1,n}(\b{u}_{\lambda,n}(\tau))\|^2_{\cL^2(U_1, \bHs)} \: \mathrm{d}\tau.
		\end{multline}
		Next, we want to apply the standard It\^{o} formula to the regularized energy functional
		$\mathcal{E}_\lambda: Z_n \times \b{Z}_n \to \mathbb{R}^+$
		\begin{equation*} \label{eq:regenergy}
			\mathcal{E}_\lambda(\varphi_{\lambda,n}, \b{u}_{\lambda,n})
			:= \dfrac{1}{2}\int_\OO |\nabla \varphi_{\lambda,n}|^2
			+ \dfrac{1}{2}\int_\OO |\b{u}_{\lambda,n}|^2 + \int_\OO F_\lambda(\varphi_{\lambda,n}).
		\end{equation*}
		However, notice that $\mathcal{E}_\lambda$ exactly contains the
		kinetic energy contribution linked to the fluid velocity field which
		we just handled in \eqref{eq:unif20}. Thus, it is sufficient to apply
		the It\^{o} formula only to the portion of the energy linked to the
		order parameter $\varphi_{\lambda,n}$. Let us stress that this is only
		possible since no coupling energy terms are present. We set
		\[
		\Psi_\lambda: Z_n \to \mathbb{R}^+, \qquad \Psi_\lambda(v)
		:= \dfrac{1}{2}\int_\OO |\nabla v|^2 + \int_\OO F_\lambda(v).
		\]
		It has already been shown in \cite[Subsection 3.2]{scarpa21} that $\Psi_\lambda$ is twice
		Fr\'{e}chet-differentiable.
		Thus it is possible to apply the It\^{o} formula in its classical version
		\cite[Theorem 4.32]{dapratozab}. This yields
		\begin{align}
		\nonumber
			&\Psi_\lambda(\varphi_{\lambda,n}(t)
			+ \int_0^t \left[\|\mu_{\lambda,n}(\tau)\|^2_H
			+\mu_{\lambda,n}(\tau)\nabla\varphi_{\lambda,n}(\tau)
			\cdot\b{u}_{\lambda,n}(\tau) \right]\, \mathrm{d}\tau \\
		\nonumber
			&= \Psi_\lambda(\varphi_{0,n}) +
			\frac12
			\int_0^t \left[ \norm{\nabla G_{2,\lambda,n}(\varphi_{\lambda,n}(\tau))}_{\cL^2(U_2,\bH)}^2 +
			\sum_{k=1}^\infty\int_\OO F_\lambda''(\varphi_{\lambda,n}(\tau))
			|g_k(J_\lambda(\varphi_{\lambda,n}(\tau)))|^2\right]\, \mathrm{d}\tau \\
			&\qquad+
		\label{eq:unif21}
			\int_0^t\left( \mu_{\lambda,n}(\tau),
			G_{2,\lambda,n}(\varphi_{\lambda,n}(\tau))\,\d  W_2(\tau)\right)_{H},
		\end{align}
		where we recall that $D\Psi_\lambda(\varphi_{\lambda,n}) = \mu_{\lambda,n}$.
		Adding \eqref{eq:unif20} and \eqref{eq:unif21} together, we find
		\begin{align}
		\nonumber
			&\dfrac{1}{2}\|\b{u}_{\lambda,n}(t)\|^2_{\bHs}
			+ \Psi_\lambda(\varphi_{\lambda,n}(t))
			+ \int_0^t \left[ \|\nabla\b{u}_{\lambda,n}(\tau)\|^2_{\bHs}
			+ \|\mu_{\lambda,n}(\tau)\|^2_H \right]\mathrm{d}\tau \\
		\nonumber
			&= \dfrac{1}{2}\|\b{u}_{0,n}\|^2_{\bHs} +\Psi_\lambda(\varphi_{0,n})
			+ \dfrac{1}{2}\int_0^t \left[ \|G_{1,n}(\b{u}_{\lambda,n}(\tau))\|^2_{\cL^2(U_1, \bHs)}
			+ \norm{\nabla G_{2,\lambda,n}(\varphi_{\lambda,n}(\tau))}_{\cL^2(U_2,\bH)}^2
			\right]\,\mathrm{d}\tau \\
			\nonumber
			&\qquad
			+ \dfrac{1}{2}\int_0^t\sum_{k=1}^\infty\int_\OO F_\lambda''(\varphi_{\lambda,n}(\tau))
			|g_k(J_\lambda(\varphi_{\lambda,n}(\tau)))|^2\,\mathrm{d}\tau \\
			\label{eq:unif22}
			&\qquad
			+ \int_0^t\left(\bu_{\lambda,n}(\tau), G_1( \bu_{\lambda,n}(\tau))\,\d  W_1(\tau)\right)_{\bHs}
			+ \int_0^t\left( \mu_{\lambda,n}(\tau),
			G_{2,\lambda,n}(\varphi_{\lambda,n}(\tau))\,\d  W_2(\tau)\right)_{H}.
		\end{align}
		Fix now $p \in [2, +\infty)$. Taking $\frac{p}{2}$-powers, supremum over $[0,t]$,
		and expectations of both sides of \eqref{eq:unif22} yield
		\begin{align} \nonumber
			&\E\supp\|\b{u}_{\lambda,n}(\tau)\|^p_{\bHs}
			+ \E\supp\|\nabla \varphi_{\lambda,n}(\tau)\|^p_{\b{H}}
			+ \E\supp \|F_\lambda(\varphi_{\lambda,n})\|_{L^1(\OO)}^\frac{p}{2} \\
			\nonumber &\qquad
			+ \E\left|\int_0^t  \|\nabla\b{u}_{\lambda,n}(\tau)\|^2_{\bHs}\,\mathrm{d}\tau \right|^\frac{p}{2}
			+ \E\left|\int_0^t  \|\mu_{\lambda,n}(\tau)\|^2_{H}\,\mathrm{d}\tau \right|^\frac{p}{2} \\
			\nonumber &
			\leq C\left[\E\|\b{u}_{0,n}\|^p_{\bHs}
			+ \E\|\nabla \varphi_{0,n}\|^p_\b{H}
			+ \E\|F_\lambda(\varphi_{0,n})\|_{L^1(\OO)}^\frac{p}{2}
			+ \E\left|\int_0^t\|G_{1,n}(\b{u}_{\lambda,n}(\tau))\|^2_{\cL^2(U_1, \bHs)}
			\,\mathrm{d}\tau \right|^\frac{p}{2} \right. \\
			\nonumber &\qquad
			+ \E\left|\int_0^t
			\norm{\nabla G_{2,\lambda,n}(\varphi_{\lambda,n}(\tau))}_{\cL^2(U_2,\bH)}^2
			\,\mathrm{d}\tau \right|^\frac{p}{2}
			+\E\left|\int_0^t\sum_{k=1}^\infty\int_\OO
			|F_\lambda''(\varphi_{\lambda,n}(\tau))|
			|g_k(J_\lambda(\varphi_{\lambda,n}(\tau)))|^2\,\mathrm{d}\tau \right|^\frac{p}{2} \\
			\nonumber
			&\qquad \left.
			 +\E\sups\left|\int_0^s\left(\bu_{\lambda,n}(\tau),
			 G_1( \bu_{\lambda,n}(\tau))\,\d  W_1(\tau)\right)_{\bHs}\right|^\frac{p}{2}\right. \\
			 \label{eq:unif23}
			 &\qquad\left.+\E\sups\left| \int_0^s\left( \mu_{\lambda,n}(\tau),
			 G_{2,\lambda,n}(\varphi_{\lambda,n}(\tau))\,\d  W_2(\tau)\right)_{H}\right|^\frac{p}{2} \right],	
		\end{align}
		where $C$ only depends on $p$.
		Next, we address the terms on the right hand side of \eqref{eq:unif23}.
		By \eqref{eq:lipschitz} and Assumption \ref{hyp:G1}, we deduce
		\begin{equation} \label{eq:unif24}
			\|G_{1,n}(\b{u}_{\lambda,n}(\tau))\|^2_{\cL^2(U_1, \bHs)} \leq
			\|G_{1}(\b{u}_{\lambda,n}(\tau))\|^2_{\cL^2(U_1, \bHs)} \leq
			2C_{G_1}^2\left( 1 + \|\b{u}_{\lambda,n}(\tau))\|^2_\bHs\right).
		\end{equation}
		Since $\varphi_{\lambda,n}(\tau) \in V_1$, recalling
		assumption \ref{hyp:G2}, \eqref{eq:unif12}, and the non-expansivity of $J_\lambda$,
		we have
		\begin{align}
        \nonumber
	    \norm{\nabla G_{2,\lambda,n}(\varphi_{\lambda,n}(\tau))}_{\cL^2(U_2,\bH)}^2
			&\leq
			\norm{\nabla G_{2,\lambda}(\varphi_{\lambda,n}(\tau))}_{\cL^2(U_2,\bH)}^2\\
            \nonumber
			& =  \sum_{k=1}^{\infty}\|g_k'(J_\lambda(\varphi_{\lambda,n}(\tau)))
			\nabla J_\lambda(\varphi_{\lambda,n}(\tau))\|_{H}^2 \\
            \label{eq:unif25}
			& \leq\sum_{k=1}^{\infty}\|g_k\|_{W^{1,\infty}(-1,1)}^2
			\|\nabla \varphi_{\lambda,n}(\tau)\|_{H}^2
			 \leq L_2^2  \|\nabla \varphi_{\lambda,n}(\tau)\|_{H}^2.
		\end{align}
		Furthermore, since $F''=\gamma'-c_F$,
		by \eqref{eq:derivativereg} we have that, for all $x\in\erre$,
		\[
		F''_\lambda(x)=\gamma_\lambda'(x)-c_F
		=\gamma'(J_\lambda(x))J_\lambda'(x) - c_F
		=F''(J_\lambda(x))J_\lambda'(x) + c_F(J_\lambda'(x)-1).
		\]
		Thus, thanks to \ref{hyp:G2} and the non-expansivity of $J_\lambda$, we get
		\begin{align}
		\nonumber
			\sum_{k=1}^\infty\int_\OO |F_\lambda''(\varphi_{\lambda,n}(\tau))|
			|g_k(J_\lambda(\varphi_{\lambda,n}(\tau)))|^2
			&\leq
			\sum_{k=1}^\infty\int_\OO
			|F''(J_\lambda(\varphi_{\lambda,n}(\tau)))|
			|g_k(J_\lambda(\varphi_{\lambda,n}(\tau)))|^2\\
			\nonumber
			&\qquad+ 2c_F|\OO|\sum_{k=1}^\infty \|g_k\|_{W^{1,\infty}(-1,1)}^2\\
			&\leq |\OO|L_2^2\left(1+2c_F\right).
			\label{eq:unif26}
		\end{align}
		Finally, we address the stochastic integrals.
		Using \eqref{eq:unif24} jointly with the Burkholder-Davis-Gundy and
		Young inequalities, for every $\delta>0$ we obtain
		\begin{align}
		\nonumber
			& \E\supp\left|\int_0^t\left(\bu_{\lambda,n}(\tau),
			G_{1,n}( \bu_{\lambda,n}(\tau))\,\d  W_1(\tau)\right)_{\bHs}\right|^\frac{p}{2} \\
		\nonumber
			&\leq C\E\left| \int_0^t \|\bu_{\lambda,n}(\tau)\|^2_\bHs
			\|G_{1,n}( \bu_{\lambda,n}(\tau))\|^2_{\cL^2(U_1,\bHs)}\, \mathrm{d}\tau \right|^\frac{p}{4} \\
		\nonumber
			&\leq  C2^\frac{p}{4}C_{G_1}^\frac{p}{2}\E\left| \supp \|\bu_{\lambda,n}(\tau)\|^2_\bHs
			\int_0^t \left( 1 + \|\bu_{\lambda,n}(\tau)\|^2_\bHs\right)\mathrm{d}\tau\right|^\frac{p}{4}\\
		\label{eq:unif27}
			&\leq \delta\E\supp\|\bu_{\lambda,n}(\tau)\|^p_\bHs
			+ C_{p,\delta} \E \int_0^t \left( 1 + \|\bu_{\lambda,n}(\tau)\|^p_\bHs\right)
			\mathrm{d}\tau,
		\end{align}
		where $C$ only depends on $\delta$, $p$, and $T$.
		Moreover, by \eqref{eq:unif12} and the same inequalities, we also get
		\begin{align}
		\nonumber
		& \E\supp\left| \int_0^t\left( \mu_{\lambda,n}(\tau),
		G_{2,\lambda,n}(\varphi_{\lambda,n}(\tau))\,\d  W_2(\tau)\right)_{H}\right|^\frac{p}{2} \\
		\nonumber
		&\leq C \E\left| \int_0^t \|\mu_{\lambda,n}(\tau)\|^2_H
		\|G_{2, \lambda, n} (\varphi_{\lambda,n}(\tau))\|^2_{\cL^2(U_2,H)}\, \mathrm{d}\tau \right|^\frac{p}{4} \\
		\label{eq:unif28}
		&\leq C L_2^\frac{p}{2}|\OO|^\frac{p}{4}\E\left| \int_0^t
		\|\mu_{\lambda,n}(\tau)\|^2_H\, \mathrm{d}\tau \right|^\frac{p}{4}
		\leq
		C
		+ \delta\E\left| \int_0^t \|\mu_{\lambda,n}(\tau)\|^2_H\, \mathrm{d}\tau \right|^\frac{p}{2},
		\end{align}
		where $C$ only depends on $p$, $\delta$, and $\OO$.
		Finally, the non-expansivity of the orthogonal projectors on $W_n$ and $\b{W}_n$ imply
		\begin{equation} \label{eq:unif29}
			\|\b{u}_{0,n}\|^p_{\bHs} \leq \|\b{u}_0\|_\bHs^p, \quad \|\nabla \varphi_{0,n}\|^p_\b{H} \leq \|\nabla \varphi_0\|_{\b{H}}^p
		\end{equation}
		whereas, since $F'_\lambda$ is linearly bounded, being Lipschitz-continuous, $F_\lambda$ is quadratically bounded so that
		\begin{equation} \label{eq:unif210}
			\|F_\lambda(\varphi_{0,n})\|_{L^1(\OO)}^\frac{p}{2}	\leq C\left(1 + \|\varphi_{0,n}\|_H^p\right)	\leq C\left(1 + \|\varphi_{0}\|_H^p\right),
		\end{equation}
		where $C$ depends on $\lambda$ and $p$.
		Collecting \eqref{eq:unif24}-\eqref{eq:unif210} and choosing $\delta$ small enough, from \eqref{eq:unif23} we infer that
		\begin{align}
		\nonumber
		&\E\supp\|\b{u}_{\lambda,n}(t)\|^p_{\bHs} +
		\E\supp\|\nabla \varphi_{\lambda,n}(\tau)\|^p_{\b{H}} +
		\E\supp \|F_\lambda(\varphi_{\lambda,n})\|_{L^1(\OO)}^\frac{p}{2} \\
		\nonumber
		&\qquad+
		\E\left|\int_0^t  \|\nabla\b{u}_{\lambda,n}(\tau)\|^2_{\bHs}\,\mathrm{d}\tau \right|^\frac{p}{2}
		+ \dfrac{1}{2}\E\left|\int_0^t  \|\mu_{\lambda,n}(\tau)\|^2_{H}\,\mathrm{d}\tau \right|^\frac{p}{2} \\
		\label{eq:unif211}
		&\leq C \left[1+\E\|\varphi_{0}\|_{V}^p
		+ \E\|\b{u}_{0}\|^p_{\bHs}
		+\E\int_0^t\|\b{u}_{\lambda,n}(\tau)\|^p_{ \bHs}\,\mathrm{d}\tau +
		\E\int_0^t\norm{\nabla \varphi_{\lambda,n}(\tau)}_{\bH}^p\,\mathrm{d}\tau\right].
		\end{align}
        Here $C$ depends on $\lambda$ and $p$.
		An application of the Gronwall lemma entails the existence of $C_2,C_3,C_4$, depending on $\lambda$, $p$ and $T$, such that
		\begin{align} \label{eq:estimate2}
			\|\b{u}_{\lambda,n}\|_{L^p_\cP(\Omega;C^0([0,T];\bHs))} + 	
			\|\b{u}_{\lambda,n}\|_{L^p_\cP(\Omega;L^2(0,T;V_1))} \leq C_2, \\
			\label{eq:estimate3}
			\|\varphi_{\lambda,n}\|_{L^p_\cP(\Omega;C^0([0,T];V_1))} \leq C_3, \\
			\label{eq:estimate4}
			\|\mu_{\lambda,n}\|_{L^p_\cP(\Omega;L^2(0,T;H))} +
			\|F_\lambda(\varphi_{\lambda,n})\|_{L^{p/2}_\cP(\Omega;C^0([0,T];L^1(\OO)))} \leq
			C_4.
		\end{align}
		\paragraph{\textit{Further estimates.}}
		The Lipschitz-continuity of $F'_\lambda$ and the fact that $F'_\lambda(0) = 0$ entail
		\[
		|F'_\lambda(\varphi_{\lambda,n}(t))| \leq C|\varphi_{\lambda,n}(t)|,
		\]
		for some $C$ only depending on $\lambda$. Therefore, thanks to \eqref{eq:estimate1} we also get the estimate
		\begin{equation} \label{eq:estimate5}
			\|F'_\lambda(\varphi_{\lambda,n})\|_{L^p_\cP(\Omega;L^2(0,T;H))} \leq C_5.
		\end{equation}
		Additionally, by comparison in \eqref{eq:4_disc}, we get
		\begin{equation} \label{eq:estimate6}
			\|\varphi_{\lambda,n}\|_{L^p_\cP(\Omega;L^2(0,T;V_2))} \leq C_6.
		\end{equation}
		Here $C_5$ or $C_6$ depend on $\lambda$, $p$, and $T$. In light of \eqref{eq:unif24}, \eqref{eq:unif25} and on account of \eqref{eq:estimate2} and \eqref{eq:estimate3}, we deduce
		\begin{align} \label{eq:estimate7}
			\|G_{1,n}(\b{u}_{\lambda,n})\|_{L^p_\cP(\Omega;L^\infty(0,T;\cL^2(U_1,\bHs)))}
			&\leq C_7,  \\
			\label{eq:estimate8}
			\|G_{2,\lambda,n}(\varphi_{\lambda,n})\|_{L^p_\cP(\Omega;L^\infty(0,T;\cL^2(U_2,V_1)))} + \|G_{2,\lambda,n}(\varphi_{\lambda,n})\|_{L^\infty(\Omega \times (0,T);\cL^2(U_2,H))} &\leq C_8,
		\end{align}
		Here, again, the constants $C_7$, $C_8$ depend on $\lambda$. As a consequence of \cite[Lemma 2.1]{fland-gat},
	 the following estimates on the It\^{o} integrals hold:
		\begin{align} \label{eq:estimate9}
			\left\| \int_0^\cdot G_{1,n}(\b{u}_{\lambda,n}(\tau))\,\mathrm{d}W_1(\tau)
			\right\|_{L^p(\Omega;W^{k,p}(0,T;\bHs))} &\leq C_9, \\
			\label{eq:estimate10}
			\left\| \int_0^\cdot G_{2,\lambda,n}(\varphi_{\lambda,n}(\tau))\,\mathrm{d}W_2(\tau)
			\right\|_{L^p(\Omega;W^{k,p}(0,T;V_1))\cap L^q(\Omega; W^{k,q}(0,T; H))}
			&\leq C_{10},
		\end{align}
		for every $k \in (0,\frac{1}{2})$ and $q\geq1$, where $C_9$ and $C_{10}$ depend on $\lambda, p, q, k$, and $T$.
		Estimates \eqref{eq:estimate9} and \eqref{eq:estimate10} enable
		us to carry out two comparison arguments. Let us interpret \eqref{eq:vardisc2}
		as an equality in $V^*_1$,
		\begin{align*}
		\langle \varphi_{\lambda,n}(t),\psi\rangle_{V_1^*, V_1} &= -
		\int_0^t\!\int_\OO\left[ \bu_{\lambda,n}(s)\cdot\nabla \varphi_{\lambda,n}(s)
		+  \mu_{\lambda,n}(s)\right]\psi\,\d s \\
		&\quad+ \langle \varphi_{0,n},\psi\rangle_{V_1^*, V_1} +
		\left(\int_0^t G_{2,\lambda,n}( \varphi_{\lambda,n}(s))\,\d  W_2(s), \psi\right)_{H}
		\end{align*}
		for all $\psi \in V_1$ such that $\|\psi\|_{V_1} = 1$,
		$t\in[0,T]$, $\P$-almost-surely. It is clear that, by the H\"{o}lder inequality,
		\begin{equation} \label{eq:unif31}
			\left| \int_\OO \psi\bu_{\lambda,n}\cdot\nabla \varphi_{\lambda,n}
			+  \psi\mu_{\lambda,n} \right| \leq  \|\bu_{\lambda,n}\|_\bHs
			\|\nabla \varphi_{\lambda,n}\|_{\mathbf{L}^4(\OO)}
			+ \|\mu_{\lambda,n}\|_H
		\end{equation}
		implying (see \eqref{eq:estimate2}, \eqref{eq:estimate4}, and
		\eqref{eq:estimate6})
		\begin{equation} \label{eq:unif32}
			\left\| \int_0^t\!\int_\OO\left[ \bu_{\lambda,n}(s)\cdot\nabla \varphi_{\lambda,n}(s)
			+  \mu_{\lambda,n}(s)\right]\d s\right\|_{L^p_\cP(\Omega; H^1(0,T;V_1^*))} \leq C	
		\end{equation}
		for some $C$ depending on $\lambda, p$ and $T$.
		Then, recalling that
		\[
		|\langle \varphi_{0,n},\psi\rangle_{V_1^*, V_1}| \leq \|\varphi_{0,n}\|_H \leq \|\varphi_{0}\|_H,
		\]
		and estimate \eqref{eq:estimate10} as well as Lemma \ref{lem:gn}, we find
		\begin{equation} \label{eq:estimate11}
			\|\varphi_{\lambda,n}\|_{L^p_\cP(\Omega;W^{\beta,p}(0,T;V_1^*))} \leq C_{11}
		\end{equation}
		for some $\beta = \beta(p) \in (\frac{1}{p}, \frac{1}{2})$ if $p>2$,
		and for all $\beta\in(0,\frac12)$ if $p=2$.
		The constant $C_{11}$ may depend on $\lambda, \beta, p$, and $T$.
		\begin{remark}\label{rem:beta}
			Observe that $\beta$ is always well defined. Here, we apply Lemma \ref{lem:gn} with $r = 2$ and $X = V_1^*$. If $\alpha$ denotes the Sobolev fractional exponent given by Lemma \ref{lem:gn}, then the following alternative holds:
			\begin{enumerate}[(a)]
				\item if $p = 2$, then any value of $\alpha \in (0,1)$ is valid,
				and therefore we can set an arbitrary $\beta \in (0,\frac{1}{2})$;
				\item if $p > 2$, then any value of $\alpha \in (0,\frac{2}{p}]$ is valid,
				and therefore we can set an arbitrary
				$\beta \in (\frac{1}{p}, \min(\frac{2}{p}, \frac{1}{2})) \subset (\frac{1}{p}, \frac{1}{2})$.
			\end{enumerate}
		\end{remark} \noindent
	Similarly, we consider the weak formulation of the discretized Navier--Stokes equation
		\begin{align*}
		 \langle \bu_{\lambda,n}(t),\bv\rangle_{\bVsd,\bVs} &= -
			\int_0^t\left[\ip{\b A \bu_{\lambda,n}(s)}{\bv}_{\b V_\sigma^*, \b V_\sigma}
			+\ip{\b B( \bu_{\lambda,n}(s),  \bu_{\lambda,n}(s))}{\bv}_{\b V_\sigma^*, \b V_\sigma}
			\right]\,\d s \\
			&\quad-\int_0^t\int_\OO \mu_{\lambda,n}(s)\nabla \varphi_{\lambda,n}(s)\cdot\bv\,\d s
			+ (\bu_{0,n},\bv)_{\bHs} +
			\left(\int_0^t G_{1,n}( \bu_{\lambda,n}(s))\,\d  W_1(s), \bv\right)_{\bHs}
		\end{align*}
		for all $\b{v} \in \bVs$ such that $\|\b{v}\|_{\bVs} = 1$, $t\in[0,T]$, $\P$-almost-surely. Then, we have
		\[
		|\langle \bu_{0,n},\bv\rangle_{\bVsd, \bVs}| \leq \|\bu_{0,n}\|_\bHs \leq \|\bu_{0}\|_\bHs.
		\]
		Owing to \eqref{eq:estimate2} and the continuity of $\b{A}$, we infer
		\[
		\left\|\int_0^\cdot\b A \bu_{\lambda,n}(s)\,\d s\right\|_{L^p_\cP(\Omega; H^1(0,T;\bVsd))}
		 \leq C,
		\]
		for some $C$ depending on $\lambda, p$ and $T$, but independent of $n$. Next, we recall the well-known inequality
		\begin{equation} \label{eq:trilinear}
			\|\b{B}( \bu_{\lambda,n},  \bu_{\lambda,n})\|_{\bVsd}
			\leq \|\bu_{\lambda,n}\|_\bHs^{2-\frac{d}{2}}\|\bu_{\lambda,n}\|_\bVs^\frac{d}{2}.
		\end{equation}
		Therefore, we find
		\[
		\left\|\int_0^\cdot\b B( \bu_{\lambda,n}(s),  \bu_{\lambda,n}(s))\,\d s
		\right\|_{L^{\frac p2}_\cP(\Omega; W^{1,\frac{4}{d}}(0,T;\bVsd))} \leq C.
		\]
		Furthermore, since by the H\"{o}lder, Gagliardo--Nirenberg and Young inequalities, we have
		\begin{equation}
		\label{eq:kortewegest}
				\left|\int_\OO \mu_{\lambda,n}\nabla\varphi_{\lambda,n}\cdot \b{v}
				\right|  \leq \|\mu_{\lambda,n}\|_H
				\|\nabla\varphi_{\lambda,n}\|_{\andrea{\b{L}^3(\OO)}}\|\bv\|_{\andrea{\b{L}^6(\OO)}}
				\leq		C\|\mu_{\lambda,n}\|_H\|\varphi_{\lambda,n}\|_{L^6(\OO)}^\frac{1}{2}
				\|\varphi_{\lambda,n}\|_{V_2}^\frac{1}{2}
		\end{equation}
		for both $d =2$ and $d = 3$. Thus we get
		\[
		\left\|\int_0^\cdot \mu_{\lambda,n}(s)\nabla \varphi_{\lambda,n}(s)\,\d s
		\right\|_{L^\frac p2_\cP(\Omega; W^{1,\frac{4}{3}}(0,T;\bVsd))} \leq C.
		\]
		Summing up, also owing to \eqref{eq:estimate8} and Lemma \ref{lem:gn}, we arrive at
		\begin{equation} \label{eq:estimate12}
			\|\bu_{\lambda,n}\|_{L^{\frac p2}_\cP(\Omega;
			W^{\gamma,p}(0,T;\bVsd))} \leq C_{12},
		\end{equation}
		for some $\gamma = \gamma(p) \in (\frac{1}{p}, \frac{1}{2})$ if $p>2$,
		and for all $\gamma\in(0,\frac12)$ if $p=2$.
		Here $C_{12}$ depends on $\lambda, \gamma, p$, and $T$.
		\begin{remark}
        \label{rem:gamma}
            Observe that $\gamma$ is always well defined.
			In this case, we apply Lemma \ref{lem:gn} with $r = \frac{4}{3}$
			and $X = \bVsd$. Let $\alpha$ denote once again the fractional
			Sobolev exponent given by Lemma \ref{lem:gn}. Given any $p > 2 > \frac{4}{3}$,
			we have that any value of $\alpha \in (0,\frac{4}{3p}]$ is valid,
			and therefore we can set an arbitrary
			$\gamma \in (\frac{1}{p}, \min(\frac{4}{3p}, \frac{1}{2})) \subset (\frac{1}{p}, \frac{1}{2})$.
			If $p=2$ then we get any value of $\alpha \in (0,\frac{2}{3}]$. Hence we can choose any $\gamma\in(0,\frac12)$.
		\end{remark} \noindent
		In the following, we assume that, given $p \geq 2$, the exponents $\beta = \beta(p)$ and $\gamma = \gamma(p)$ are fixed. Notice that if $p > 2$, then trivially $\beta$ and $\gamma$ are both greater than 1.
	\subsection{Passage to the limit as $n\to+\infty$}
	\label{ssec:lim_n}
	Owing to the previously proven uniform estimates, we now pass to the limit as $n\to +\infty$ keeping $\lambda \in (0,1)$ fixed. Let $p > 2$.
	\begin{lem} \label{lem:tight1}
		The family of laws of $(\bu_{\lambda,n})_{n \in \enne}$ is tight in the space $Z_\b{u}:=L^2(0,T;\bHs) \cap C^0([0,T]; D(\b{A}^{-\delta}))$ for any $\delta \in (0,\frac{1}{2})$. The family of laws of $(\varphi_{\lambda,n})_{n \in \enne}$ is tight in the space $Z_\varphi:=L^2(0,T;V_1) \cap C^0([0,T];H)$.
	\end{lem}
	\begin{proof}
		To prove the claims, we follow a standard argument (refer, for instance, to \cite[Subsection 3.3]{scarpa21} or \cite[Proposition 1]{Medjo21}). We first recall we have that the embeddings (see \cite[Corollary 5]{simon})
		\begin{align*}
			L^\infty(0,T;V_1) \cap W^{\beta,p}(0,T;V_1^*) & \embed C^0([0,T];H), & \quad L^\infty(0,T;\bHs) \cap W^{\gamma,p}(0,T;\bVsd) & \embed C^0([0,T];D(\b{A}^{-\delta})), \\
			L^2(0,T;V_2) \cap W^{\beta,p}(0,T;V_1^*) & \embed L^2(0,T;V_1), & \quad L^2(0,T;\bVs) \cap W^{\gamma,p}(0,T;\bVsd) & \embed L^2(0,T;\bHs),
		\end{align*}
		are compact (the intersection spaces are endowed with their canonical norm), since $\beta, \gamma > \frac{1}{p}$, $p > 2$. Here, $\delta \in (0, \frac{1}{2})$. Let us prove the first one, the other three cases being similar. For any $R > 0$, let $B_R$ denote the closed ball of radius $R$ in $L^\infty(0,T;V_1) \cap W^{\beta,p}(0,T;V_1^*)$. Then the Markov inequality, jointly with estimates \eqref{eq:estimate3} and \eqref{eq:estimate12}, implies
		\[
		\begin{split}
			\P\left\{ \varphi_{\lambda,n} \in B_R^C \right\} & = \P\left\{ \|\varphi_{\lambda,n}\|_{L^\infty(0,T;V_1) \cap W^{\beta,p}(0,T;V_1^*)} > R\right\} \\
			& \leq \dfrac{1}{R^p}\E\|\varphi_{\lambda,n}\|_{L^\infty(0,T;V_1) \cap W^{\beta,p}(0,T;V_1^*)}^p \\
			& \leq \dfrac{C_{\lambda}}{R^p}.
		\end{split}
		\]
		for some $C_{\lambda} > 0$ depending on $\lambda$> This yields
		\[
		\lim_{n\to+\infty} \sup_{n \in \enne} \, \P\left\{ \varphi_{\lambda,n} \in B_R^C \right\} = 0,
		\]
		so that the first claim is proven. The remaining claims can be proven analogously, replacing the spaces accordingly and exploiting the corresponding estimates.
	\end{proof} \noindent
	We now set
	\[
	G_{1,n}(\bu_{\lambda,n}) \cdot W_1 := \int_0^\cdot G_{1,n}(\b{u}_{\lambda,n}(\tau))\,\mathrm{d}W_1(\tau), \qquad G_{2,\lambda,n}(\varphi_{\lambda,n}) \cdot W_2 := \int_0^\cdot G_{2,\lambda,n}(\varphi_{\lambda,n}(\tau))\,\mathrm{d}W_2(\tau).
	\]
	With a little modification in the proof of Lemma \ref{lem:tight1}, we can also prove
	\begin{lem} \label{lem:tight2}
		The family of laws of $(G_{1,n}(\bu_{\lambda,n}) \cdot W_1)_{n \in \enne}$ is tight in the space $C^0([0,T];\bVsd)$. The family of laws of $(G_{2,\lambda,n}(\varphi_{\lambda,n}) \cdot W_2)_{n \in \enne}$ is tight in the space $C^0([0,T];H)$.
	\end{lem}
	\begin{proof}
		By \cite[Theorem 2.2]{fland-gat}, since $\beta p > 1$ and $\gamma p > 1$, we have that the embeddings
		\[
		W^{\beta,p}(0,T;V_1) \embed C^0([0,T]; H), \qquad W^{\gamma,p}(0,T;\bHs) \embed C^0([0,T]; \bVsd)
		\]
		are compact.
		The argument of the proof of Lemma \ref{lem:tight1}, recalling estimates \eqref{eq:estimate9} and \eqref{eq:estimate10}, is enough to conclude the claims.
	\end{proof} \noindent
	Next, we consider the constant sequences of cylindrical Wiener processes
	\[
	W_{1,n} \equiv W_1, \qquad W_{2,n} \equiv W_2.
	\]
	\begin{lem} \label{lem:tight3}
		The family of laws of $(W_{1,n})_{n \in \mathbb{N}}$ is tight in $C^0([0,T]; U_1^0)$. The family of laws of $(W_{2,n})_{n \in \mathbb{N}}$ is tight in $C^0([0,T]; U_2^0)$.
	\end{lem}
	\begin{proof}
		It directly follows from the fact that every measure on a complete separable metric space is tight.
	\end{proof} \noindent
	Finally, we consider the sequences of approximated initial conditions.
	\begin{lem} \label{lem:tight4}
		The family of laws of $(\bu_{0,n})_{n \in \enne}$ is tight in $\bVsd$. The family of laws of $(\varphi_{0,n})_{n \in \enne}$ is tight in $H$.
	\end{lem}
	\begin{proof}
		It is a third iteration of the proof of Lemma \ref{lem:tight1}, exploting the compact embeddings
		\[
		\bHs \embed \bVsd, \qquad V_1 \embed H,
		\]
		and the Markov inequality on closed balls of $\bHs$ and $V_1$, respectively.   
	\end{proof} \noindent
	As an immediate consequence of Lemmas \ref{lem:tight1}-\ref{lem:tight4}, we get that the family of laws of \[ (\bu_{\lambda,n}, \varphi_{\lambda,n},  G_{1,n}(\bu_{\lambda,n}) \cdot W_{1,n}, G_{2,\lambda,n}(\varphi_{\lambda,n}) \cdot W_{2,n}, W_{1,n}, W_{2,n}, \bu_{0,n}, \varphi_{0,n})_{n \in \enne} \]
	is tight in the product space
	\[
	Z_\b{u} \times Z_\varphi \times C^0([0,T]; \bVsd) \times   C^0([0,T]; H) \times C^0([0,T]; U_1^0) \times C^0([0,T]; U_2^0) \times \bVsd \times H.
	\]
	Owing to the Prokhorov and Skorokhod theorems (see \cite[Theorem 2.7]{ike-wata} and \cite[Theorem 1.10.4, Addendum 1.10.5]{vaa-well}), there exists a probability space $(\tom, \widetilde{\mathscr{F}},\tP)$ and a sequence of random variables $X_n:  (\tom, \tF)\to(\Omega, \mathscr{F})$ such that the law of $X_n$ is $\P$ for every $n \in \enne$, namely $\tP \circ X_n^{-1} = \P$ (so that composition with $X_n$ preserves laws), and the following convergences hold
	\begin{align*}
		\widetilde{\bu}_{\lambda,n} := \bu_{\lambda,n} \circ X_n \to \widetilde{\bu}_{\lambda} \quad & \text{in }Z_\b{u}=L^2(0,T;\bHs) \cap C^0([0,T]; D(\b{A}^{-\delta})), \,\, \P\text{-a.s.}; \\
		\widetilde{\varphi}_{\lambda,n}:= \varphi_{\lambda,n} \circ X_n \to \widetilde{\varphi}_{\lambda} \quad & \text{in }Z_\varphi:=L^2(0,T;V_1) \cap C^0([0,T];H), \,\, \P\text{-a.s.};\\
		\widetilde{I}_{\lambda,n} := (G_{1,n}(\bu_{\lambda,n}) \cdot W_{1,n}) \circ X_n \to \widetilde{I}_{\lambda} \quad & \text{in } C^0([0,T]; \bVsd), \,\, \P\text{-a.s.}; \\
		\widetilde{J}_{\lambda,n} := (G_{2,\lambda,n}(\varphi_{\lambda,n}) \cdot W_{2,n}) \circ X_n \to \widetilde{J}_{\lambda} \quad & \text{in }C^0([0,T];H), \,\, \P\text{-a.s.}; \\
		\widetilde{W}_{1,n} := W_{1,n} \circ X_n \to \widetilde{W}_{1}\quad & \text{in } C^0([0,T]; U_1^0) , \,\, \P\text{-a.s.}; \\
		\widetilde{W}_{2,n} := W_{2,n} \circ X_n \to \widetilde{W}_{2}\quad & \text{in }C^0([0,T]; U_2^0) , \,\, \P\text{-a.s.}; \\
		\widetilde{\bu}_{0,n} := \bu_{0,n} \circ X_n \to \widetilde{\bu}_0  \quad& \text{in } \bVsd , \,\, \P\text{-a.s.};\\
		\widetilde{\varphi}_{0,n} := \varphi_{0,n} \circ X_n \to \widetilde{\varphi}_0 \quad & \text{in }H, \,\, \P\text{-a.s.},
	\end{align*}
	for some limiting processes $\widetilde{\bu}_{\lambda}, \widetilde{\varphi}_{\lambda}, \widetilde{I}_{\lambda}, \widetilde{J}_{\lambda}, \widetilde{W}_{1}, \widetilde{W}_{2}, \widetilde{\bu}_0, \widetilde{\varphi}_0$ belonging to the specified spaces. Let us recall that, for the sake of what follows, if $(S, \mathcal{M}, \nu)$ is a finite positive measure space and $X$ is any Banach space, then the Bochner space $L^r(S; X)$ is reflexive if and only if $L^r(S, \nu)$ and $X$ are reflexive (see, for instance, \cite[Corollary 2, p. 100]{diestel-uhl}). By the previously proven uniform estimates and the preservation of laws under $X_n$, up to a subsequence which we do not relabel, the Vitali convergence theorem, the Eberlein-Smulian theorem and the Banach-Alaoglu theorem entail
	\begin{align*}
		\widetilde{\bu}_{\lambda,n} \to \widetilde{\bu}_{\lambda}
		\quad & \text{in } L^q(\tom;L^2(0,T;\bHs) \cap C^0([0,T]; D(\b{A}^{-\delta})))  \text{ if }q < p, \\
		\widetilde{\bu}_{\lambda,n} \rightharpoonup \widetilde{\bu}_{\lambda}
		\quad & \text{in } L^p(\tom;L^2(0,T;\bVs)),\\
		\widetilde{\bu}_{\lambda,n} \overset{\ast}{\rightharpoonup} \widetilde{\bu}_{\lambda}
		\quad & \text{in } L^p_w(\tom;L^\infty(0,T;\bHs)) \cap
		L^{\frac p2}(\tom;W^{\gamma,p}(0,T;\bVsd)),\\
		\widetilde{\varphi}_{\lambda,n} \to \widetilde{\varphi}_{\lambda}
		\quad & \text{in } L^q(\tom;L^2(0,T;V_1) \cap C^0([0,T];H)) \text{ if }q < p, \\
		\widetilde{\varphi}_{\lambda,n} \rightharpoonup \widetilde{\varphi}_{\lambda}
		\quad & \text{in } L^p(\tom;L^2(0,T;V_2)),\\
		\widetilde{\varphi}_{\lambda,n} \overset{\ast}{\rightharpoonup} \widetilde{\varphi}_{\lambda}
		\quad & \text{in } L^p_w(\tom;L^\infty(0,T;V_1)) \cap L^p(\tom; W^{\beta,p}(0,T;V_1^*)), \\
		\widetilde{I}_{\lambda,n} \to \widetilde{I}_{\lambda}
		\quad & \text{in } L^q(\tom;C^0([0,T];\bVsd)) \text{ if }q < p, \\
		\widetilde{J}_{\lambda,n} \to \widetilde{J}_{\lambda}
		\quad & \text{in } L^q(\tom;C^0([0,T];H)) \text{ if }q < p, \\
		\widetilde{W}_{1,n}  \to \widetilde{W}_{1}
		\quad & \text{in } L^q(\tom;C^0([0,T];U_1^0)) \text{ if }q < p, \\
		\widetilde{W}_{2,n}  \to \widetilde{W}_{2}
		\quad & \text{in } L^q(\tom;C^0([0,T];U_2^0)) \text{ if }q < p, \\
		\widetilde{\bu}_{0,n} \to \widetilde{\bu}_0
		\quad & \text{in } L^q(\tom;\bVsd) \text{ if }q < p, \\
		\widetilde{\varphi}_{0,n} \to \widetilde{\varphi}_0
		\quad & \text{in } L^q(\tom;H) \text{ if }q < p.
	\end{align*}
	Let us now define
	\[
	\widetilde{\mu}_{\lambda,n} := \mu_{\lambda,n} \circ X_n.
	\]
	By uniform boundedness and weak compactness, there exists some $\widetilde{\mu}_{\lambda}$ such that
	\begin{align*}
		\widetilde{\mu}_{\lambda,n} \rightharpoonup \widetilde{\mu}_{\lambda} & \quad  \text{in } L^p(\tom;L^2(0,T;H)).
	\end{align*}
	Let us notice that it is possible to take the probability space $(\tom, \tF,\tP)$ large enough so that it does not depend on $\lambda$. Taking into account the previous considerations and further straightforward weak convergences, the limit processes fulfill the following regularity properties:
	\begin{align*}
		\widetilde{\bu}_\lambda & \in L^{\frac p2}(\tom; W^{\gamma,p}(0,T;\bVsd)) \cap
		L^p(\tom; C^0([0,T];D(\b{A}^{-\delta})) \cap L^2(0,T;\bVs))
		\cap L^p_w(\tom;L^\infty(0,T;\bHs)) ; \\
		\widetilde{\varphi}_\lambda & \in L^p(\tom; W^{\beta,p}(0,T;V_1^*)\cap C^0([0,T];H)
		\cap L^2(0,T;V_2))
		\cap L^p_w(\tom;L^\infty(0,T;V_1))  ; \\
		\widetilde{\mu}_\lambda & \in L^p(\tom;L^2(0,T;H)); \\
		\widetilde{I}_\lambda & \in L^p(\tom;C^0([0,T];\bVsd)); \\
		\widetilde{J}_\lambda & \in L^p(\tom;C^0([0,T];H)); \\
		\widetilde{W}_1 & \in L^p(\tom;C^0([0,T];U_1^0)); \\
		\widetilde{W}_2 & \in L^p(\tom;C^0([0,T];U_2^0)); \\
		\widetilde{\bu}_0 & \in L^p(\tom;\bHs);\\
		\widetilde{\varphi}_0 & \in L^p(\tom;\mathcal{B} \cap V_1).
	\end{align*}
	From this starting point, we now address several issues.
	\paragraph{\textit{The nonlinearities}} First of all, by Lipschitz-continuity of $F_\lambda'$, it follows that
	\[
	F_\lambda'(\widetilde{\varphi}_{\lambda,n}) \to F_\lambda'(\widetilde{\varphi}_{\lambda}) \quad \text{in } L^p(\tom; L^2(0,T;H)).
	\]
	Moreover, since $G_{1,n}$ is uniformly Lipschitz-continuous (recall Proposition \ref{prop:lipschitz}) and
	\[
	\begin{split}
		& \|G_{1,n}(\widetilde{\bu}_{\lambda,n})-G_{1}(\widetilde{\bu}_{\lambda})\|_{L^p(\tom, L^2(0,T;\cL^2(U_1,Y)))} \\
		& \qquad \qquad \leq \|G_{1,n}(\widetilde{\bu}_{\lambda,n})-G_{1,n}(\widetilde{\bu}_{\lambda})\|_{L^p(\tom, L^2(0,T;\cL^2(U_1,Y)))} + \|G_{1,n}(\widetilde{\bu}_{\lambda})-G_{1}(\widetilde{\bu}_{\lambda})\|_{L^p(\tom, L^2(0,T;\cL^2(U_1,Y)))},
	\end{split}
	\]
	we conclude
	\[
	G_{1,n}(\widetilde{\bu}_{\lambda,n}) \to G_{1}(\widetilde{\bu}_{\lambda}) \quad \text{in }L^q(\tom; L^2(0,T;\cL^2(U_1,Y))) \text{ if } q < p,
	\]
	A very similar computation also shows
	\[
	G_{2,\lambda,n}(\widetilde{\varphi}_{\lambda,n}) \to G_{2,\lambda}(\widetilde{\varphi}_{\lambda}) \quad \text{in }L^q(\tom; L^2(0,T;\cL^2(U_2,H))) \text{ if } q < p.
	\]
	Next, we address the Korteweg term representing the capillary force. Let us prove that \andrea{
	\[
	\widetilde\mu_{\lambda,n}\nabla\widetilde\varphi_{\lambda,n} \rightharpoonup
	\widetilde\mu_{\lambda}\nabla\widetilde\varphi_{\lambda}
	\quad\text{in } \b{L}^1(\tom\times(0,T)\times\OO).
	\]
	Indeed, for any $\b{w} \in \b{L}^\infty(\tom\times(0,T)\times\OO)$,
	\[
	\begin{split}
		& \left| \tE\int_{\OO \times (0,T)} (\widetilde\mu_{\lambda,n}\nabla\widetilde\varphi_{\lambda,n}
		-\widetilde\mu_{\lambda}\nabla\widetilde\varphi_{\lambda} ) \cdot \b{w} \right| \\
		& \leq \left| \tE\int_{\OO \times (0,T)}
		\widetilde\mu_{\lambda,n}(\nabla\widetilde\varphi_{\lambda,n} -\nabla\widetilde\varphi_{\lambda} )
		\cdot \b{w} \right| +
		 \left|\tE \int_{\OO \times (0,T)}
		 (\widetilde\mu_{\lambda,n} -\widetilde\mu_{\lambda})\nabla\widetilde\varphi_{\lambda}
		 \cdot \b{w} \right| \\
		& \leq \|\b{w}\|_{\b{L}^\infty(\tom\times(0,T)\times\OO)}
		\|\widetilde\mu_{\lambda,n}\|_{L^2(\tom\times(0,T)\times\OO)}
		\|\nabla\widetilde\varphi_{\lambda,n}
		-\nabla\widetilde\varphi_{\lambda}\|_{\b{L}^2(\tom\times(0,T)\times\OO)}
		+ \left| \tE\int_{\OO \times (0,T)} (\widetilde\mu_{\lambda,n} -\widetilde\mu_{\lambda})
		\nabla\widetilde\varphi_{\lambda} \cdot \b{w} \right|
	\end{split}
	\]
	and both terms tend to zero as $n \to +\infty$ by the above convergences
	(note that $\nabla \widetilde\varphi_{\lambda} \cdot \b{w} \in L^2(\tom\times(0,T)\times\OO)$). Here, $\tE$ stands for the expectation with respect to the probability $\tP$.} As far as the other nonlinear term appearing in the Navier-Stokes equations, we have, as a straightforward application of \eqref{eq:trilinear},
	\[
	\b{B}(\bu_{\lambda,n}, \bu_{\lambda,n}) \to \b{B}(\bu_{\lambda}, \bu_{\lambda})
	\text{ in } L^q(\tom;L^\frac{4}{d}(0,T;\bVsd)) \text{ if } q < \frac p2.
	\]
	Finally, we address the convective term. Observe that
	\[
	\begin{split}
		\widetilde{\b{u}}_{\lambda,n} \cdot \nabla \widetilde{\varphi}_{\lambda,n}
		- \widetilde{\b{u}}_{\lambda} \cdot \nabla \widetilde{\varphi}_{\lambda}
		& =(\widetilde{\b{u}}_{\lambda,n} - \widetilde{\b{u}}_{\lambda})
		\cdot \nabla \widetilde{\varphi}_{\lambda,n}
		+ \widetilde{\b{u}}_{\lambda} \cdot (\nabla \widetilde{\varphi}_{\lambda,n}
		-\nabla \widetilde{\varphi}_{\lambda}).
	\end{split}
	\]
	Thus, by the H\"older inequality, it holds that
	\[
	\widetilde{\b{u}}_{\lambda,n} \cdot \nabla \widetilde{\varphi}_{\lambda,n}
	\rightharpoonup
	\widetilde{\b{u}}_{\lambda} \cdot \nabla \widetilde{\varphi}_{\lambda}
	\quad \text{in }L^{\frac{p}2}(\tom; L^1(0,T;L^\frac{3}{2}(\OO))\cap L^2(0,T;L^1(\OO))).
	\]
	\paragraph{\textit{The stochastic integrals}} Let us now identify $\widetilde{I}_{\lambda}$ and $\widetilde{J}_{\lambda}$. The procedure is standard, for instance see \cite[Section 8.4]{dapratozab}. We introduce a family of filtrations on $(\tom, \tF, \tP)$, namely we set
	\[
	\tF_{\lambda,n,t} := \sigma\left\{\widetilde{\bu}_{\lambda,n}(s),\, \widetilde{\varphi}_{\lambda,n}(s),\,  \widetilde{I}_{\lambda,n}(s),\,  \widetilde{J}_{\lambda,n}(s),\, \widetilde{W}_{1,n}(s),\, \widetilde{W}_{2,n}(s), \, \widetilde\bu_{0,n}, \, \widetilde{\varphi}_{0,n}, \, s \in [0,t]\right\},
	\]
	for any $t \geq 0$, $n \in \enne$ and $\lambda \in (0,1)$, in such a way that both $\widetilde{W}_{1,n}$ and $\widetilde{W}_{2,n}$ are adapted. In particular, by preservation of laws and the definitions of Wiener process and stochastic integral, we readily have that $W_{i,n}$ is a $Q_{i}^0$-Wiener process on $U_i^0$ and
	\[
	\widetilde{I}_{\lambda, n} = \int_0^t G_{1,n}(\widetilde{\bu}_{\lambda,n}(\tau)) \: \mathrm{d}\widetilde{W}_{1,n}(\tau), \qquad \widetilde{J}_{\lambda, n} = \int_0^t G_{2,\lambda,n}(\widetilde{\varphi}_{\lambda,n}(\tau)) \: \mathrm{d}\widetilde{W}_{2,n}(\tau),
	\]
	are respectively a $\bVsd$-valued and an $H$-valued martingale. Let us iterate the same procedure on the limit processes: we define
	\[
	\tF_{\lambda,t} := \sigma\left\{\widetilde{\bu}_{\lambda}(s),\,  \widetilde{\varphi}_{\lambda}(s),\, \widetilde{I}_{\lambda}(s),\,  \widetilde{J}_{\lambda}(s),\, \widetilde{W}_{1}(s),\, \widetilde{W}_{2}(s), \, \widetilde\bu_0,\,\widetilde{\varphi}_0,\, s \in [0,t]\right\}.
	\]
	It is easy to infer, by the proven convergences, that both $\widetilde{W}_{1}(0)$ and $\widetilde{W}_{2}(0)$ are zero. Let now $t > 0$, $s \in [0,t]$ and set
	\begin{align*}
		Z_{\bu,s} &:= L^2(0,s;\bHs) \cap C^0([0,s];D(\b{A}^{-\delta})), \\
		Z_{\varphi,s} &:= L^2(0,s;V_1) \cap C^0([0,s];H), \\
		\X_s &:= Z_{\bu,s} \times Z_{\varphi,s} \times C^0([0,s]; \bVsd) \times C^0([0,s]; H)  \times  C^0([0,s]; U_1^0) \times C^0([0,s]; U_2^0) \times \bVsd \times H.
	\end{align*}
	Let $\psi: \X_s \to \erre$ be a bounded and continuous function. By definition of martingale, we have
	\begin{equation} \label{eq:limit_n0}
		\tE \left[ \left(\widetilde{W}_{i,n}(t) -\widetilde{W}_{i,n}(s) \right) \psi\left( \widetilde{\varphi}_{\lambda,n},\, \widetilde{\bu}_{\lambda,n},\, \widetilde{I}_{\lambda,n},\,  \widetilde{J}_{\lambda,n},\, \widetilde{W}_{1,n},\, \widetilde{W}_{2,n},\, \widetilde\bu_{0,n},\,\widetilde{\varphi}_{0,n} \right)  \right] = 0
	\end{equation}
	for $i = 1,2$. Here, the arguments of $\psi$ are intended to be restricted over $[0,s]$ when necessary and $\tE$ denotes the expectation with respect to $\tP$. Letting $n \to +\infty$ in \eqref{eq:limit_n0}, an application of the dominated convergence theorem, owing to the proven convergences and the properties of $\psi$, entails
	\begin{equation} \label{eq:limit_n1}
		\tE \left[ \left(\widetilde{W}_{i}(t) -\widetilde{W}_{i}(s) \right) \psi\left( \widetilde{\varphi}_{\lambda,n},\, \widetilde{\bu}_{\lambda,n},\, \widetilde{I}_{\lambda},\,  \widetilde{J}_{\lambda},\, \widetilde{W}_{1},\, \widetilde{W}_{2},\, \widetilde\bu_0,\,\widetilde{\varphi}_0  \right)  \right] = 0,
	\end{equation}
	which expresses the fact that $\widetilde{W}_{i}$ is a $U_i^0$-valued $(\tF_{\lambda,t})_t$-martingale for $i = 1,2$. The characterization of $Q$-Wiener processes given in \cite[Theorem 4.6]{dapratozab} leads us to compute the quadratic variation of $\widetilde{W}_{i}$. To this end, notice that \eqref{eq:limit_n1} means that, for every $v, w \in U_i^0$
	\begin{multline*}
		\tE \left[ \left(\left( \widetilde{W}_{i,n}(t), v \right)_{U_i^0}
		\left( \widetilde{W}_{i,n}(t), w \right)_{U_i^0}
		-\left( \widetilde{W}_{i,n}(s), v \right)_{U_i^0}
		\left( \widetilde{W}_{i,n}(s), w \right)_{U_i^0}
		\right. \right. \\ \left. \left.
		- (t-s)\left(Q_i^0v,w\right)_{U_i^0}\right)
		\psi\left( \widetilde{\varphi}_{\lambda,n},\, \widetilde{\bu}_{\lambda,n},\,
		\widetilde{I}_{\lambda,n},\,  \widetilde{J}_{\lambda,n},\, \widetilde{W}_{1,n},\,
		\widetilde{W}_{2,n},\, \widetilde\bu_{0,n},\,\widetilde{\varphi}_{0,n} \right)  \right] = 0,
	\end{multline*}
	and using once more the dominated convergence theorem, we get
	\begin{multline*}
		\tE \left[ \left(\left( \widetilde{W}_{i}(t), v \right)_{U_i^0}\left( \widetilde{W}_{i}(t), w \right)_{U_i^0} -\left( \widetilde{W}_{i}(s), v \right)_{U_i^0}\left( \widetilde{W}_{i}(s), w \right)_{U_i^0} \right. \right. \\ \left. \left. - (t-s)\left(Q_i^0v,w\right)_{U_i^0}\right) \psi\left( \widetilde{\varphi}_{\lambda},\, \widetilde{\bu}_{\lambda},\, \widetilde{I}_{\lambda},\,  \widetilde{J}_{\lambda},\, \widetilde{W}_{1},\, \widetilde{W}_{2},\, \widetilde\bu_0,\,\widetilde{\varphi}_0 \right)  \right] = 0,
	\end{multline*}
	namely
	\begin{equation*}
		\left\llangle \widetilde{W}_{i}\right\rrangle(t) = tQ_i^0, \qquad t \in [0,T],
	\end{equation*}
	which is enough to conclude that $\widetilde{W}_{i}$ is a $Q_i^0$-Wiener process, adapted to $(\tF_{\lambda,t})_t$, owing to \cite[Theorem 4.6]{dapratozab}. We are now in a position to study the stochastic integrals. Arguing exactly as in \eqref{eq:limit_n0}-\eqref{eq:limit_n1}, we find that $\widetilde{I}_{\lambda}$ (resp. $\widetilde{J}_{\lambda}$) is a $\bVsd$-valued (resp. an $H$-valued) martingale. As far as the quadratic variations are concerned, an application of \cite[Theorem 4.27]{dapratozab} yields
	\begin{equation*}
		\left\llangle \widetilde{I}_{\lambda,n}\right\rrangle(t) = \int_0^t G_{1,n}(\widetilde{\bu}_{\lambda,n}(\tau)) \circ G_{1,n}(\widetilde{\bu}_{\lambda,n}(\tau))^*\: \mathrm{d}\tau , \qquad \left\llangle \widetilde{J}_{\lambda,n}\right\rrangle(t) = \int_0^t G_{2,\lambda,n}(\widetilde{\varphi}_{\lambda,n}(\tau)) \circ G_{2,\lambda,n}(\widetilde{\varphi}_{\lambda,n}(\tau))^*\: \mathrm{d}\tau,
	\end{equation*}
	for every $t \in [0,T]$. Let us outline the argument for the first sequence (similar considerations hold for the second one). Once again, fixing $\b{v}, \b{w} \in \bVs$, we have
	\begin{multline*}
		\tE \left[ \left(\left\langle \widetilde{I}_{\lambda,n}(t), \b{v} \right\rangle_{\bVsd, \bVs}\left\langle \widetilde{I}_{\lambda,n}(t), \b{w} \right\rangle_{\bVsd, \bVs} - \left\langle\widetilde{I}_{\lambda,n}(s), \b{v} \right\rangle_{\bVsd, \bVs}\left\langle \widetilde{I}_{\lambda,n}(s), \b{w} \right\rangle_{\bVsd, \bVs} \right. \right. \\ \left. \left. -\andrea{\int_0^t \left( G_{1}(\widetilde{\bu}_{\lambda}(\tau)) \circ G_{1}(\widetilde{\bu}_{\lambda}(\tau))^*\b{v}, \b{w}\right)_\bHs\d\tau}\right) \psi\left( \widetilde{\varphi}_{\lambda,n},\, \widetilde{\bu}_{\lambda,n},\, \widetilde{I}_{\lambda,n},\,  \widetilde{J}_{\lambda,n},\, \widetilde{W}_{1,n},\, \widetilde{W}_{2,n}
		,\, \widetilde\bu_{0,n},\,\widetilde{\varphi}_{0,n}\right)  \right] = 0,
	\end{multline*}
	and, as $n \to +\infty$, the dominated convergence theorem implies that
	\begin{multline*}
		\tE \left[ \left(\left\langle \widetilde{I}_{\lambda}(t), \b{v} \right\rangle_{\bVsd, \bVs}\left\langle \widetilde{I}_{\lambda}(t), \b{w} \right\rangle_{\bVsd, \bVs} - \left\langle\widetilde{I}_{\lambda}(s), \b{v} \right\rangle_{\bVsd, \bVs}\left\langle \widetilde{I}_{\lambda}(s), \b{w} \right\rangle_{\bVsd, \bVs} \right. \right. \\ \left. \left. -\int_0^t \left( G_{1}(\widetilde{\bu}_{\lambda}(\tau)) \circ G_{1}(\widetilde{\bu}_{\lambda}(\tau))^*\b{v}, \b{w}\right)_\bHs \mathrm{d}\tau\right) \psi\left( \widetilde{\varphi}_{\lambda},\, \widetilde{\bu}_{\lambda},\, \widetilde{I}_{\lambda},\,  \widetilde{J}_{\lambda},\, \widetilde{W}_{1},\, \widetilde{W}_{2},
		\, \widetilde\bu_{0},\,\widetilde{\varphi}_{0} \right)  \right] = 0,
	\end{multline*}
	Notice that in the above equality the dualities are necessary. The quadratic variation of $\widetilde{I}_{\lambda}$ is therefore
	\[
	\left\llangle \widetilde{I}_{\lambda} \right\rrangle(t) = \int_0^t G_{1}(\widetilde{\bu}_{\lambda}(\tau)) \circ G_{1}(\widetilde{\bu}_{\lambda}(\tau))^*\: \mathrm{d}\tau, \qquad t \in [0,T].
	\]
	Let us identify $\widetilde{I}_{\lambda}$ with the martingale
	\[
	\widetilde{M}_{\lambda}(t) := \int_0^t G_{1}(\widetilde{\bu}_{\lambda}(\tau))\: \mathrm{d}\widetilde{W}_1(\tau),
	\]
	which is a $\bVsd$-valued $(\tF_{\lambda,t})_t$-martingale having the same quadratic variation of $ \widetilde{I}_{\lambda}$. By \cite[Theorem 3.2]{Pard}, we can write
	\begin{equation} \label{eq:limit_n2}
		\begin{split}
			\left\llangle \widetilde{M}_{\lambda} - \widetilde{I}_{\lambda} \right\rrangle & = \left\llangle \widetilde{M}_{\lambda} \right\rrangle + \left\llangle \widetilde{I}_{\lambda} \right\rrangle -2 \left\llangle \widetilde{M}_{\lambda}, \widetilde{I}_{\lambda} \right\rrangle \\
			& = 2\int_0^\cdot G_{1}(\widetilde{\bu}_{\lambda}(\tau)) \circ G_{1}(\widetilde{\bu}_{\lambda}(\tau))^*\: \mathrm{d}\tau - 2\int_0^\cdot G_{1}(\widetilde{\bu}_{\lambda}(\tau)) \: \mathrm{d}\left\llangle\widetilde{W}_1, \widetilde{I}_{\lambda} \right\rrangle (\tau).
		\end{split}
	\end{equation}
	Thus, we now compute the cross quadratic variation appearing on the right hand side in \eqref{eq:limit_n2}. To this end, notice that by \cite[Theorem 3.2]{Pard}, we have
	\[
	\begin{split}
		\left\llangle\widetilde{I}_{\lambda,n}, \widetilde{W}_{1,n} \right\rrangle & = \int_0^\cdot G_{1,n}(\widetilde{\bu}_{\lambda,n}(\tau)) \circ \iota^{-1}_1 \: \mathrm{d}\left\llangle\widetilde{W}_{1,n}, \widetilde{W}_{1,n}\right\rrangle (\tau) \\
		& = \int_0^\cdot G_{1,n}(\widetilde{\bu}_{\lambda,n}(\tau)) \circ \iota^{-1}_1 \circ Q_1^0 \: \mathrm{d}\tau \\
		& = \int_0^\cdot G_{1,n}(\widetilde{\bu}_{\lambda,n}(\tau)) \circ \iota^{-1}_1 \circ \iota_1 \circ \iota^*_1 \: \mathrm{d}\tau \\
		& = \int_0^\cdot G_{1,n}(\widetilde{\bu}_{\lambda,n}(\tau)) \circ \iota^*_1 \: \mathrm{d}\tau,
	\end{split}
	\]
	where we also used the fact that $Q_1^0 = \iota_1 \circ \iota^*_1$, where $\iota_1: U_1 \to U_1^0$ is the classical Hilbert-Schmidt embedding. This implies that
	\[
	\left\llangle \widetilde{W}_{1,n}, \widetilde{I}_{\lambda,n} \right\rrangle = \int_0^\cdot \iota_1 \circ G_{1,n}(\widetilde{\bu}_{\lambda,n}(\tau))^* \: \mathrm{d}\tau.
	\]
	A further application of the dominated convergence theorem entails that, as $n \to +\infty$,
	\begin{equation} \label{eq:limit_n3}
		\left\llangle \widetilde{W}_{1}, \widetilde{I}_{\lambda} \right\rrangle = \int_0^\cdot \iota_1 \circ G_{1}(\widetilde{\bu}_{\lambda}(\tau))^* \: \mathrm{d}\tau.
	\end{equation}
	The identification follows injecting \eqref{eq:limit_n3} in \eqref{eq:limit_n2}.
	\paragraph{\textit{Identification of the limit solution.}} We are now left to prove that the limiting processes solve the regularized Allen-Cahn-Navier-Stokes system \eqref{eq:1_lam}-\eqref{eq:6_lam}. Testing \eqref{eq:1_disc} by some $\b{v} \in \bVs$ and integrating the obtained identity with respect to time yield
	\begin{multline*}
		(\bu_{\lambda,n}(t),\bv)_{\bHs} +
		\int_0^t\left[\ip{\b A \bu_{\lambda,n}(s)}{\bv}_{\b V_\sigma^*, \b V_\sigma}
		+\ip{\b B( \bu_{\lambda,n}(s),  \bu_{\lambda,n}(s))}{\bv}_{\b V_\sigma^*, \b V_\sigma}
		+\int_\OO \mu_{\lambda,n}(s)\nabla \varphi_{\lambda,n}(s) \cdot \bv
		\right]\,\d s \\ = (\widetilde\bu_{0,n},\bv)_{\bHs} +
		\left(\int_0^t G_{1,n}( \bu_{\lambda,n}(s))\,\d  W_1(s), \bv\right)_{\bHs}
		\qquad\forall\,t\in[0,T]\,,\quad \P\text{-a.s.}.
	\end{multline*}
	Letting $n\to+\infty$, owing to above convergences and using the dominated convergence theorem, we obtain
	\begin{multline} \label{eq:limit_n4}
		(\widetilde{\bu}_{\lambda}(t),\bv)_{\bHs} +
		\int_0^t\left[\ip{\b A \widetilde{\bu}_{\lambda}(s)}{\bv}_{\b V_\sigma^*, \b V_\sigma}
		+\ip{\b B( \widetilde{\bu}_{\lambda}(s),  \widetilde{\bu}_{\lambda}(s))}{\bv}_{\b V_\sigma^*, \b V_\sigma}
		+\int_\OO \mu_{\lambda}(s)\nabla \varphi_{\lambda}(s) \cdot \bv
		\right]\,\d s\\
		= (\widetilde\bu_{0},\bv)_{\bHs} +
		\left(\int_0^t G_{1}(\widetilde{\bu}_{\lambda}(s))\,\d  \widetilde{W}_1(s), \bv\right)_{\bHs}
		\qquad\forall\,t\in[0,T]\,,\quad \P\text{-a.s.}
	\end{multline}
	Next, we identify the limit chemical potential. Testing \eqref{eq:4_disc} by some $v \in V_1$, passing to the limit as $n \to +\infty$ yields and exploiting the proven convergences entail
	\begin{equation} \label{eq:weakmu}
		\int_\OO \widetilde{\mu}_\lambda v = -\int_\OO \Delta  \widetilde{\varphi}_\lambda v + \int_\OO F'_\lambda(\widetilde{\varphi}_\lambda)v,
	\end{equation}
	almost everywhere in $[0,T]$ and $\P$-almost surely. Finally, consider the approximating Allen-Cahn equation. Testing \eqref{eq:3_disc} by some $v \in V_1$ and
    passing to the limit as $n \to +\infty$, we get
	\[
	( \widetilde{\varphi}_{\lambda}(t), v)_H +
		\int_0^t\!\int_\OO\left[ \widetilde{\bu}_{\lambda}(s)\cdot\nabla \widetilde{\varphi}_{\lambda}(s) +  \widetilde{\mu}_{\lambda}(s)\right]v\,\d s= ( \widetilde{\varphi}_{0},v)_{H} +
		\left(\int_0^t G_{2,\lambda}( \widetilde{\varphi}_{\lambda}(s))\,\d  W_2(s), v\right)_{H}.
	\]
	Therefore, system \eqref{eq:1_lam}-\eqref{eq:6_lam} is satisfied (in the weak sense) once we identify (the law of) the initial state. By the properties of $X_n$, we know that
	\andrea{
	\[
	\widetilde{\varphi}_{0,n} \laweq {\varphi}_{0,n}, \qquad \widetilde{\bu}_{0,n} \laweq {\bu}_{0,n}
	\]
	}for any $n \in \enne$, and by uniqueness of the distributional limit (jointly with the above convergences) we conclude
	\andrea{
	\[
	\widetilde{\varphi}_{0} \laweq {\varphi}_{0}, \qquad \widetilde{\bu}_{0} \laweq {\bu}_{0}.
	\]
	}The initial conditions are therefore attained in law.
	\subsection{Uniform estimates with respect to $\lambda$} Here, we prove further uniform estimates, now independent of the Yosida parameter $\lambda$. The symbol $K$ (possibly numbered) denotes a positive constant, always independent of $\lambda$, which may change from line to line.
	\paragraph{\textit{First estimate.}} Notice that the constant $C_1$ in \eqref{eq:estimate1} does not depend on $\lambda$. By lower semicontinuity and preservation of laws of $X_n$, we infer
	\begin{equation} \label{eq:lam_est1}
		\|\widetilde{\varphi}_{\lambda}\|_{L^p_\cP(\Omega;C^0([0,T];H))} + \|\widetilde{\varphi}_{\lambda}\|_{L^{p}_\cP(\Omega;L^2(0,T;V_1))} \leq K_1.
	\end{equation}
	\paragraph{\textit{Second estimate.}} Let us collect, in \eqref{eq:unif23}, all controls which are already uniform with respect to $\lambda$, that is, the bounds on the diffusion coefficients \eqref{eq:unif24} and \eqref{eq:unif25}, the bounds on the stochastic terms \eqref{eq:unif27} and \eqref{eq:unif28}, and the initial data bounds given in \eqref{eq:unif29}. This can be summarized as follows
	(we can express the result in the new variables since $X_n$ preserves laws)
	\begin{align}
	\nonumber
	&\tE\supp\|\widetilde{\bu}_{\lambda,n}(t)\|^p_{\bHs}
	+\tE\supp\|\nabla \widetilde{\varphi}_{\lambda,n}(\tau)\|^p_{\b{H}}
	+\tE\supp \|F_\lambda(\widetilde{\varphi}_{\lambda,n})\|_{L^1(\OO)}^\frac{p}{2} \\
	\nonumber
	&\qquad+\tE\left|\int_0^t  \|\nabla\widetilde{\bu}_{\lambda,n}(\tau)\|^2_{\bHs}
	\,\mathrm{d}\tau \right|^\frac{p}{2}
	+\tE\left|\int_0^t  \|\widetilde{\mu}_{\lambda,n}(\tau)\|^2_{H}\,\mathrm{d}\tau \right|^\frac{p}{2} \\
	\nonumber
	&\leq C\left[\tE\|F_\lambda({\varphi}_{0,n})\|_{L^1(\OO)}^\frac{p}{2}+
	1 + \tE\left|\int_0^t\|\widetilde{\bu}_{\lambda,n}(\tau)\|^2_{ \bHs}\,\mathrm{d}\tau \right|^\frac{p}{2}
	+ \tE\left|\int_0^t\norm{\nabla \widetilde{\varphi}_{\lambda,n}(\tau)}_{\bH}^2
	\,\mathrm{d}\tau \right|^\frac{p}{2} \right. \\
	&\qquad\left. +\tE\left|\int_0^t\sum_{k=1}^\infty\int_\OO |F_\lambda''(\widetilde{\varphi}_{\lambda,n}(\tau))|
		|g_k(J_\lambda(\widetilde{\varphi}_{\lambda,n}(\tau)))|^2\,\mathrm{d}\tau \right|^\frac{p}{2}\right].
	\label{eq:unif41}
	\end{align}
	where $C$ depends on $p$ but is independent of $\lambda$.
	Next, we would like to take the limit as $n\to +\infty$ in \eqref{eq:unif41}. On the left hand side, the previously proven uniform estimates, convergences and weak lower semicontinuity of the norms are enough to pass to the limit. Moreover, it is easily seen, by Lipschitz-continuity of $F'_\lambda$, that $F_\lambda({\varphi}_{0,n}) \to F_\lambda({\varphi}_0)$ in $L^{\frac p2}(\tom; L^1(\OO))$ by the dominated convergence theorem. Finally, in order to pass to the limit in the last term at right hand side, we bound each term of the sequence as follows:
	\[
	\begin{split}
		\int_\OO |F_\lambda''(\widetilde{\varphi}_{\lambda,n}(\tau))|
		|g_k(J_\lambda(\widetilde{\varphi}_{\lambda,n}(\tau)))|^2 & = \int_\OO |F_\lambda''(\widetilde{\varphi}_{\lambda,n}(\tau))|
		|g_k(J_\lambda(\widetilde{\varphi}_{\lambda,n}(\tau)))|^2 \\
		& \leq |\OO|\|g_k \circ J_\lambda\|_{L^\infty(\mathbb{R})}^2\sup_{x \in \mathbb{R}} |F''_\lambda(x)| \\
		& = |\OO|\|g_k\|_{L^\infty(-1,1)}^2\sup_{x \in \mathbb{R}} |F''_\lambda(x)| .
	\end{split}
	\]
	Thanks to the proven convergences, it is straightforward to conclude that (cfr.~\cite{scarpa21})
	\[
	|F_\lambda''(\widetilde{\varphi}_{\lambda,n}(\tau))|
	|g_k(J_\lambda(\widetilde{\varphi}_{\lambda,n}(\tau)))|^2 \to |F_\lambda''(\widetilde{\varphi}_{\lambda}(\tau))|
	|g_k(J_\lambda(\widetilde{\varphi}_{\lambda}(\tau)))|^2
	\]
	almost everywhere in $\tom \times \OO \times (0,T)$.
	Therefore, applying the dominated convergence theorem
	and the weak lower semicontinuity of the norms, we find
	\begin{align}
	\nonumber
	&\tE\supp\|\widetilde{\bu}_{\lambda}(t)\|^p_{\bHs}
	+\tE\supp\|\nabla \widetilde{\varphi}_{\lambda}(\tau)\|^p_{\b{H}}
	+\tE\supp \|F_\lambda(\widetilde{\varphi}_{\lambda})\|_{L^1(\OO)}^\frac{p}{2} \\
	\nonumber
	&\qquad+\tE\left|\int_0^t  \|\nabla\widetilde{\bu}_{\lambda}(\tau)\|^2_{\bHs}
	\,\mathrm{d}\tau \right|^\frac{p}{2}
	+\tE\left|\int_0^t  \|\widetilde{\mu}_{\lambda}(\tau)\|^2_{H}\,\mathrm{d}\tau \right|^\frac{p}{2} \\
	\nonumber
	&\leq C\left[\tE\|F_\lambda({\varphi}_{0})\|_{L^1(\OO)}^\frac{p}{2}+
	1 + \tE\left|\int_0^t\|\widetilde{\bu}_{\lambda}(\tau)\|^2_{ \bHs}\,\mathrm{d}\tau \right|^\frac{p}{2}
	+ \tE\left|\int_0^t\norm{\nabla \widetilde{\varphi}_{\lambda}(\tau)}_{\bH}^2
	\,\mathrm{d}\tau \right|^\frac{p}{2} \right. \\
	&\qquad\left. +\tE\left|\int_0^t\sum_{k=1}^\infty\int_\OO |F_\lambda''(\widetilde{\varphi}_{\lambda}(\tau))|
		|g_k(J_\lambda(\widetilde{\varphi}_{\lambda}(\tau)))|^2\,\mathrm{d}\tau \right|^\frac{p}{2}\right].
	\label{eq:unif42}
	\end{align}
	We now need to find uniform bounds with respect to $\lambda$ for the two terms involving $F_\lambda$. Notice first that, as customary,
	\[
	\|F_\lambda({\varphi}_{0})\|_{L^1(\OO)} \leq \|F({\varphi}_{0})\|_{L^1(\OO)},
	\]
	which is finite by the hypotheses on the initial datum. Concerning the other term, we have
	\[ \small
	\begin{split}
		\int_0^t\sum_{k=1}^\infty\int_\OO |F_\lambda''(\widetilde{\varphi}_{\lambda}(\tau))|
		|g_k(J_\lambda(\widetilde{\varphi}_{\lambda}(\tau)))|^2\,\mathrm{d}\tau & = \int_0^t\sum_{k=1}^\infty\int_\OO |\beta_\lambda'(\widetilde{\varphi}_{\lambda}(\tau)) - c_F|
		|g_k(J_\lambda(\widetilde{\varphi}_{\lambda}(\tau)))|^2\,\mathrm{d}\tau \\
		& \leq c_FL_2^2|\OO|t + \int_0^t\sum_{k=1}^\infty\int_\OO |\beta'(J_\lambda(\widetilde{\varphi}_{\lambda}(\tau)))J'_\lambda(\widetilde{\varphi}_{\lambda}(\tau))|
		|g_k(J_\lambda(\widetilde{\varphi}_{\lambda}(\tau)))|^2\,\mathrm{d}\tau \\
		& \leq c_FL_2^2|\OO|t +\int_0^t\sum_{k=1}^\infty\int_\OO |F''(J_\lambda(\widetilde{\varphi}_{\lambda}(\tau))) + c_F|
		|g_k(J_\lambda(\widetilde{\varphi}_{\lambda}(\tau)))|^2\,\mathrm{d}\tau \\
		& \leq 2c_FL_2^2|\OO|t + \|F''g^2_k\|_{L^\infty(-1,1)}|\OO|t \\
		& \leq (2c_F + 1)L_2^2|\OO|t,
	\end{split}
	\]
	where we made use of \ref{hyp:G2} and we exploited the non-expansivity of $J_\lambda$. Collecting the two results in \eqref{eq:unif42}, we get
	\begin{align}
	\nonumber
	&\tE\supp\|\widetilde{\bu}_{\lambda}(t)\|^p_{\bHs}
	+\tE\supp\|\nabla \widetilde{\varphi}_{\lambda}(\tau)\|^p_{\b{H}}
	+\tE\supp \|F_\lambda(\widetilde{\varphi}_{\lambda})\|_{L^1(\OO)}^\frac{p}{2} \\
	\nonumber
	&\qquad+\tE\left|\int_0^t  \|\nabla\widetilde{\bu}_{\lambda}(\tau)\|^2_{\bHs}
	\,\mathrm{d}\tau \right|^\frac{p}{2}
	+\tE\left|\int_0^t  \|\widetilde{\mu}_{\lambda}(\tau)\|^2_{H}\,\mathrm{d}\tau \right|^\frac{p}{2} \\
	&\leq C_p\left[1 + \tE\int_0^t\|\widetilde{\bu}_{\lambda}(\tau)\|^p_{ \bHs}\,\mathrm{d}\tau
	+ \tE\int_0^t\norm{\nabla \widetilde{\varphi}_{\lambda}(\tau)}_{\bH}^p
	\,\mathrm{d}\tau \right],
	\label{eq:unif43}
	\end{align}
	and an application of the Gronwall lemma to \eqref{eq:unif43} gives
		\begin{align} \label{eq:lam_est2}
			\|\widetilde{\bu}_{\lambda}\|_{L^p_\cP(\tom;L^\infty([0,T];\bHs))}
			+ \|\widetilde{\bu}_{\lambda}\|_{L^p_\cP(\tom;L^2(0,T;V_1))} \leq K_2, \\
			\label{eq:lam_est3}
			\|\widetilde{\varphi}_{\lambda}\|_{L^p_\cP(\tom;L^\infty([0,T];V_1))} \leq K_3, \\
			\label{eq:lam_est4}
			\|\widetilde\mu_{\lambda,n}\|_{L^p_\cP(\tom;L^2(0,T;H))} +
			\|F_\lambda(\widetilde\varphi_{\lambda,n})\|_{L^{\frac p2}_\cP(\tom;C^0([0,T];L^1(\OO)))} \leq K_4.
		\end{align}
		\paragraph{\textit{Further estimates.}} Choosing $v = \beta_\lambda(\widetilde{\varphi}_\lambda)$ in \eqref{eq:weakmu} yields:
		\begin{equation*}
			\int_\OO \widetilde{\mu}_\lambda F'_\lambda(\widetilde{\varphi}_\lambda) +  c_F\int_\OO \left[ \widetilde{\mu}_\lambda\widetilde{\varphi}_\lambda - \widetilde{\varphi}_\lambda F'_\lambda(\widetilde{\varphi}_\lambda) \right] = \int_\OO \beta'_\lambda(\widetilde{\varphi}_\lambda) \nabla  \widetilde{\varphi}_\lambda \cdot \nabla\widetilde{\varphi}_\lambda  + \int_\OO |F'_\lambda(\widetilde{\varphi}_\lambda)|^2,
		\end{equation*}
		and exploiting the monotonicity of $\beta_\lambda$, the H\"{o}lder and the Young inequalities, after an integration over $[0,t]$, we get
		\[
		\dfrac{1}{2}\|F'_\lambda(\widetilde{\varphi}_\lambda)\|_{L^2(0,T;H)}^2 \leq  \dfrac{3}{2}\|\widetilde{\mu}_\lambda\|_{L^2(0,T;H)}^2 + \dfrac{3	c_F^2}{2}\|\widetilde{\varphi}_\lambda\|_{L^2(0,T;H)}^2.
		\]
		Therefore, by estimates \eqref{eq:lam_est1} and \eqref{eq:lam_est4}, we find
		\begin{equation} \label{eq:lam_est5}
			\|F'_\lambda(\widetilde{\varphi}_\lambda)\|_{L^p_\cP(\tom; L^2(0,T;H))} \leq K_5.
		\end{equation}
		Again, by comparison in \eqref{eq:4_lam}, we also obtain
		\begin{equation} \label{eq:lam_est6}
			\|\widetilde{\varphi}_\lambda\|_{L^p_\cP(\tom; L^2(0,T;V_2))} \leq K_6.
		\end{equation}
		The remaining estimates can be obtained following line by line the work already showed in Subsection \ref{ssec:est}. In this way, we also recover the following: given any $k \in (0,\frac{1}{2})$ and $p \geq 2$, there exist $\beta = \beta(p)$ and $\gamma = \gamma(p)$, satisfying $\beta p > 1$ and $\gamma p > 1$ if $p > 2$ (see Remarks \ref{rem:beta} and \ref{rem:gamma}), such that
		\begin{align}
			\label{eq:lam_est7}
			\left\| \int_0^\cdot G_{1}(\widetilde{\bu}_{\lambda}(\tau))\,\mathrm{d}W_1(\tau)
			\right\|_{L^p_\cP(\tom;W^{k,p}(0,T;\bHs))} & \leq K_{7}, \\
			\label{eq:lam_est8}
			\left\| \int_0^\cdot G_{2,\lambda}(\widetilde{\varphi}_{\lambda}(\tau))\,\mathrm{d}W_2(\tau)
			\right\|_{L^p_\cP(\tom;W^{k,p}(0,T;V_1))} & \leq K_{8}, \\
			\label{eq:lam_est9}
			\|\widetilde{\varphi}_{\lambda}\|_{L^p_\cP(\tom;W^{\beta,p}(0,T;V_1^*))} & \leq K_9 \\
			\label{eq:lam_est10}
			\|\widetilde{\bu}_{\lambda}\|_{L^{\frac p2}_\cP(\tom;W^{\gamma,p}(0,T;\bVsd))}
			& \leq K_{10}.
		\end{align}
	\subsection{Passage to the limit as $\lambda \to 0^+$} 	\label{ssec:lim_lam}
	We are now in a position to let $\lambda \to 0^+$ (along a suitable subsequence). The argument is similar to the one of Subsection \ref{ssec:lim_n}, thus we will omit some details for the sake of brevity. Iterating the proofs of Lemmas \ref{lem:tight1}-\ref{lem:tight3}, we learn that the family of laws of \[ (\widetilde{\bu}_{\lambda}, \widetilde{\varphi}_{\lambda},  G_{1}(\widetilde{\bu}_{\lambda}) \cdot \widetilde{W}_{1,\lambda}, G_{2,\lambda}(\widetilde{\varphi}_{\lambda}) \cdot \widetilde{W}_{2,\lambda}, \widetilde{W}_{1,\lambda}, \widetilde{W}_{2,\lambda}, \widetilde{\bu}_{0,\lambda}, \widetilde{\varphi}_{0,\lambda})_{\lambda \in (0,1)} \]
	is again tight in the product space
	\[
	Z_\b{u} \times Z_\varphi \times C^0([0,T]; \bVsd) \times   C^0([0,T]; H) \times C^0([0,T]; U_1^0) \times C^0([0,T]; U_2^0) \times \bVsd \times H.
	\]
	Here, we recall that $\widetilde{W}_{i,\lambda} \equiv \widetilde{W}_i$ and we set $\widetilde\bu_{0,\lambda} \equiv \widetilde{\bu}_0$ and $\widetilde{\varphi}_{0,\lambda} \equiv \widetilde{\varphi}_0$ for $i = 1,2$ and any $\lambda \in (0,1)$.
	Owing to the Prokhorov and Skorokhod theorems (see \cite[Theorem 2.7]{ike-wata} and \cite[Theorem 1.10.4, Addendum 1.10.5]{vaa-well}), there exists a probability space $(\hom, \hF,\hP)$ and a family of random variables $Y_\lambda:  (\hom, \hF)\to(\tom, \tF)$ such that the law of $Y_\lambda$ is $\tP$ for every $\lambda \in (0,1)$, namely $\hP \circ Y_\lambda^{-1} = \tP$ (so that composition with $Y_\lambda$ preserves laws), and the following convergences hold as $\lambda \to 0^+$:
	\begin{align*}
		\widehat{\bu}_{\lambda}:= \widetilde{\bu}_\lambda \circ Y_\lambda \to \widehat{\bu}
		\quad & \text{in } L^q(\hom;L^2(0,T;\bHs) \cap C^0([0,T]; D(\b{A}^{-\delta})))  \text{ if }q < p, \\
		\widehat{\bu}_{\lambda} \rightharpoonup \widehat{\bu}
		\quad & \text{in } L^p(\hom;L^2(0,T;\bVs)),\\
		\widehat{\bu}_{\lambda}  \overset{\ast}{\rightharpoonup} \widehat{\bu}
		 \quad & \text{in } L^p_w(\hom;L^\infty(0,T;\bHs)) \cap
		 L^{\frac p2}(\hom;W^{\gamma,p}(0,T;\bVsd)),\\
		\widehat{\varphi}_{\lambda}:= \widetilde{\varphi}_\lambda \circ Y_\lambda \to \widehat{\varphi}
		\quad & \text{in } L^q(\hom;L^2(0,T;V_1) \cap C^0([0,T];H)) \text{ if }q < p, \\
		\widehat{\varphi}_{\lambda} \rightharpoonup \widehat{\varphi}
		\quad & \text{in } L^p(\hom;L^2(0,T;V_2)),\\
		\widehat{\varphi}_{\lambda} \overset{\ast}{\rightharpoonup} \widehat{\varphi}
		\quad & \text{in } L^p_w(\hom;L^\infty(0,T;V_1)) \cap L^p(\hom; W^{\beta,p}(0,T;V_1^*)), \\
		\widehat{I}_{\lambda}:= (G_1(\widetilde{\varphi}_\lambda)
		\cdot \widetilde{W}_{1,\lambda}) \circ Y_\lambda \to \widehat{I}
		\quad & \text{in } L^q(\hom;C^0([0,T];\bVsd)) \text{ if }q < p, \\
		\widehat{J}_{\lambda}:= (G_{2,\lambda}(\widetilde{\varphi}_\lambda)
		\cdot \widetilde{W}_{1,\lambda}) \circ Y_\lambda\to \widehat{J}
		\quad & \text{in } L^q(\hom;C^0([0,T];H)) \text{ if }q < p, \\
		\widehat{W}_{1,\lambda} := \widetilde{W}_{1,\lambda} \circ Y_\lambda  \to \widehat{W}_{1}
		\quad & \text{in } L^q(\hom;C^0([0,T];U_1^0)) \text{ if }q < p, \\
		\widehat{W}_{2,\lambda} := \widetilde{W}_{2,\lambda} \circ Y_\lambda  \to \widehat{W}_{2}
		\quad & \text{in } L^q(\hom;C^0([0,T];U_2^0)) \text{ if }q < p, \\
		\widehat{\bu}_{0,\lambda} := \widetilde{\bu}_{0,\lambda} \circ Y_\lambda  \to \widehat{\bu}_0
		\quad & \text{in } L^q(\tom;\bVsd) \text{ if }q < p, \\
		\widehat{\varphi}_{0,\lambda} := \widetilde{\varphi}_{0,\lambda}
		\circ Y_\lambda\to \widehat{\varphi}_0  \quad & \text{in } L^q(\tom;H) \text{ if }q < p,
	\end{align*}
	for some limiting processes satisfying
	\begin{align*}
		\widehat{\bu} & \in L^{\frac p2}(\hom; W^{\gamma,p}(0,T;\bVsd))
		\cap
		L^p(\hom;C^0([0,T];D(\b{A}^{-\delta}))\cap L^2(0,T;\bVs)) \cap L^p_w(\hom;L^\infty(0,T;\bHs)); \\
		\widehat{\varphi} & \in L^p(\hom; W^{\beta,p}(0,T;V_1^*)\cap C^0([0,T];H) \cap L^2(0,T;V_2))
		\cap L^p_w(\hom;L^\infty(0,T;V_1)); \\
		\widehat{\mu} & \in L^p(\hom;L^2(0,T;H)); \\
		\widehat{I} & \in L^p(\hom;C^0([0,T];\bVsd)); \\
		\widehat{J} & \in L^p(\hom;C^0([0,T];H)); \\
		\widehat{W}_1 & \in L^p(\hom;C^0([0,T];U_1^0)); \\
		\widehat{W}_2 & \in L^p(\hom;C^0([0,T];U_2^0)); \\
		\widehat{\bu}_0 & \in L^p(\hom;\bHs);\\
		\widehat{\varphi}_0 & \in L^p(\hom;\mathcal{B} \cap V_1).
	\end{align*}
	Again, by estimate \eqref{eq:lam_est4}, we also have the following weak convergence of the redefined chemical potentials
	\[
	\widehat{\mu}_\lambda := \widetilde{\mu}_\lambda \circ Y_\lambda \rightharpoonup \widehat{\mu}
	\quad\text{in }L^p(\hom;L^2(0,T;H)).
	\]
	Mimicking the arguments illustrated in Subsection \ref{ssec:lim_n}, we now address several issues.
	\paragraph{\textit{The nonlinearities}.} First of all, we show that
	\[
		F'_\lambda(\widehat{\varphi}_\lambda) \to F'(\widehat{\varphi}) \text{ in }L^p(\hom;L^2(0,T;H)).
	\]
	This comes from the weak-strong closure of maximal monotone operators (see, for instance, \cite[Proposition 2.1]{barbu-monot}) combined with the strong convergence for $\widehat{\varphi}_\lambda$ proved above (recall that $F'_\lambda(x) = \beta_\lambda(x) + c_Fx)$. Next, the diffusion coefficients. As for $G_1$, it is easy by Lipschitz continuity to deduce
	\[
	G_{1}(\widehat{\bu}_{\lambda}) \to G_{1}(\widehat{\bu}) \quad \text{in }L^q(\hom; L^2(0,T;\cL^2(U_1,Y))) \text{ if } q < p.
	\]
	Moreover, arguing similarly (recall also Proposition \ref{prop:lipschitz}), we get
	\[
	\begin{split}
		& \|G_{2,\lambda}(\widehat{\varphi}_{\lambda})-G_{2}(\widehat{\varphi})\|_{L^p(\tom, L^2(0,T,\cL^2(U_2,H)))} \\
		& \qquad \qquad \leq \|G_{2,\lambda}(\widehat{\varphi}_{\lambda})-G_{2,\lambda}(\widehat{\varphi})\|_{L^p(\tom, L^2(0,T,\cL^2(U_2,H)))} + \|G_{2,\lambda}(\widehat{\varphi})-G_{2}(\widehat{\varphi})\|_{L^p(\tom, L^2(0,T,\cL^2(U_2,H)))},
	\end{split}
	\]
	and we conclude
	\[
	G_{2,\lambda}(\widehat{\varphi}_{\lambda}) \to G_{2}(\widehat{\varphi}) \quad \text{in }L^q(\hom; L^2(0,T;\cL^2(U_2,H))) \text{ if } q < p.
	\]
	Regarding the convective term and the Korteweg force, on account of the obtained convergences, we deduce that
	\begin{align*}
		\widehat\mu_{\lambda}\nabla\widehat\varphi_{\lambda} \rightharpoonup \widehat\mu\nabla\widehat\varphi &\quad\text{in }\b{L}^1(\OO \times (0,T)); \\
		\b{B}(\widehat\bu_\lambda, \widehat{\bu}_\lambda) \to \b{B}(\widehat\bu, \widehat{\bu})
		&\quad\text{in }\b{L}^q(\tom; L^\frac{4}{d}(0,T; \bVsd) \text{ if } q < \frac p2; \\
		\widehat{\bu}_\lambda \cdot \nabla\widehat{\varphi}_\lambda
		\rightharpoonup \widehat{\bu} \cdot \nabla\widehat{\varphi}
		&\quad\text{in }
		L^{\frac{p}{2}}(\hom; L^2(0,T;L^1(\OO))\cap L^1(0,T; L^{\frac32}(\OO))).
	\end{align*}
	
	\paragraph{\textit{The stochastic integrals.}} Following line by line the argument presented in Subsection \ref{ssec:lim_n}, it is possible to identify the limits $\widehat{I}$ and $\widehat{J}$. Indeed, we have
	\[
	\widehat{I}(t) = \int_0^t G_{1}(\widehat{\bu}(\tau)) \: \mathrm{d}\widehat{W}_{1}(\tau), \qquad \widehat{J}(t) = \int_0^t G_{2}(\widehat{\varphi}(\tau)) \: \mathrm{d}\widehat{W}_{2}(\tau),
	\]
	which are a $\bVsd$ and an $H$-valued martingale, respectively, adapted with respect to a suitable filtration $(\hF_t)_t$.
	\paragraph{\textit{Identification of the limit solution.}} Again, a multiple application of the dominated convergence theorem allows us to
infer that the limit processes form a martingale solution of the original problem. The existence of a martingale solution is proved.
	\subsection{The energy inequality} We are left to prove the energy inequality. To this end, we simply pass to the limit in a suitable approximating energy inequality. Let us add \eqref{eq:unif20} and \eqref{eq:unif21} together and take expectations. Recalling that stochastic integrals are martingales, we obtain the identity
	\begin{align}
	\nonumber
		&\dfrac{1}{2}\E\|\b{u}_{\lambda,n}(t)\|^2_{\bHs}
		+ \dfrac{1}{2}\E\|\nabla \varphi_{\lambda,n}\|_\b{H}^2
		+ \E\|F_\lambda(\varphi_{\lambda,n})\|_{L^1(\OO)}
		+ \E\int_0^t \left[ \|\nabla\b{u}_{\lambda,n}(\tau)\|^2_{\bHs}
		+ \|\mu_{\lambda,n}(\tau)\|^2_H\right]\mathrm{d}\tau \\
		\nonumber
		&= \dfrac{1}{2}\E\|\b{u}_{0,n}\|^2_{\bHs}
		+ \dfrac{1}{2}\E\|\nabla \varphi_{0,n}\|_\b{H}^2
		+ \E\|F_\lambda(\varphi_{0,n})\|_{L^1(\OO)}
		+ \dfrac{1}{2}\E\int_0^t \|G_{1,n}(\b{u}_{\lambda,n}(\tau))\|^2_{\cL^2(U_1, \bHs)} \: \mathrm{d}\tau \\
		&\qquad+\frac12 \E\int_0^t
		\left[ \norm{\nabla G_{2,\lambda,n}(\varphi_{\lambda,n}(\tau))}_{\cL^2(U_2,\bH)}^2
		+ \sum_{k=1}^\infty\int_\OO F_\lambda''(\varphi_{\lambda,n}(\tau))
		|g_k(J_\lambda(\varphi_{\lambda,n}(\tau)))|^2\right]\, \mathrm{d}\tau.
		\label{eq:energyineq1}
	\end{align}
	Thank to \eqref{eq:unif24} and \eqref{eq:unif25}, from \eqref{eq:energyineq1} we infer
	\begin{align}
	\nonumber
		&\dfrac{1}{2}\E\|\b{u}_{\lambda,n}(t)\|^2_{\bHs}
		+ \dfrac{1}{2}\E\|\nabla \varphi_{\lambda,n}\|_\b{H}^2
		+ \E\|F_\lambda(\varphi_{\lambda,n})\|_{L^1(\OO)}
		+ \E\int_0^t \left[ \|\nabla\b{u}_{\lambda,n}(\tau)\|^2_{\bHs}
		+ \|\mu_{\lambda,n}(\tau)\|^2_H\right]\mathrm{d}\tau \\
	\nonumber
		&\leq C_{G_1}^2t
		+ \dfrac{1}{2}\E\|\b{u}_{0,n}\|^2_{\bHs}
		+ \dfrac{1}{2}\E\|\nabla \varphi_{0,n}\|_\b{H}^2
		+ \E\|F_\lambda(\varphi_{0,n})\|_{L^1(\OO)}
		+ C_{G_1}^2\E\int_0^t \|\b{u}_{\lambda,n}(\tau)\|^2_{\bHs} \: \mathrm{d}\tau  \\
		&\qquad
		+ \frac{L_2^2}2\E\int_0^t \norm{\nabla \varphi_{\lambda,n}(\tau)}_{\bH}^2 \, \mathrm{d}\tau
		+ \frac12\E\int_0^t\sum_{k=1}^\infty\int_\OO F_\lambda''(\varphi_{\lambda,n}(\tau))
		|g_k(J_\lambda(\varphi_{\lambda,n}(\tau)))|^2\, \mathrm{d}\tau.
		\label{eq:energyineq2}
	\end{align}
	Exploiting the preservation of laws by $X_n$, and letting $n \to +\infty$, we find
	\begin{align}
	\nonumber
	&\dfrac{1}{2}\tE\|\widetilde\bu_{\lambda}(t)\|^2_{\bHs}
	+ \dfrac{1}{2}\tE\|\nabla \widetilde\varphi_{\lambda}\|_\b{H}^2
	+ \tE\|F_\lambda(\widetilde\varphi_{\lambda})\|_{L^1(\OO)}
	+ \tE\int_0^t \left[ \|\nabla\widetilde\bu_{\lambda}(\tau)\|^2_{\bHs}
	+ \|\widetilde\mu_{\lambda}(\tau)\|^2_H\right]\mathrm{d}\tau \\
	\nonumber
	&\leq C_{G_1}^2 t
	+ \dfrac{1}{2}\tE\|\widetilde\bu_{0}\|^2_{\bHs}
	+ \dfrac{1}{2}\tE\|\nabla \widetilde\varphi_{0}\|_\b{H}^2
	+ \tE\|F_\lambda(\widetilde\varphi_{0})\|_{L^1(\OO)}
	+ C_{G_1}^2\tE\int_0^t \|\widetilde\bu_{\lambda}(\tau)\|^2_{\bHs} \: \mathrm{d}\tau \\
	&\qquad
	+\frac{L_2^2}2\tE\int_0^t \norm{\nabla \widetilde\varphi_{\lambda}(\tau)}_{\bH}^2 \, \mathrm{d}\tau
	+ \frac12\tE\int_0^t \sum_{k=1}^\infty\int_\OO F_\lambda''(\widetilde\varphi_{\lambda}(\tau))
		|g_k(J_\lambda(\widetilde\varphi_{\lambda}(\tau)))|^2\, \mathrm{d}\tau.
	\label{eq:energyineq3}
	\end{align}
	Here we have used the lower semicontinuity of the norms and the dominated convergence theorem. A second passage to the limit entails the claimed inequality. Indeed, exploiting preservation of laws by $Y_\lambda$ in \eqref{eq:energyineq3} as well as \ref{hyp:G2}, and letting $\lambda \to 0^+$, we get
	\begin{align}
	\nonumber
	&\dfrac{1}{2}\hE\|\widehat\bu(t)\|^2_{\bHs}
	+ \dfrac{1}{2}\hE\|\nabla \widehat\varphi\|_\b{H}^2
	+ \hE\|F(\widehat\varphi)\|_{L^1(\OO)}
	+ \hE\int_0^t \left[ \|\nabla\widehat\bu(\tau)\|^2_{\bHs}
	+ \|\widehat\mu(\tau)\|^2_H\right]\mathrm{d}\tau \\
	\nonumber
	&\leq \left(C_{G_1}^2 + \frac{L_2^2}2|\OO|\right)t
	+ \dfrac{1}{2}\hE\|\widehat\bu_{0}\|^2_{\bHs}
	+ \dfrac{1}{2}\hE\|\nabla \widehat\varphi_{0}\|_\b{H}^2
	+ \hE\|F(\widehat\varphi_{0})\|_{L^1(\OO)} \\
	&\qquad+ C_{G_1}^2\hE\int_0^t \|\widehat\bu(\tau)\|^2_{\bHs} \: \mathrm{d}\tau
	+ \frac{L_2^2}2\hE\int_0^t \norm{\nabla \widehat\varphi(\tau)}_{\bH}^2 \, \mathrm{d}\tau.
	\label{eq:energyineq4}
	\end{align}
	Observe that, passing in the limit in the third term on the left hand side of \eqref{eq:energyineq3} is possible by lower semicontinuity since recalling that
	\[
	|J_\lambda \hphi_\lambda - \hphi| \leq |J_\lambda \hphi_\lambda - \hphi_\lambda| + |\hphi_\lambda - \hphi| \leq \lambda|\beta_\lambda(\hphi_\lambda)| + |\hphi_\lambda - \hphi|,
	\]
	it follows $J_\lambda \hphi_\lambda \to \hphi$ almost everywhere in $\hom \times \OO \times (0,T)$. \andrea{Fixed any $t > 0$}, the energy inequality follows taking the supremum \andrea{over $[0,t]$} in both sides of \eqref{eq:energyineq4}.
	\subsection{Recovery of the pressure.} \label{ssec:pressure} It is possible to recover a pressure through a generalization of the classical De Rham theorem to stochastic processes (see \cite{Simon03}). The result is of independent interest and we report it hereafter for reader's convenience.
	\begin{thm}[{\cite[Theorem 4.1]{Simon03}}] \label{th:derham}
		Let $\OO$ be a bounded Lipschitz domain of $\mathbb{R}^d$ and let $(\Omega, \cF, \P)$ be a complete probability space. Let $s_1 \in \andrea{\mathbb{R}}$ and $r_0,\, r_1 \in [1,+\infty]$. Let
		\[
		\b{h} \in L^{r_0}(\Omega; W^{s_1,r_1}(0,T;\andrea{(\b{H}^1_0(\OO))^*}))
		\]
		be such that, for all $\b{v} \in \left[C_0^\infty(\OO) \right]^d$ satisfying $\div \b{v} = 0$,
		\[
		\left\langle \b{h}, \b{v}  \right\rangle_{\left( \left[ C_0^\infty(\OO) \right]^d\right)^*,
		\left[ C_0^\infty(\OO) \right]^d} = 0
		\qquad \text{in }\left(C^\infty_0(0,T)\right)^*,\quad\P\text{-a.s.}
		\]
		Then there exists a unique (up to a constant)
		\[
		\pi \in L^{r_0}(\Omega; W^{s_1,r_1}(0,T;H))
		\]
		such that
		\[
		\nabla \pi = \b{h} \qquad \text{in }
		\left(\left[\mathcal{C}^\infty_0((0,T)\times\OO)\right]^d\right)^*, \quad\P\text{-a.s.}
		\]
		and
		\[
		\int_\OO \pi = 0 \qquad \text{in }\left(\mathcal{C}^\infty_0(0,T)\right)^*, \quad\P\text{-a.s.}
		\]
		Furthermore, there exists a positive constant $C = C(\OO)$, independent of $\b{h}$, such that
		\[
		\|\pi\|_{W^{s_1,r_1}(0,T;H)} \leq C(\OO)\|\b{h}\|_{W^{s_1,r_1}(0,T;\andrea{(\b{H}^1_0(\OO))^*})}
		\quad \P\text{-a.s.}
		\]
	\end{thm} \noindent
	Let us now find suitable values for the parameters $r_0, r_1$ and $s_1$.
	By choosing
	$\bv \in \left[ \mathcal{C}_0^\infty(\OO) \right]^d$ with $\div \b{v} = 0$ in \eqref{eq:ACNS1},
	after elementary rearrangements and integration by parts we obtain that
	\begin{align*}
	&\left\langle \partial_t(\hbu - G_1(\hbu)\cdot \widehat{W}_1)(t), \b{v}  \right\rangle_\andrea{(\b{H}^1_0(\OO))^*, \b{H}^1_0(\OO)}
	+\int_\OO\nabla\hbu(t): \nabla\b v
	\\
	& \qquad \qquad + \left\langle \b B(\hbu(t), \hbu(t)), \b{v}  \right\rangle_\andrea{(\b{H}^1_0(\OO))^*, \b{H}^1_0(\OO)} -\int_\OO\hmu(t)\nabla\hphi(t)\cdot\b v = 0
	\end{align*}
	for almost every $t\in(0,T)$, $\widehat\P$-almost surely.
	Hence,
	by setting
	\[
	\b{h} := \partial_t(\hbu - G_1(\hbu)\cdot\widehat{W}_1)
	+\b L\hbu + \b B(\hbu, \hbu) - \hmu\nabla\hphi,
	\]
	one has in particular,
	for all $\b{v} \in \left[ \mathcal{C}_0^\infty(\OO) \right]^d$ with $\div \b{v} = 0$, that
	\[
		\left\langle \b{h}, \b{v}  \right\rangle_{\left( \left[ \mathcal{C}_0^\infty(\OO) \right]^d\right)^*,
		\left[ \mathcal{C}_0^\infty(\OO) \right]^d} = 0
		\qquad \text{in }\left(\mathcal{C}^\infty_0(0,T)\right)^*,\quad \andrea{\hP\text{-a.s.}}
	\]
	Let us recover the regularity of $\b h$. Observing that
	$\hbu - G_1(\hbu)\cdot \widehat{W}_1 \in L^p_\cP(\Omega; L^\infty(0,T; \bH))$
	and that $\partial_t : L^\infty(0,T;\bH) \to W^{-1,\infty}(0,T;\bH)$ is
	linear and continuous, we have
	\[
	\partial_t(\hbu - G_1(\hbu)\cdot\widehat{W}_1) \in
	L^p_{\cP}(\hom; W^{-1,\infty}(0,T;\bHs)) \subset L^p_\cP(\hom; W^{-1,\infty}(0,T;\andrea{(\b{H}^1_0(\OO))^*})).
	\]
	Furthermore, recalling that
	$L^1(0,T;\andrea{(\b{H}^1_0(\OO))^*}) \embed W^{-1,\infty}(0,T;\andrea{(\b{H}^1_0(\OO))^*})$
	thanks to the fundamental theorem of calculus as shown in the proof of \cite[Theorem 2.2]{Simon03},
	one has that
	\[
	\b L\hbu \in  L^p_\cP(\hom; L^2(0,T;\andrea{(\b{H}^1_0(\OO))^*}))\subset L^p_\cP(\hom; W^{-1,\infty}(0,T;\andrea{(\b{H}^1_0(\OO))^*})).
	\]
	Moreover, since for $d\in\{2,3\}$ the bilinear form
	\[
	\b{B}: \bVs \times \bVs \to \b{L}^\frac{6}{5}(\OO) \embed \andrea{(\b{H}^1_0(\OO))^*}
	\]
	is continuous, thanks to the regularity of $\hbu$ it follows that
	\[
	  \b B(\hbu,\hbu) \in
	  L^{\frac{p}{2}}_\cP(\hom; L^1(0,T;\andrea{(\b{H}^1_0(\OO))^*})) \embed L^{\frac{p}2}_\cP(\hom; W^{-1,\infty}(0,T;\andrea{(\b{H}^1_0(\OO))^*})).
	\]
	Eventually, iterating the computations in \eqref{eq:kortewegest}, we obtain
	\[
	\hmu\nabla\hphi \in L^{\frac p2}_\cP(\hom;L^\frac{4}{3}(0,T;\andrea{(\b{H}^1_0(\OO))^*})) \embed
	L^{\frac p2}_\cP(\hom; W^{-1,\infty}(0,T;\andrea{(\b{H}^1_0(\OO))^*})),
	\]
	Hence, we have shown that $\b{h} \in L^{\frac p2}(\hom; W^{-1,\infty}(0,T;\andrea{(\b{H}^1_0(\OO))^*}))$ and an application of Theorem \ref{th:derham} with $r_0 = \andrea{\frac p2}$, $s_1 = -1$ and $r_1 = +\infty$
	yields the existence of the (unique up to a constant) pressure
	$\pi \in L^{\frac p2}(\hom; W^{-1,\infty}(0,T;H))$.
	Finally, we derive an estimate for $\pi$. The continuous dependence given by
	Theorem~\ref{th:derham} implies that
	\[
	\begin{split}
		&\|\widehat{\pi}\|_{W^{-1,\infty}(0,T;H)} \\
		& \leq C\left(\|\hbu-G_1(\hbu)\cdot\widehat{W}_1\|_{L^\infty(0,T;\bHs)}
		+ \|\hbu\|_{L^2(0,T;\bVs)} + \|\hbu\|_{L^2(0,T;\bVs)}^2
		+ \|\hmu\nabla\hphi\|_{L^\frac{4}{3}(0,T;\andrea{(\b{H}^1_0(\OO))^*})}\right).
	\end{split}
	\]
	Knowing that
	\[
		\begin{split}
			\|\hmu\nabla\hphi\|_{L^\frac{4}{3}(0,T;\b{V}_1^*)} & \leq C
			\left(\|\hmu\|_{L^2(0,T;H)}^2 + \|\hphi\|_{L^2(0,T;V_2)}^2 \right) \\
			& \leq C\left(\|F'(\hphi)\|_{L^2(0,T;H)}^2 + 2\|\hphi\|_{L^2(0,T;V_2)}^2 \right),
		\end{split}	
	\]
	and exploiting the Burkholder-Davis-Gundy inequality together with assumption \ref{hyp:G1}, we arrive at
	\begin{align*}
		\|\widehat{\pi}\|_{L^{\frac p2}_\cP(\hom;W^{-1,\infty}(0,T;H))} &\leq
		C\left(1+\|\hbu\|_{L^{\frac p2}_\cP(\hom;L^\infty(0,T;\bHs))}
		+ \|\hbu\|_{L^{\frac p2}(\hom;L^2(0,T;\bVs))}
		+ \|\hbu\|_{L^p(\hom;L^2(0,T;\bVs))}^2 \right. \\
		&\qquad\quad\left. + \|\hphi\|_{L^p(\hom;L^2(0,T;V_2))}^2
		+ \|F'(\hphi)\|_{L^p(\hom;L^2(0,T;H))}^2\right).
	\end{align*}
	The proof of Theorem \ref{th:martingalesols} is complete.
	
	\section{Existence of probabilistically-strong solutions when $d = 2$} \label{sec:proof2} \noindent
	This section is devoted to proving Theorem \ref{th:probstrongsol}. To this end, we will use a standard approach, namely we shall deduce it
    from pathwise uniqueness of martingale solutions.
	\begin{prop} \label{prop:uniqueness}
		Let $d = 2$ and $p \in (2,+\infty)$. Assume \ref{hyp:potential}-\ref{hyp:G2} and
		consider two sets of initial conditions $(\bu_{0,i}, \varphi_{0,i})$ for $i = 1,2$
		complying with the hypotheses of Theorem \ref{th:martingalesols}.
		Let $(\widehat\varphi_i, \widehat\bu_i)$ denote some martingale solutions
		to \eqref{eq:ACNS1}-\eqref{eq:ACNS6}, defined
		on the same suitable filtered space $(\hom, \hF, (\hF_t)_t, \hP)$
		and with respect to a pair of Wiener processes $\widehat W_1, \widehat W_2$. Then,
		there exist a sequence of positive real numbers $(C_n)_n$ and
		a sequence of stopping times $\{\zeta_n\}_n$, with
		$\zeta_n\nearrow T$ $\widehat\P$-almost surely as $n\to\infty$, such that
		the following continuous dependence estimate holds
		\begin{align*}
			&\|(\hbu_1 - \hbu_2)^{\zeta_n}\|_{L^p_\cP(\hom;
			C^0([0,T];\bVsd)) \cap L^p_\cP(\hom; L^2(0,T;H))}
			+ \|(\hphi_1 - \hphi_2)^{\zeta_n}\|_{L^p_\cP(\hom;
			C^0([0,T];H)) \cap L^p_\cP(\hom; L^2(0,T;V_1))}\\
			&\qquad\leq C_n\left(\|\hbu_{0,1}-\hbu_{0,2}\|_{L^p(\hom; \bVsd)}
			+ \|\hphi_{0,1}-\hphi_{0,2}\|_{L^p(\hom; H)} \right).
		\end{align*}
		In particular, the martingale solution to \eqref{eq:ACNS1}-\eqref{eq:ACNS6} is pathwise unique.
	\end{prop}
	\begin{proof}
		Let us set
		\begin{align*}
			\hbu & := \hbu_1 - \hbu_2, \\
			\hphi & := \hphi_1 - \hphi_2, \\
			\hmu & := \hmu_1 - \hmu_2, \\
			\hbu_{0} & := \bu_{0,1}-\bu_{0,2}, \\
			\hphi_{0} & := \varphi_{0,1}-\varphi_{0,2}.
		\end{align*}
		For every $n\in\enne$ and $i\in\{1,2\}$ we define the stopping time
		$\zeta_n^i:\hom\to\erre$ as
		\[
		\zeta_n^i:=\inf\left\{t\in[0,T]: \sup_{s\in[0,t]}\norm{\hbu_i(s)}^2_\bHs
		+ \int_0^t\left(\norm{\hbu_i(s)}_{\bVs}^2
		+\norm{\widehat\varphi_i(s)}_{V_2}^2\right)\,\mathrm{d} s \geq n^2\right\},
		\]
		with the usual convention that $\inf\emptyset=T$, and set
		\[
		\zeta_n:=\zeta_n^1\wedge\zeta_n^2.
		\]
		Clearly, $\zeta_n\nearrow T$ almost surely as $n\to\infty$.
		Let us also introduce the functionals
		\begin{align*}
			\Psi_1: \bVsd \to \mathbb{R}, & \qquad \Psi_1(\bv) := \dfrac{1}{2}\|\nabla \ia\bv\|_\bHs^2,\\
			\Psi_2: V_1 \to \mathbb{R}, & \qquad \Psi_2(v) :=  \dfrac{1}{2}\|\nabla v\|_H^2.
		\end{align*}
		We point out, once and for all, that what follows is valid $\hP$-almost
		surely for every $t \in [0,T]$. Let us consider at first $\Psi_1$. First of all,
		let us compute its first two Fréchet derivatives. If we set
		\[
		\Psi_0: \bVs \to \erre, \qquad \Psi_0(\bv) := \dfrac{1}{2}\|\nabla \bv\|_\bHs^2,
		\]
		then we have $\Psi_1 = \Psi_0 \circ \ia$. Therefore, an application
		of the chain rule implies that $D\Psi_1: \bVsd \to \b{V}_\sigma^{**}$ is defined by
		\[
		\begin{split}
			D\Psi_1(\b{v}) &= D(\Psi_0 \circ \ia)(\bv) \\
			& = D\Psi_0(\ia\bv) \circ D\ia(\bv) \\
			& = \b{A}\ia\bv \circ \ia \\
			& = \bv \circ \ia.
		\end{split}
		\]
		Here, of course, we exploited the facts that $D\Psi_0 = \b{A}$
		and that $\ia \in \mathcal{L}(\bVsd, \bVs)$. The above identity must be understood as follows
		\[
		\left\langle 	D\Psi_1(\b{v}), \b{w} \right\rangle_{\b{V}_\sigma^{**}, \bVsd} = \left\langle \bv, \ia\b{w}\right\rangle_{\bVsd,\bVs} = \left(\bv, \ia \b{w}\right)_\bHs.
		\]
		Moreover, by the properties of the inverse of the Stokes operator, it holds
		\begin{equation} \label{eq:inversestokes}
		\left\langle \bv, \ia\b{w}\right\rangle_{\bVsd,\bVs} = \left( \nabla \ia\bv, \nabla\ia\b{w} \right)_\bHs = \left\langle \b{w}, \ia\b{v}\right\rangle_{\bVsd,\bVs}
		\end{equation}
		for every $\bv, \b{w} \in \bVsd$.
		Notice that $D\Psi_1 \in \mathcal{L}(\bVsd, \b{V}_\sigma^{**})$ and thus $D^2\Psi_1(\bv) = D\Psi_1$	 for every $\bv \in \bVsd$.
		Applying the It\^{o} lemma \cite[Theorem 4.32]{dapratozab} to $\Psi_1(\hbu)$ and stopping at time $\zeta_n$, we obtain
		\begin{align}
		\nonumber
		&\dfrac{1}{2}\|\nabla \ia \hbu(t\wedge\zeta_n)\|_\bHs^2
		+ \int_0^{t\wedge\zeta_n} \left[ \left\langle\hbu(\tau), \ia\left[ \b{B}(\hbu_1(\tau),
		\hbu_1(\tau)) - \b{B}(\hbu_2(\tau), \hbu_2(\tau))
		\right]\right\rangle_{\bVsd, \bVs} \right]\mathrm{d}\tau \\
		\nonumber
		&\qquad+ \int_0^{t\wedge\zeta_n} \left[\| \hbu(\tau)\|_\bHs^2 -\left\langle \hbu(\tau),
		\ia\left[ \hmu_1(\tau)\nabla\hphi_1(\tau) -  \hmu_2(\tau)\nabla \hphi_2(\tau)
		\right]\right\rangle_{\bVsd, \bVs}\right]\mathrm{d}\tau \\
		\nonumber
		&= \dfrac{1}{2}\|\nabla \ia\hbu_0\|^2_{\bHs}
		+\int_0^{t\wedge\zeta_n}
		\left\langle \hbu(\tau), \ia\left[\big( G_{1}( \hbu_1(\tau)) - G_1(\hbu_2(\tau))
		\big) \d  \widehat W_1(\tau)\right]\right\rangle_{\bVsd, \bVs}\\
		&\qquad+ \dfrac{1}{2}\int_0^{t\wedge\zeta_n} \|\ia G_{1}(\hbu_1(\tau))
		- \ia G_1(\hbu_2(\tau))\|^2_{\cL^2(U_1, \bHs)} \: \mathrm{d}\tau.
		\label{eq:uniq1}
		\end{align}
		For the ease of notation, throughout computations
		we may omit the evaluation of the functions at the time $\tau\in[0,\zeta_n(\omega)]$,
		for $\widehat\P$-almost every $\omega\in\widehat\Omega$.
		We address the various terms in \eqref{eq:uniq1} separately.
		First of all, notice that, by \eqref{eq:inversestokes},
		\[
		\begin{split}
			\left\langle\hbu, \ia\left[ \b{B}(\hbu_1, \hbu_1) - \b{B}(\hbu_2, \hbu_2)\right]\right\rangle_{\bVsd, \bVs} & =
			\left\langle \b{B}(\hbu_1, \hbu_1), \ia\hbu\right\rangle_{\bVsd, \bVs} - \left\langle\b{B}(\hbu_2, \hbu_2), \ia\hbu\right\rangle_{\bVsd, \bVs} \\
			& = (\bu \otimes \bu_1, \nabla \ia \bu)_\bHs +  (\bu_2 \otimes \bu, \nabla \ia \bu)_\bHs,
		\end{split}
		\]
		on account of the incompressibility condition
		\[
		(\bu_i \cdot \nabla)\bu_i = -\div(\bu_i \otimes \bu_i)
		\]
		for $i = 1,2$.
		Then, using the H\"{o}lder, Young and Ladyzhenskaya inequalities, together with the definition of $\zeta_n$, we find
		\begin{align}
		\nonumber
				& \left| (\bu \otimes \bu_1, \nabla \ia \bu)_\bHs
				+  (\bu_2 \otimes \bu, \nabla \ia \bu)_\bHs \right|\\
		\nonumber
				& \qquad \qquad \leq \left( \|\hbu_1\|_{\b{L}^4(\OO)}
				+ \|\hbu_2\|_{\b{L}^4(\OO)} \right)\|\hbu\|_\bHs\|\nabla \ia \bu\|_{\b{L}^4(\OO)}\\
		\nonumber
				& \qquad \qquad \leq C
				\left( \|\hbu_1\|_{\bHs}^{\frac12}+ \|\hbu_2\|_\bHs^{\frac12} \right)
				\left( \|\hbu_1\|_{\bVs}^{\frac12}+ \|\hbu_2\|_\bVs^{\frac12} \right)
				\|\hbu\|_\bHs^\frac{3}{2}\|\nabla \ia \bu\|^\frac{1}{2}_\bHs\\
				& \qquad \qquad \leq \dfrac{1}{6}\|\hbu\|_\bHs^2 +
				Cn^2
				\left( \|\hbu_1\|_{\bVs}^2 + \|\hbu_2\|_\bVs^2 \right)\|\nabla\ia\hbu\|_\bHs^2.
		 \label{eq:uniq11}
		\end{align}
		Here, we also used the well-known fact that $\|\b{A}\bu\|_\bHs$ is an equivalent norm in $\b{H}^2(\OO) \cap \bVs$. Next, we address the coupling term. We make use of the customary formula
		\[
		\hmu_i\nabla\hphi_i = -\div(\nabla \hphi_i \otimes \nabla \hphi_i) + \nabla\left( \dfrac{1}{2}\nabla|\hphi_i|^2 + F(\hphi_i) \right)
		\]
		for $i = 1,2$. The above makes sense in $\bVsd$, since the chemical potential is not regular enough. Therefore, integrating by parts, we recover the identities
		\[
		\begin{split}
			\left\langle \hbu, \ia\left[ \hmu_1\nabla\hphi_1-  \hmu_2\nabla \hphi_2\right]\right\rangle_{\bVsd, \bVs} & = \left\langle \hmu_1\nabla\hphi_1-  \hmu_2\nabla \hphi_2, \ia\hbu\right\rangle_{\bVsd, \bVs} \\
			& = (\nabla \hphi_1 \otimes \nabla \hphi_1 - \nabla \hphi_2 \otimes \nabla \hphi_2, \nabla \ia \hbu)_\b{H} \\
			& = (\nabla \hphi_1 \otimes \nabla \hphi, \nabla \ia \hbu)_\bHs + ( \nabla \hphi \otimes \nabla \hphi_2, \nabla \ia \hbu)_\bHs.
		\end{split}
		\]
		On the other hand, by H\"{o}lder, Young and Ladyzhenskaya inequalities, we obtain
		\begin{align}
		\nonumber
		& \left| (\nabla \hphi_1 \otimes \nabla \hphi, \nabla \ia \hbu)_\b{H}
		+ ( \nabla \hphi \otimes \nabla \hphi_2, \nabla \ia \hbu)_\b{H} \right| \\
		\nonumber
		& \qquad \leq \left( \|\nabla \hphi_1\|_{\b{L}^4(\OO)}
		+ \|\nabla \hphi_2\|_{\b{L}^4(\OO)}\right)\|\nabla \hphi\|_\b{H}\|\nabla \ia \hbu\|_{\b{L}^4(\OO)} \\
		\nonumber
		& \qquad \leq \left( \|\hphi_1\|_{\b{L}^\infty(\OO)}^\frac{1}{2}
		\|\hphi_1\|_{V_2}^\frac{1}{2} + \|\hphi_2\|_{\b{L}^\infty(\OO)}^\frac{1}{2}\|\hphi_2
		\|_{V_2}^\frac{1}{2}\right)\|\nabla \hphi\|_\b{H}\|\hbu\|_\bHs^\frac{1}{2}
		\|\nabla \ia \hbu\|_{\bHs}^\frac{1}{2} \\
		& \qquad \leq\dfrac{1}{6}\|\hbu\|_\bHs^2 +
		\dfrac{1}{4}\|\nabla \hphi\|_\b{H}^2 +
		C\left( 1 + \|\hphi_1\|_{V_2}^2 + \|\hphi_2\|_{V_2}^2\right)\|\nabla \ia \hbu\|_{\bHs}^2.
		\label{eq:uniq12}
		\end{align}
		By Assumption \ref{hyp:G1} we also get (recall that $Y = \bVsd$),
		\begin{align}
			\nonumber
			\|\b A^{-1}G_{1}(\hbu_1) - \b A^{-1}G_1(\hbu_2)\|^2_{\cL^2(U_1, \bHs)}
			&=\|G_{1}(\hbu_1) - G_1(\hbu_2)\|^2_{\cL^2(U_1, \bVsd)} \\
			&\leq L_1^2\|\hbu\|_\bVsd \leq CL_1^2\|\nabla \ia \hbu\|^2_\bHs,
			\label{eq:uniq13}
		\end{align}
		since $\|\nabla \ia \bu\|_{\bHs}$ is an equivalent norm in $\bVsd$.
		Collecting \eqref{eq:uniq11}-\eqref{eq:uniq13}, we infer from \eqref{eq:uniq1} that
		\begin{align}
			\nonumber
			&\dfrac{1}{2}\|\nabla \ia \hbu(t\wedge\zeta_n)\|_\bHs^2
			+ \int_0^{t\wedge\zeta_n} \left[\dfrac{4}{6}\|\hbu(\tau)\|_\bHs^2
			- \dfrac{1}{4}\|\nabla \hphi(\tau)\|_H^2\right]\mathrm{d}\tau \\
			\nonumber
			&= \dfrac{1}{2}\|\nabla \ia\hbu_0\|^2_{\bHs}
			+\sup_{s\in[t\wedge\zeta_n]}\left|\int_0^s\left\langle \hbu(\tau),
			\ia\left[\big( G_{1}( \hbu_1(\tau)) - G_1(\hbu_2(\tau)) \big)
			\d  \widehat W_1(\tau)\right]\right\rangle_{\bVsd, \bVs}\right| \\
			&\qquad+ Cn^2\int_0^{t\wedge\zeta_n} \left( 1 + \|\hbu_1(\tau) \|^2_\bVs
			+ \|\hbu_2(\tau) \|^2_\bVs + \|\hphi_1(\tau)\|_{V_2}^2
			+ \|\hphi_2(\tau)\|_{V_2}^2 \right)\|\nabla \ia \hbu(\tau)\|_\bHs^2 \: \mathrm{d}\tau.
			\label{eq:uniqfinal1}
		\end{align}
		Before dealing with the stochastic integral in \eqref{eq:uniqfinal1},
		we consider $\Psi_2$. Applying the It\^{o} lemma to $\Psi_2(\hphi)$ yields,
		thanks to \cite[Theorem 4.2.5]{LiuRo},
		\begin{align}
			\nonumber
			&\dfrac{1}{2}\|\hphi(t\wedge\zeta_n)\|_{H}^2
			+ \int_0^{t\wedge\zeta_n} \left[ (\hphi(\tau), \hmu(\tau))_H
			+ \left(\hphi(\tau), \hbu_1(\tau) \cdot \nabla \hphi_1(\tau)
			- \hbu_2(\tau) \cdot \nabla \hphi_2(\tau) \right)_H \right] \d \tau\\
			\nonumber
			&= \dfrac{1}{2}\|\hphi_0\|_{H}^2
			+ \dfrac{1}{2}\int_0^{t\wedge\zeta_n} \| G_2(\hphi_1(\tau))
			-  G_2(\hphi_2(\tau))\|^2_{\cL^2(U_2,H)} \,\d\tau \\
			&\qquad+ \int_0^{t\wedge\zeta_n} (\hphi(\tau), \left[ G_2(\hphi_1(\tau))
			- G_2(\hphi_2(\tau))\right] \d \widehat W_2(\tau))_H .
			\label{eq:uniq2}
		\end{align}
		Observe now that, by the mean value theorem and \ref{hyp:potential},
		\begin{equation} \label{eq:uniq21}
			\begin{split}
				(\hphi, \hmu)_H & = \|\nabla \hphi\|^2_\b{H} + (F'(\hphi_1)-F'(\hphi_2), \hphi)_H \\
				& \geq \|\nabla \hphi\|^2_\b{H} - c_F\|\hphi\|_H^2.
			\end{split}
		\end{equation}
		Moreover, we have
		\begin{align}
		\nonumber
		\left|\left(\hphi, \hbu_1 \cdot \nabla \hphi_1  - \hbu_2 \cdot \nabla \hphi_2 \right)_H\right|
		& = \left|\left(\hphi, \hbu \cdot \nabla \hphi_1\right)_H
		+ \left(\hphi, \hbu_2 \cdot \nabla \hphi \right)_H \right| \\
		\nonumber
		& = \left|\left(\hphi, \hbu \cdot \nabla \hphi_1\right)_H\right| \\
		\nonumber
		& \leq \|\hbu\|_\bHs\|\hphi\|_{L^4(\OO)}\|\nabla \hphi_1\|_{\b{L}^4(\OO)} \\
		\nonumber
		& \leq \|\hbu\|_\bHs\|\hphi\|_H^\frac{1}{2}\|\hphi\|_{V_1}^\frac{1}{2}
		\|\hphi_1\|_{L^\infty(\OO)}^\frac{1}{2}\|\hphi_1\|_{V_2}^\frac{1}{2} \\
		& \leq \dfrac{1}{6}\|\hbu\|_\bHs^2 + \dfrac{1}{4}\|\nabla \hphi\|^2_\b{H}
		+ C\left( 1 + \|\hphi_1\|^2_{V_2} \right)\|\hphi\|^2_H.
		\label{eq:uniq22}
		\end{align}
		By \ref{hyp:G2}, we easily deduce
		\begin{equation} \label{eq:uniq23}
			 \| G_2(\hphi_1) -  G_2(\hphi_2)\|^2_{\cL^2(U_2,H)} \leq L_2^2\|\hphi\|^2.
		\end{equation}
		On account of \eqref{eq:uniq21}-\eqref{eq:uniq23}, from \eqref{eq:uniq2} we arrive at
		\begin{align}
			\nonumber
			&\dfrac{1}{2}\|\hphi(t\wedge\zeta_n)\|_{H}^2
			+ \int_0^{t\wedge\zeta_n} \left[ \dfrac{3}{4}\|\nabla \hphi(\tau)\|^2_\b{H}
			- \dfrac{1}{6}\|\hbu(\tau)\|^2_\bHs \right] \d \tau\\
			\nonumber
			&\leq \dfrac{1}{2}\|\hphi_0\|_{H}^2
			+ \sup_{s\in[0,t\wedge\zeta_n]}
			\left|\int_0^s (\hphi(\tau), \left[ G_2(\hphi_1(\tau))
			- G_2(\hphi_2(\tau))\right] \d \widehat W_2(\tau))_H\right|\\
			&\qquad+
			C\int_0^{t\wedge\zeta_n} \left( 1 + \|\hphi_1(\tau)\|^2_{V_2} \right)\| \hphi(\tau)\|^2_H \,\d\tau .
			\label{eq:uniqfinal2}
 		\end{align}
 		Adding \eqref{eq:uniqfinal1} and \eqref{eq:uniqfinal2} together, we obtain
		\begin{align}
			\nonumber
			&\dfrac{1}{2}\|\nabla \ia \hbu(t\wedge\zeta_n)\|_\bHs^2 +
			\dfrac{1}{2}\|\hphi(t\wedge\zeta_n)\|_{H}^2
			+\frac12 \int_0^{t\wedge\zeta_n} \left[\|\nabla \hphi(\tau)\|^2_\b{H}
			+\|\hbu(\tau)\|^2_\bHs \right] \d \tau\\
			\nonumber
			&\leq \dfrac{1}{2}\|\nabla \ia\hbu_0\|^2_{\bHs} +\dfrac{1}{2}\|\hphi_0\|_{H}^2
			+\sup_{s\in[t\wedge\zeta_n]}\left|\int_0^s\left\langle \hbu(\tau),
			\ia\big( G_{1}( \hbu_1(\tau)) - G_1(\hbu_2(\tau)) \big)
			\d  \widehat W_1(\tau)\right\rangle_{\bVsd, \bVs}\right|\\
			\nonumber
			&\qquad+ \sup_{s\in[0,t\wedge\zeta_n]}
			\left|\int_0^s (\hphi(\tau), \left[ G_2(\hphi_1(\tau))
			- G_2(\hphi_2(\tau))\right] \d \widehat W_2(\tau))_H\right|\\
			&\qquad+
			Cn^2\int_0^{t\wedge\zeta_n}
			\left[ 1 + \sum_{i=1,2}\left( \|\hbu_i(\tau) \|^2_\bVs
			+\|\hphi_i(\tau)\|_{V_2}^2 \right)\right]
			\left(\|\nabla \ia \hbu(\tau)\|_\bHs^2+\| \hphi(\tau)\|^2_H\right) \: \mathrm{d}\tau,
			\label{eq:uniq_luca}
 		\end{align}
		so that the Gronwall Lemma and the definition of $\zeta_n$ yield
		\begin{align}
			\nonumber
			&\|\nabla \ia \hbu(t\wedge\zeta_n)\|_\bHs^2 +
			\|\hphi(t\wedge\zeta_n)\|_{H}^2
			+\int_0^{t\wedge\zeta_n} \left[\|\nabla \hphi(\tau)\|^2_\b{H}
			+\|\hbu(\tau)\|^2_\bHs \right] \d \tau\\
			\nonumber
			&\leq e^{C(T+n^4)}\left(\|\nabla \ia\hbu_0\|^2_{\bHs} +\|\hphi_0\|_{H}^2 \right)\\
			\nonumber
			&\qquad+2e^{C(T+n^4)}
			\sup_{s\in[t\wedge\zeta_n]}\left|\int_0^s\left\langle \hbu(\tau),
			\ia\big( G_{1}( \hbu_1(\tau)) - G_1(\hbu_2(\tau)) \big)
			\d  \widehat W_1(\tau)\right\rangle_{\bVsd, \bVs}\right|\\
			\label{eq:uniq_luca2}
			&\qquad+2e^{C(T+n^4)} \sup_{s\in[0,t\wedge\zeta_n]}
			\left|\int_0^s (\hphi(\tau), \left[ G_2(\hphi_1(\tau))
			- G_2(\hphi_2(\tau))\right] \d \widehat W_2(\tau))_H\right|.
 		\end{align}
		Take now $\frac{p}{2}$-powers, the supremum (with respect to time) and expectations (with respect to $\hP$):
		let us deal with the stochastic integrals on the right hand side of \eqref{eq:uniq_luca2}.
		The Burkholder-Davis-Gundy inequality combined with the Young inequality and \ref{hyp:G1}
		entail, for every $\delta>0$, that
 		\begin{align}
		\nonumber
 		& \hE\sup_{s\in[t\wedge\zeta_n]}\left|\int_0^s\left\langle \hbu(\tau),
			\ia\big( G_{1}( \hbu_1(\tau)) - G_1(\hbu_2(\tau)) \big)
			\d  \widehat W_1(\tau)\right\rangle_{\bVsd, \bVs}\right|^\frac{p}{2}\\
		\nonumber
 		& \qquad \leq C\hE\left( \int_0^{t\wedge\zeta_n}
		\|\nabla\ia\hbu(s)\|^2_{\bHs}
		\| G_{1}( \hbu_1(s)) - G_1(\hbu_2(s))\|^2_{\cL^2(U_1,\bVsd)}\,\d s\right)^\frac{p}{4} \\
		\nonumber
 		& \qquad \leq
		C\hE\left( \sup_{s\in[0,t\wedge\zeta_n]}
		\|\nabla\ia\hbu(s)\|^2_{\bHs}\int_0^{t\wedge\zeta_n}
		\|\hbu(\tau)\|^2_{\bVsd}\,\d \tau\right)^\frac{p}{4} \\
 		& \qquad \leq
		\delta\hE\sup_{s\in[0,t\wedge\zeta_n]}
		\|\nabla\ia\hbu(s)\|^p_{\bHs} +
		C_\delta
		\widehat\E\int_0^{t\wedge\zeta_n}
		\|\nabla\ia\hbu(\tau)\|^p_\bHs\,\d \tau,
		\label{eq:uniq31}
 		\end{align}
 		while the same inequalities and \ref{hyp:G2} also yield
 		\begin{align}
		\nonumber
 		& \hE\sup_{s\in[t\wedge\zeta_n]}\left| \int_0^\tau (\hphi(s),
		\left[ G_2(\hphi_1(s)) - G_2(\hphi_2(s))\right] \d \widehat W_2(s))_H \right|^\frac{p}{2} \\
		\nonumber
 		& \qquad \leq C\hE\left(\int_0^{t\wedge\zeta_n}
		\|\hphi(s)\|^2_H\| G_2(\hphi_1(s)) - G_2(\hphi_2(s))\|^2_{\cL^2(U_2,H)}\,\d s \right)^\frac{p}{4} \\
 		& \qquad \leq \delta\hE\sup_{s\in[t\wedge\zeta_n]}\|\hphi(s)\|^p_H
		+ C_\delta\E\int_0^{t\wedge\zeta_n} \|\hphi(\tau)\|^p_H\,\d \tau.
		\label{eq:uniq32}
 		\end{align}
 		Taking \eqref{eq:uniq31} and \eqref{eq:uniq32} into account in \eqref{eq:uniq_luca2}
		and choosing $\delta$ small enough, an application of the Gronwall lemma entails the claimed
		continuous dependence estimate. In turn,
		upon choosing $\hbu_{0,1} = \hbu_{0,2}$ and $\hphi_{0,1} = \hphi_{0,2}$,
		this also yields
		$\hbu_1=\hbu_2$ and $\hphi_1=\hphi_2$ on the stochastic interval
		$[\![0,\zeta_n]\!]$ for every $n\in\enne$.
		Hence pathwise uniqueness of the solution follows since $\zeta_n\nearrow T$ almost surely.
	\end{proof} \noindent
	The existence of a probabilistically-strong solution follows from standard results (see, for instance, \cite[Theorem 2.1]{Rockner08}), which also turns out to be unique. The existence and uniqueness (up to a constant) of a pressure $\pi \in L^\andrea{\frac p2}(\Omega; W^{-1,\infty}(0,T;H))$ can be deduced arguing as in Subsection \ref{ssec:pressure}. The proof of Theorem \ref{th:probstrongsol} is finished.

\bigskip\noindent
{\bf Acknowledgments.} The second and third authors are members of Gruppo Nazionale per l'Analisi Matematica, la Probabilit\`{a} e le loro Applicazioni
(GNAMPA), Istituto Nazionale di Alta Matematica \linebreak (INdAM). The second author has been partially funded by MIUR-PRIN
Grant 2020F3NCPX \linebreak ``Mathematics  for industry 4.0 (Math4I4)''.

	\appto{\bibsetup}{\raggedright}
	\printbibliography[title=Bibliography]
\end{document}